\documentclass[11pt,a4paper]{article}
\usepackage{amsmath}
\usepackage{amssymb}
\usepackage{amsthm}
\usepackage[latin1]{inputenc}
\usepackage[T1]{fontenc}
\usepackage{tabularx}
\usepackage{oldgerm}
\usepackage[all,2cell,ps]{xy}
\UseAllTwocells
\UsePSspecials{}
\CompileMatrices
\newtheorem{thm}{Theorem}[subsection]
\newtheorem{prop}[thm]{Proposition}
\newtheorem{lemma}[thm]{Lemma}
\newtheorem{cor}[thm]{Corollary}

\newtheorem{defi}[thm]{Definition}
\newtheorem{remark}[thm]{Remark}
\def\df[#1]{%
  \ar@{}[#1]|(.35)*{}="A" \ar@{}[#1]|(.65)*{}="B" 
  \ar@{=>}"A";"B" }
\def\tdf[#1]{%
  \ar@{}[#1]|(.25)*{}="A" \ar@{}[#1]|(.55)*{}="B" 
  \ar@{=>}"A";"B" }

\newcommand{\on}{\operatorname}
\newcommand{\inlim}[2]{\underset{#1}{\varinjlim}\,{#2}}  

\newcommand{\prlim}[2]{\underset{#1}{\varprojlim}\,{#2}}

\newcommand{\twodilim}[2]{\underset{#1}{\on{2\!}\varinjlim}\,{#2}}

\newcommand{\ind}{\operatorname{Ind}}

\newcommand{\Ob}{\operatorname{Ob}}
\newcommand{\Mor}[1]{\operatorname{Mor}{\mathcal{#1}}}

\newcommand{\Hom}[1]{\operatorname{Hom}_{_{#1}}}

\newcommand{\Perv}{\operatorname{Perv}}

\newcommand{\coker}{\operatorname{coker}}
\newcommand{\coim}{\operatorname{coim}}
\newcommand{\im}{\operatorname{im}}

\newcommand{\ol}{\overline}
\newcommand{\ul}{\underline}

\newcommand{\sct}{\scriptscriptstyle}

\newcommand{\geqs}{\geqslant}
\newcommand{\leqs}{\leqslant}
\newcommand{\mc}{\mathcal}
\newcommand{\mr}{\mathrm}

\newcommand{\me}{\matheus}
\newcommand{\ra}{\rightarrow}
\newcommand{\lra}{\longrightarrow}
\newcommand{\N}{\mathbb{N}}

\newcommand{\Z}{\mathbb{Z}}
\newcommand{\R}{\mathbb{R}}
\newcommand{\C}{\mathbb{C}}

\newcommand{\Db}{\operatorname{D}^{\operatorname{b}}}
\newcommand{\I}{\operatorname{I}}
\newcommand{\Dr}{\operatorname{R}\!}
\renewcommand{\SS}{\operatorname{SS}}
\newcommand{\supp}{\operatorname{supp}}
\newcommand{\Cc}{\C\text{-c}}
\newcommand{\Rc}{\R\text{-c}}
\newcommand{\murh}{\mu\!\on{RH}}
\newcommand{\musol}{\mu\!\on{Sol}}
\newcommand{\thom}{\on{T}\!\me{H}om}
\newcommand{\ctimes}[1]{\underset{#1}{\times}}

\DeclareMathAlphabet{\matheus}{U}{eus}{m}{n}
\newdir{ (}{{}*!/-5pt/@^{(}}
\newdir{ )}{{}*!/-5pt/@_{(}}
\numberwithin{equation}{subsection}
\title{\Huge{\textbf{Microlocal perverse sheaves}}}
\author{Ingo Waschkies}
\date{}
\begin{document}
\maketitle

\begin{abstract}
In this paper we give an explicit construction of the stack of microlocal perverse 
sheaves on the projective cotangent bundle of a complex manifold. Microlocal
perverse sheaves will be represented as complexes of analytic ind-sheaves which
have recently been studied by Kashiwara-Schapira \cite{KS2}. This description allows
us to formulate the microlocal Riemann-Hilbert correspondance in order to establish
an equivalence of stacks with the stack of regular holonomic microdifferential modules. 
\end{abstract}

\tableofcontents

\newpage

\section{Introduction}

In \cite{K1} Kashiwara constructed the stack of microdifferential 
modules on a complex contact manifold, generalizing the stack of modules over
the ring of microdifferential operators $\me{E}_X$ on the projective cotangent
bundle $P^*X$ of a complex variety $X$. In his discussion Kashiwara asked for
the construction of a stack of microlocal (complex) perverse sheaves that 
should be equivalent to the stack of regular holonomic modules by a microlocal 
Riemann-Hilbert correspondance. Such a stack should be defined over any field $k$,
but the Riemann-Hilbert morphism only makes sense over $\C$.

There have been several attempts to construct a local version of such a stack.
In \cite{An1} Andronikof defined a prestack on $P^*X$ and announced the 
microlocal Riemann-Hilbert correspondance on the stalks. However, at this time 
there did not exist tools to define a global microlocal Riemann-Hilbert morphism. 
Another topological construction was proposed in \cite{GMV}, but to our knowledge 
this project has neither been completed nor published.

Our approach makes use of the theory of
analytic ind-sheaves, recently introduced  in 
\cite{KS2} by Kashiwara and Schapira. Hence, microlocal perverse sheaves on a 
$\C^\times$-conic open subset $U\subset T^*X$ will be ind-sheaves (or more precisely 
objects of the derived category of 
ind-sheaves) on $U$ contrary to the construction of
\cite{An1} in which microlocal perverse sheaves on $U\subset T^*X$ were represented 
by complexes of sheaves on the base space $X$. The theory of ind-sheaves provides us 
with a nice represantitive of the stack associated to the prestack of \cite{An1} and 
allows us to use the machinery developped in \cite{KS2}. The essential tool 
in this description is Kashiwara's functor of ind-microlocalization 
$\mu:\,\Db(k_X)\ra \Db(\I(k_{T^*X}))$ of \cite{K4}. This functor enables us to 
define explicitly a global Riemann-Hilbert morphism when $k=\C$.

In the future, we will hopefully show that we can actually patch (a twisted 
version of) this stack on a complex contact manifold and prove the Riemann-Hilbert 
theorem in the complex case. 

In more detail, the contents of this paper are as follows.

\textbf{Section 2} gives a criterion for a subprestack of the prestack of derived 
categories of ind-sheaves on a manifold to be a stack. It is a generalization of a 
proof of \cite{KS3} showing that the prestack of perverse sheaves is a substack of the 
derived category of sheaves with $\C$-constructible cohomology. Then we investigate 
abelian stacks on a topological spaces. Roughly speaking, an additive stack on a 
topological space is abelian if and only if its stalks are abelian categories and we 
have a ``lifting property'' for kernels and cokernels. This will be applied in 
Section 7 to construct the abelian stack of microlocal perverse sheaves as a substack 
of the prestack of derived ind-sheaves on the contangent bundle of a complex manifold.

\textbf{Section 3} first recalls the theory of microlocalization of \cite{KS3} on 
a real manifold $X$. We do not review in detail the theory of the micro-support of 
sheaves but concentrate on the definition of the microlocal category $\Db(k_X,S)$ where 
$S\subset T^*X$ is an arbitrary subset. It is defined as the localization of the 
category $\Db(k_X)$ by the objects $\me{F}\in\Db(k_X)$ whose micro-support does
not intersect $S$. For any $\me{F},\me{G}\in\Db(k_X)$ we get a natural morphism
$$ \on{Hom}_{\Db(k_X,S)}(\me{F},\me{G}) \lra \on{H}^0(S,\mu hom(\me{F},\me{G})). $$
In the case where $S=\{pt\}$ the category $\Db(k_X,S)$ has been
intensively studied in \cite{KS3}, and in particular it is proved that the morphism
above is an isomorphism.\\
We will show that this result is still valid in the category $\Db(k_X,\{x\}\times 
\dot\delta)$ where $x$ is a point of $X$, $\delta\subset T^*_xX$ a closed cone and 
$\dot \delta=\delta\setminus \{0\}$. Later we will be mainly interested in 
the case where $\delta$ is a complex line. The main tool is the refined microlocal
cut-off lemma for non-convex sets, which we recall adding a few comments.
We will also need the cut-off functor in Section 6.

\textbf{Section 4} extends the definitions and results of Section 3 first to 
$\R$-constructible then to $\C$-constructible sheaves.
There are two natural ways to define the microlocalization of the derived category of 
$\R$-constructible sheaves. We either localize the category 
$\Db_{\Rc}(k_X)$ by sheaves whose micro-support does not intersect $S$ or we 
take the full subcategory of $\Db(k_X,S)$ whose objects are represented by 
$\R$-constructible sheaves. Following \cite{An2} we will use the first definition. 
One important question is whether or not the two definitions coincide. The main 
result of this section is that this is the case when $S=\{x\}\times \dot\delta$.\\
Now suppose that $X$ is a complex manifold. Recall that $\C$-constructible sheaves 
may be defined by a microlocal property: an object $\me{F}\in\Db_{\Rc}(k_X)$ is 
$\C$-constructible if and only if its micro-support is $\C^\times$-conic. It is 
then natural to define a microlocally complex constructible sheaf as a microlocal
$\R$-constructible sheaf whose micro-support is $\C^\times$-conic on $S$.

In \textbf{Section 5} we show that the constructions of Section 4 are locally 
``invariant under quantized contact transformations''. 

\textbf{Section 6} is devoted to the study of microlocally $\C$-constructible 
sheaves in the category $\Db(k_X,\C^\times p)$. In Section 5 we have shown that 
the category $\Db_{\Cc}(k_X,\C^\times p)$ is invariant by quantized contact 
transformation. Hence we are reduced to study microlocally $\C$-constructible 
sheaves in generic position, i.e. complexes of sheaves whose micro-support is 
contained in $T^*_ZX$ for a complex (not necessarily smooth) hypersurface $Z$ 
in a neighborhood of $p$. We give a complete proof that microlocally 
$\C$-constructible sheaves in generic position may be represented by 
$\C$-constructible sheaves (as anounced in \cite{An2}).

Following \cite{An1}, we define in \textbf{Section 7} the category of microlocal 
perverse sheaves as a full subcategory of $\Db_{\Cc}(k_X,\C^\times p)$. An object 
$\me{F}\in\Db_{\Cc}(k_X,\C^\times p)$ is perverse if for any non-singular 
point $q\in\SS(\me{F})$ in a neighborhood of $\C^\times p$ the complex $\me{F}$  
is isomorphic in $\Db(k_X,\C^\times q)$ to a constant sheaf $M_Y[d_Y]$ supported 
on a closed submanifold $Y\subset X$. This definition is natural in view of the 
microlocal characterization of perverse sheaves of \cite{KS3} and also leads to
definition of a prestack of microlocal perverse sheaves on $P^*X$. Then we prove 
that the category $\Db_{\on{perv}}(k_X,\C^\times p)$ is abelian as has been announced 
in \cite{An1}. Our proof gives a refined result which allows us to conclude that 
the stack associated to this prestack is abelian. This stack is the stack of 
microlocal perverse sheaves on $P^*X$. However we will need a more explicit 
representation in order to define the microlocal Riemann-Hilbert correspondance.

In \textbf{Section 8} we finally define microlocal perverse sheaves as objects 
of the derived category of ind-sheaves on conic open subsets of $T^*X$. 
In Section 7 we have constructed the category of microlocal perverse sheaves at any
$p\in P^*X$ (or on $\C^\times p\subset \dot T^*X$) which will be equivalent to the 
stalk of the stack $\mu Perv$ of microlocal perverse sheaves. 
The idea of the construction of $\mu Perv$ is to use the fact Kashiwara's functor 
$\mu$ of ind-microlocalization induces a fully faithful functor from 
$\Db_{\on{perv}}(k_X,\C^\times p)$ into the stalk of the prestack of bounded derived
categories of ind-sheaves on $\C^\times$-conic subsets of $T^*X$. Then we can define 
a microlocal perverse sheaf on a conic open subset $U\subset T^*X$ as an object of 
$\Db(\I(k_{U}))$ that is isomorphic to a microlocal perverse sheaf of 
$\Db_{\on{perv}}(k_X,\C^\times p)$ at any point of $p\in U$. Finally, we show that 
the stack of microlocal perverse sheaves is canonically equivalent to the stack 
associated to the prestack of the last section. 

In \textbf{Section 9} we establish the microlocal Riemann-Hilbert correspondance:
$$ \xymatrix@C=3cm{ {\mu Perv(\Omega)} \ar@<+3pt>[r]^{\murh} & 
       {\mc{H}ol\mc{R}eg(\me{E}_X|_{\gamma^{-1}\Omega}).} 
        \ar@<+3pt>[l]^{\musol} } $$
It is defined by the formulas
$$ \musol(\me{M})= 
  \Dr\me{IH}om_{\beta(\me{E}_X|_{\gamma^{-1}\Omega})}(\beta(\me{M}),
    \mu \me{O}_X|_{\gamma^{-1}\Omega})$$
$$ \murh(\me{F})=\gamma^{-1}_{\Omega}\Dr \gamma_{\Omega *}
  (\Dr\me{H}om(\me{F},\mu \me{O}^t|_{\gamma^{-1}(\Omega)})) $$
where $\gamma_{\Omega}$ is the restriction of the natural map 
$\gamma:\dot T^*X\ra P^*X$ to $\gamma^{-1}\Omega$, $\Dr \me{IH}om,\beta$ are
well-known functors from the theory of ind-sheaves (cf. \cite{KS2}) and 
the ``ring'' $\me{O}^t\in\Db(\I(k_X))$ is the complex of tempered holomorphic 
functions on $X$ of loc. cit.

\textbf{Appendix A} gives a short introduction to stacks with emphasis on
the special properties resulting from the fact that we work on a
topological space. 

In \textbf{Appendix B} we give a short summary of the properties of Kashiwara's functor 
$\mu$.

First of all, I thank my thesis director P. Schapira both for having suggested this 
subject to me, and for always having been ready with precious help, guidance and 
encouragement throughout the last three years. Secondly, my gratitude goes out to 
M. Kashiwara with whom I had many invaluable conversations. I would particularly 
like to thank him for having shared with me his unpublished work on the microlocalization 
of ind-sheaves which provided me with certain key ideas on which this paper is based.
It goes without saying, of course, that I could never have been able to complete this 
work without either of them.
Finally, I would like to thank A. D'Agnolo, P. Polesello, F. Ivorra and D.-C. Cisinski 
for many useful discussions.

\section{Abelian substacks of a prestack}

Perverse sheaves on a complex manifold $X$ are local objects - they form an 
abelian stack which is a subprestack of the prestack of (derived) sheaves on $X$
(see \cite{BBD} for the general theory of perverse sheaves, see also 
\cite{KS3}, Chapter X, for a microlocal approach to perverse sheaves). In Section
2.1 we will generalize the method used in \cite{KS3} in order to prove that this 
subprestack is actually a stack. In particular we will show that a similar method
can be applied to find substacks of the prestack of (derived) ind-sheaves.
The abelian structure of the stack of perverse sheaves is defined by a $t$-structure
on the triangulated prestack of derived categories of sheaves with $\C$-constructible
cohomology. However, the category of microlocal perverse sheaves will not be defined
as the heart of a $t$-structure. Our strategy is based on the idea that a stack is
``almost'' abelian, if its stalks are abelian categories. Roughly speaking, an additive 
stack is abelian if and only if its stalks are abelian and kernels and cokernels can be
lifted to small neighborhoods. We will investigate this statement more precisely 
in Sections 2.2 and 2.3. 

\subsection{A criterion for substacks}

The basic definitions from the theory of stacks (on a topological space) are recalled
in Appendix B. The results on proper stacks and ind-sheaves that we will use can be 
found in \cite{KS2}.
\begin{defi}\label{localprop}
  Consider a prestack $\me{C}$ on a topological space $X$. We say that
  a full subprestack $\me{C'}\subset\me{C}$ is \textbf{defined by a local
  property} (with respect to $\me{C}$) if the following conditions
  are satisfied:
  \begin{itemize}
    \item[(i)] the prestack $\me{C'}$ is stable by isomorphisms,
      i.e. if $U\subset X$ is open, $A\in\Ob{\me{C'}(U)}$ then
      any object $B\in\Ob{\me{C}(U)}$ isomorphic to $A$ is also an
      object of $\Ob{\me{C'}(U)}$, 
    \item[(ii)] if $U\subset X$ is open and $A\in\Ob{\me{C}(U)}$ then
      $A\in\Ob{\me{C'}(U)}$ if and only if there is an open
      covering $U=\bigcup_{i\in I}{U_i}$ such that $A|_{_{U_i}}\in
      \Ob{\me{C'}(U_i)}$ for all $i\in I$.
   \end{itemize}
\end{defi}
\begin{remark}
  \em Consider a full subprestack $\me{C}'\subset \me{C}$ and a point $p\in X$. Then 
  the natural functor $\me{C}'_p\ra \me{C}_p$ is fully faithful. 
  Therefore the (full) subprestack $\me{C}'$ is defined by a local property if and only 
  if for any object $A\in\Ob{\me{C}(U)}$ the statements (a) and (b) below
  are equivalent:
  \begin{itemize}
             \item[(a)] $A\in\Ob{\me{C}'(U)}.$ 
             \item[(b)] For every $p\in X$ the object $A$ is in the essential image
                 of the functor $\me{C}'_p\ra \me{C}_p$, i.e. there exists an object
                 $B\in\Ob\me{C}'_p$ such that $A$ is isomorphic to $B$ in $\me{C}_p$.
  \end{itemize} \em 
\end{remark}
\begin{lemma}\label{lemma41}
  Let $\me{C}$ be a triangulated prestack. Assume moreover that 
  \begin{itemize}
    \item[(1)] for any $V\subset U$ the restriction functor $i_{VU}^{-1}$ 
       has a fully faithful left adjoint $i_{VU!}$\begin{footnote}{Recall that 
       $i_{VU!}$ is fully faithful if and only if the adjunction morphism $\on{Id}\ra 
       i_{VU!}i_{VU}^{-1}$ is an isomorphism. Also note that for any three open 
       subsets $W\subset V\subset U$ the isomorphism 
       $i_{WV}^{-1}i_{VU}^{-1}\simeq i_{WU}^{-1}$ induces an isomorphism
       $i_{WV!}i_{VU!}\simeq i_{WU!}$.}\end{footnote},
    \item[(2)] these functors satisfy the base change theorem, i.e. for any
        Cartesian square of open subsets
        $$\xymatrix{ U_{12}  \ar[r]\ar[d]  & U_1 \ar[d]  \\
                     U_2   \ar[r] \ar@{}[ur]|\square  &  V } $$
        we have $i_{U_{12}U_2!}i_{U_{12}U_1}^{-1}\simeq 
        i_{U_2V}^{-1}i_{U_1V!}$, where
        $U_{12}=U_1\cap U_2$.
  \end{itemize}
  Consider the union of two open sets $U=U_1\cup U_2$ and suppose that we 
  are given
  \begin{itemize}
     \item[(i)] objects $A_1\in\Ob{\me{C}(U_1)}$ and
          $A_2\in\Ob{\me{C}(U_2)}$,
     \item[(ii)] an isomorphism $f_{21}:\,A_1|_{_{U_{12}}}
        \overset{\sim}{\lra} A_2|_{_{U_{12}}}$ in $\me{C}(U_{12})$.
  \end{itemize}
  Then there exist an object $A\in\Ob{\me{C}(U)}$ and isomorphisms 
  $f_{1}:\,A|_{_{U_1}}\overset{\sim}{\ra} A_1$, $f_2:\,A|_{_{U_2}}
  \overset{\sim}{\ra}A_2$ that are compatible with $f_{21}$ on $U_{12}$, 
  i.e. the following diagram commutes:
  $$ \xymatrix{
       A_1|_{_{U_{12}}} \ar[rr]^{f_{21}}_{\sim} & & A_2|_{_{U_{12}}} \\
        & A|_{_{U_{12}}} \ar[ur]_{f_2|_{_{U_{12}}}}^{\sim} 
        \ar[ul]^{f_1|_{_{U_{12}}}}_{\sim} &  }$$
\end{lemma}
\begin{proof}
Let us simplify the notations for the restriction functors by suppressing $U$.
Hence we have the functors
\begin{align*}
   i_1^{-1}:\,\me{C}(U)\ra &\me{C}(U_1) \qquad\qquad  
   i_2^{-1}:\,\me{C}(U)\ra \me{C}(U_2) \qquad\qquad 
   i_{12}^{-1}:\,\me{C}(U)\ra\me{C}(U_{12})\\
   i_{12,1}^{-1}&:\,\me{C}(U_1)\ra \me{C}(U_{12}) \qquad\qquad
   i_{12,2}^{-1}:\,\me{C}(U_2)\ra\me{C}(U_{12}) 
\end{align*}
We use a similar notation with lower-case symbol $!$ for the left adjoints of 
these functors. Now define $A$ by chosing a distinguished triangle
$$ \xymatrix@C=1.2cm{
    i_{12!}(A_1|_{_{U_{12}}}) \ar[r]^(.45){\binom {g_1}{g_2}} &
    i_{1!}A_1\oplus i_{2!}A_2 \ar[r]^(.65){(h_1,h_2)} &  A \ar[r]^{+} & }$$
where the first morphism defining the triangle is given by 
$$ \xymatrix@R=0.3cm{ g_1:\ i_{12!}(A_1|_{_{U_{12}}})\ar[r]^{\sim} &
   i_{1!}i_{12,1!}(A_1|_{_{U_{12}}}) \ar[r] &
   i_{1!}(A_1), \\
   g_2:\ i_{12!}(A_1|_{_{U_{12}}}) \ar[r]^(.5){-f_{21}}_(.5){\sim} &
   i_{12!}(A_2|_{_{U_{12}}}) \ar[r]^(.45){\sim} &
   i_{2!}i_{12,2!}(A_2|_{_{U_{12}}}) \ar[r] &
   i_{2!}(A_2).} $$
Define $f_1$ and $f_2$ by
$$ \xymatrix@R=0cm@C=1.3cm{
  f_1:\ A_1\simeq i_1^{-1}i_{1!} A_1  \ar[r]^(.65){h_1|_{_{U_1}}} & A|_{_{U_1}}, \\
  f_2:\ A_2\simeq i_2^{-1}i_{2!} A_2 \ar[r]^(.65){h_2|_{_{U_2}}} & A|_{_{U_2}}. 
  } $$
Let us check that $f_1$ is an isomorphism. It is sufficent to show that $g_2|_{_{U_1}}$
is an isomorphism. Since by the base change theorem we have
$$ i_1^{-1}i_{2!}(A_2)\simeq i_{12,1!}(A_2|_{_{U_{12}}}) \simeq
   i_1^{-1}i_{12!}(A_2|_{_{U_{12}}}),$$ 
we get the result. A similar argument shows that $f_2$ is an isomorphism. To prove 
that these isomorphisms are compatible with $f_{21}$ one uses the fact that 
$h_1g_1+h_2g_2=0$.
\end{proof}
\begin{remark}
  \emph{In the situation of the preceding Lemma \ref{lemma41}, suppose that
  we are given a full (but not necessarily triangulated) subprestack 
  $\me{C}'\subset\me{C}$ that is defined by a local property. Then the 
  lemma holds in $\me{C}'$, i.e. if the objects $A_1$,$A_2$ are in $\me{C}'$ then 
  the object $A$ lies also in $\me{C}'$. Indeed, we may patch the given objects 
  $A_1, A_2$ of $\me{C}'$ to an object $A$ in the prestack $\me{C}$ using Lemma 
  \ref{lemma41}. Then the axioms (cf. Definition \ref{localprop}) immediately imply 
  that $A$ is an object of $\me{C}'$.\\
  Note that if moreover $\me{C}'$ is separated, then the object $A$ is unique up
  to unique isomorphism\begin{footnote}{More precisely, if $A'$ is another object with
  isomorphisms $f_i':\, A'|_{_{U_i}}\overset{\sim}{\ra}A_i$ for $i\in\{1,2\}$ such
  that $f_{21}\circ f_1'|_{_{U_{12}}}=f_2'|_{_{U_{12}}}$ then there exists a unique
  isomorphism $\varphi:\,A\overset{\sim}{\ra} B$ such that 
  $f_i'\circ\varphi|_{_{U_i}}=f_i$ for $i\in\{1,2\}$.}\end{footnote}.}
\end{remark}
Let us apply Lemma \ref{lemma41} to the prestack of bounded derived categories of
ind-sheaves. Denote by $\Db(\I(k_*))$ the prestack $U\mapsto\Db(\I(k_U))$ on 
a locally compact space $X$ with a countable base of open sets. This prestack has 
the following properties:
\begin{itemize}
  \item[(i)] it is a triangulated prestack,
  \item[(ii)] if $V\subset U\subset X$ are open subsets and if we denote as usual by 
     $i_{{VU}}^{-1}:\,\on{D^b}(I(k_{_U}))\ra \on{D^b}(I(k_{_V}))$ the restriction 
     functor, then this functor has a fully faithful left adjoint $\Dr i_{{VU!!}}$,
  \item[(iii)] the base change theorem is satisfied (cf. Lemma \ref{lemma41}). 
\end{itemize}
The proof is based on the fact that $\I(k_*)$ is a proper stack (cf. \cite{KS2}). 
Hence we get
\begin{cor}
  Consider the prestack $\Db(I(k_*))$ on a locally compact topological space $X$ 
  and two open subsets $U_1,U_2\subset X$. Set $U=U_1\cup U_2$.\\
  Suppose that we are given the following data
  \begin{itemize}
    \item[(1)] two objects $\me{F}_1\in\Ob{\on{D^b}(I(k_{_{U_1}}))}$, 
         $\me{F}_2\in\Ob{\on{D^b}(I(k_{_{U_2}}))}$,
    \item[(2)] an isomorphism $f_{21}:\,\me{F}_1|_{_{U_{12}}}
         \overset{\sim}{\lra} \me{F}_2|_{_{U_{12}}}$.
  \end{itemize}
  Then there exist an object $\me{F}\in\Ob{\on{D^b}(I(k_{_{U}}))}$ and 
  isomorphisms $f_{1}:\,\me{F}|_{_{U_1}}\overset{\sim}{\ra} \me{F}_1$, 
  $f_2:\,\me{F}|_{_{U_2}}\overset{\sim}{\ra}\me{F}_2$ compatible with 
  $f_{21}$ on $U_{12}$.
\end{cor}
Now let us state Proposition 10.2.9. of \cite{KS3} in a slightly more general
context and change the proof so that we may adapt it later to the case of ind-sheaves.
\begin{prop}\label{propo}
  Let $X$ be a locally compact paracompact space with a countable base of open sets. 
  Consider a proper stack\begin{footnote}{For the definition of a proper stack 
  see \cite{KS2}. A proper stack $\me{A}$ and the associated prestack of bounded 
  derived categories $\Db(\me{A})$ satisfy the hypothesis of Lemma \ref{lemma41}. 
  Moreover for each open subset $U\subset X$, the abelian category $\me{A}(U)$ 
  admits filtered exact colimits and the restriction functors commute to such 
  colimits.}\end{footnote} 
  $\me{A}$ such that for every open subset $U\subset X$ the category $\me{A}(U)$ has
  enough injective objects. Denote by $\Db(\me{A})$ the associated prestack of bounded 
  derived categories.\\ 
  Let $\me{C}\subset\Db(\me{A})$ be a separated full subprestack that is 
  defined by a local property.\\
  Then $\me{C}$ is a stack. 
\end{prop}
\begin{proof}
We have to show that $\me{C}$ satisfies the patching condition (cf. Definition 
\ref{patchcond} in Appendix B). Since $X$ is 
paracompact we have to verify the patching condition only for countable coverings.\\
Since the conditions of Lemma \ref{lemma41} are satisfied for the prestack 
$\Db(\me{A})$ and $\me{C}$ is defined by a local property, we know that $\me{C}$ 
satisfies the patching condition for any covering of type $U=U_1\cup U_2$.\\
Using the fact that $\me{C}$ is separated, we can easily verify by induction that the 
patching condition is satisfied for finite coverings.\\
Therefore, using again the fact that $\me{C}$ is separated, it is sufficent to prove 
that objects may be patched in $\me{C}$ for open coverings of type 
$U=\bigcup_{_{n\in\N}}{U_n}$ where $U_{n}\subset U_{n+1}$.\\
Consider a family of objects $A_n\in\Ob{\on{D^b}(\me{A}(U_n))}$
and isomorphisms $f_{n-1}:\,A_{n-1}\overset{\sim}{\ra}A_n|_{_{U_{n-1}}}$ (the other
isomorphisms are uniquely determined by the cocycle condition).\\ 
Denote by $i_n:\,U_n\hookrightarrow U$ 
the inclusion map. Then we lift the morphisms of th system to 
$$ \xymatrix@C=2cm{
    g_{n-1}:\ i_{n-1!}A_{n-1} \ar[r]^{i_{n-1!}(f_{n-1})}_{\sim} &  
       i_{n-1!}A_{n}|_{_{U_{n-1}}} \ar[r] &  i_{n!}A_n } $$ 
in $\Db(\me{A}(U))$. Hence we get a family of morphisms
$\{g_{n-1}:\,i_{n-1!}A_{n-1}\ra i_{n!}A_n\}_{n\geqs 1}$ in $\Db(\me{A}(U))$.
Note that $g_{n-1}|_{U_{n-1}}$ is an isomorphism by the base change theorem.\\
By hypothesis $\me{A}(U)$ has enough injective objects and therefore we have an 
equivalence of categories $\on{K}^b(\on{Inj}(\me{A}(U)))\simeq \Db(\me{A}(U))$.
Here $\on{Inj}(\me{A}(U))$ denotes the full abelian subcategory of $\me{A}(U)$ whose 
objects are the injective objects of $\me{A}(U)$ and  $\on{K}^b(\on{Inj}(\me{A}(U)))$ 
is the triangulated category of bounded complexes of $\on{Inj}(\me{A}(U))$ where 
morphisms of complexes are considered up to homotopy. 
Hence there exist objects $I_n$ and morphisms 
$h_{n-1}:\,I_{n-1}\ra I_{n}$ in $\on{C}^b(\on{Inj}(\me{A}(U)))$ such that the diagram
$\{h_{n-1}:\,I_{n-1}\ra I_{n}\}_{n\geqs 1}$ is isomorphic (in $\Db(\me{A}(U))$) to 
the diagram $\{g_{n-1}:\,i_{n-1!}A_{n-1}\ra i_{n!}A_n\}_{n\geqs 1}$.\\ 
Let $A=\varinjlim I_n$ and consider $A$ as an object in $\Db(\me{A}(U))$. 
Then $A|_{_{U_n}}$ is quasi-isomorphic in $\on{C}^b(\me{A}(U))$ to $I_n$ 
because for $m\geqs n$ the morphism $I_n\ra I_m|_{U_{n}}$ is a 
quasi-isomorphism. Hence there are natural isomorphisms $A|_{_{U_n}}\simeq A_n$ in 
$\Db(\me{A}(U))$, and a simple diagram chase shows that they are compatible with the 
morphisms $f_n$.
\end{proof}
\begin{remark}\label{trueitis}
  \em Note that in Proposition \ref{propo} the hypothesis that the categories 
  $\Db(\me{A}(U))$ possess enough injective objects can be weakened. During the proof, 
  we actually only use the fact that any diagram in $\Db(\me{A}(U))$ of type 
  $\{A_n\ra A_{n+1}\}_{n\geqs 0}$ can be lifted to a diagram 
  $\{I_n\ra I_{n+1}\}_{n\geqs 0}$ in $\on{C}^b(\me{A}(U))$. We do not use the 
  fact that the objects $I_n$ are injective.
\end{remark}
We need a proposition of \cite{KS2} (part of Theorem 11.2.6).
\begin{prop}\label{injectiveobjects}
  Let $\mc{A}$ be an abelian category with a system of strict
  generators. Denote by $\ind{\mc{A}}$ the category of ind-objects of $\me{A}$ 
  and let $S\subset \Ob{\on{Ind}(\mc{A})}$ be a small subset. Then there exists 
  an essentially small\begin{footnote}{A category is essentially small if it is 
  equivalent to a small category.}\end{footnote} full abelian subcategory 
  $\mc{B}\subset\mc{A}$ such that
  \begin{itemize}
     \item[(i)] $\mc{B}$ is stable by subobject, quotient and
       extension in $\mc{A}$,
     \item[(ii)] $\ind{\mc{B}}\subset \ind{\mc{A}}$ is stable by
       subobject, quotient and extension and contains $S$,
     \item[(iii)] $\ind{\mc{B}}$ has enough injectives.
  \end{itemize}
\end{prop}
Hence we get the following corollary.
\begin{cor}\label{injec}
  Let $\mc{S}$ be a small diagram in $\Db(\on{Ind}(\me{A}))$, i.e. 
  $\mc{S}\subset\Mor{\on{D^b}\on{Ind}(\mc{A})}$ is a set of 
  morphisms. Then there exists an essentially small full abelian 
  subcategory $\mc{B}\subset\mc{A}$ such that
  \begin{itemize}
     \item[(i)] $\mc{B}$ is stable by subobject, quotient and
       extension in $\mc{A}$,
     \item[(ii)] $\ind{\mc{B}}\subset \ind{\mc{A}}$ is stable by
       subobject, quotient and extension,
     \item[(iii)] $\ind{\mc{B}}$ has enough injectives,
     \item[(iv)] $\mc{S}$ is contained in the image of the natural 
        functor
          $$\on{D^b}(\on{Ind}(\mc{B})) \lra  \on{D^b}(\on{Ind}(\mc{A})). $$
  \end{itemize}
\end{cor}
\begin{proof}
Every morphism $(f:\,A\ra A')\in\mc{S}$ may be represented by a diagram
$$ \xymatrix{
      A \ar[r] & A'' & A'. \ar[l] }$$
Chose such a diagram for every $f\in\mc{S}$. Consider the set of 
objects of $\on{Ind}\mc{A}$ appearing in some complex in some diagram and 
apply Proposition \ref{injectiveobjects} to this set. 
\end{proof}
Combining Proposition \ref{propo} with Remark \ref{trueitis} and Corollary \ref{injec}
we get:
\begin{thm}
  Let $X$ be a paracompact locally compact topological space with a countable base 
  of open sets and consider a separated full subprestack 
  $\me{C}\subset\on{D^b}(I(k_{{*}}))$ that is defined by a local property.\\
  Then $\me{C}$ is a stack.
\end{thm}

\subsection{Limits and colimits in stacks}

Recall that if $\me{C}$ is a prestack on a topological space $X$ we denote
by $\rho_{_{VU}}$ the restriction functor for two open subsets 
$V\subset U\subset X$ and $\rho_p^U:\,\me{C}(U)\ra\me{C}_p$ the
canonical functor into the stalk at $p$.
\begin{defi}
 Let $\mc{I}$ be a small category. We say that $\me{C}$ admits limits (resp. 
 colimits) indexed by $\mc{I}$ if for every open subset $U\subset X$ the 
 category $\me{C}(U)$ admits limits (resp. colimits) indexed by $\mc{I}$ such
 that the restriction functors commute to these limits (resp. colimits).
\end{defi}
Let $\mc{I}$ be a finite category. It is easy to see that if $\me{C}$ 
admits limits (resp. colimits) indexed by $\mc{I}$, then for every $p\in X$ 
the category $\me{C}_p$ admits limits (resp. colimits) indexed by $\mc{I}$ and
the functor $\rho_p^U$ commutes to such limits (resp. colimits).\\
However, the converse is not true. We cannot know simply by looking at the stalks 
whether or not a prestack admits limits or colimits indexed by $\mc{I}$ (even if
$\me{C}$ is a stack).\\ 
If $\me{C}$ is separated we can at least see from the stalks whether or not 
a given object represents a limit or colimit indexed by a finite category. 
By duality we only need to consider the case of finite colimits.
\begin{lemma}\label{preplemma}
  Let $\me{C}$ be a separated prestack on a topological space $X$.\\
  Consider a finite category $\mc{I}$, an open subset $U\subset X$ 
  and a functor $\alpha:\,\mc{I}\ra\me{C}(U)$.\\
  Suppose given an object $L\in\Ob{\me{C}(U)}$ and morphisms
  $\sigma_i:\,\alpha(i)\ra L$ such that for any morphism $s:\,i\ra j$
  of $\Mor{I}$ we have $\sigma_j\circ\alpha(s)=\sigma_i$.\\ 
  Then the two following assertions are equivalent:
  \begin{itemize}
     \item[(i)] $(L,\{\sigma_i\}_{i\in I})$ is a colimit of
       $\alpha$ in $\me{C}(U)$ and for any open subset $V\subset U$ the pair
       $(L|_{_V},\{\sigma_i|_{_V}\}_{i\in I})$ is a colimit of 
       $\rho_{_{VU}}\alpha$ in $\me{C}(V)$.
     \item[(ii)] $(L,\{(\sigma_i)_p\}_{i\in I})$ is a colimit
       of $\rho^U_p\alpha$ in $\me{C}_p$ for all $p\in U$.
  \end{itemize}
\end{lemma}
\begin{proof}
Let $V\subset U$ be an open subset.\\
The object $L\in\Ob{\me{C}(U)}$ and the morphisms $\sigma_i|_{_V}$ define a
natural morphism of sheaves ($\me{C}$ is separated) for any object 
$A\in\Ob{\me{C}(V)}$
\begin{equation}\label{formulaa}
   \mc{H}\!\on{om}_{_{\me{C}|_{_V}}}(L|_{_V},A) \lra
   \prlim{i\in I}{\mc{H}\!\on{om}_{_{\me{C}|_{_V}}}(\alpha(i)|_{_V},A)}
\end{equation}
Since $\mc{I}$ is a finite category we have for every $p\in V$
$$(\prlim{i\in I}{\mc{H}\!\on{om}_{_{\me{C}|_{_V}}}(\alpha(i)|_{_V},A)})_p
  \simeq \prlim{i\in
    I}{\mc{H}\!\on{om}_{_{\me{C}|_{_V}}}(\alpha(i)|_{_V},A)_p}. $$
Hence the morphism \eqref{formulaa} induces in the stalks 
\begin{equation}\label{formulab}
   \mc{H}\!\on{om}_{_{\me{C}_p}}(L,A)\simeq  
   \mc{H}\!\on{om}_{_{\me{C}|_{_V}}}(L,A)_p \lra
   \prlim{i\in I}{\mc{H}\!\on{om}_{_{\me{C}|_{_V}}}(\alpha(i)|_{_V},A)_p}
\end{equation}  
Assertion (i) is clearly equivalent to the fact that the morphism \eqref{formulaa} 
is an isomorphism for all $V\subset U$ and any $A\in\me{C}(V)$. \\
Assertion (ii) is equivalent to the fact that the morphism \eqref{formulab} is an
isomorphism for all $V\subset U$, $A\in\me{C}(V)$ and $p\in V$.\\
Since $\me{C}$ is separated the morphism \eqref{formulaa} is an isomorphism if and 
only if for every $p\in V$ the morphism \eqref{formulab} is an isomorphism, which
proves the lemma.
\end{proof}
\begin{remark}
  \em Consider a prestack $\me{C}$, $\eta^\dagger:\,\me{C}\ra\me{C}^{\dagger}$ the 
  natural functor into the associated separated prestack (resp. 
  $\eta^{\ddagger}:\,\me{C}\ra \me{C}^{\ddagger}$ the natural functor into 
  the associated stack), $\mc{I}$ a finite category and $\alpha:\mc{I} 
  \ra\me{C}(U)$ a functor.\\ 
  Suppose that there is an object $L\in\me{C}(U)$ and morphisms 
  $\sigma_i:\,\alpha(i)\ra L$ such that condition (ii) of the Lemma \ref{preplemma} 
  is verified.\\
  Then Lemma \ref{preplemma} immediately implies that 
  $(\eta^\dagger(L)|_{_V},\{\eta^\dagger(\sigma_i)|_{_V}\})$ (resp. 
  $(\eta^\ddagger (L)|_{_V},\{\eta^\ddagger(\sigma_i)|_{_V}\})$) is a colimit 
  of $\rho_{VU}^\dagger\eta^\dagger\alpha$ (resp. $\rho_{VU}^\ddagger
  \eta^\ddagger\alpha$) in $\me{C}^\dagger(V)$ (resp. $\me{C}^\ddagger(V)$ 
  for all $V\subset U$.
\end{remark}
Now suppose that we are given a stack $\me{C}$, a finite category $\mc{I}$ and a 
functor $\alpha:\,\mc{I}\ra\me{C}(U)$. In order to check that there exists a colimit of
$\alpha$ in $\me{C}(U)$ we can apply Lemma \ref{preplemma}. However, in practical 
situations (as in Section 7.2) it is often difficult to establish the existence
of an object $L$ defined on $U$ that verfies condition (ii) of Lemma \ref{preplemma}.
Therefore we will use a refinement of Lemma \ref{preplemma} adapted to stacks which 
states that it is sufficent to prove the existence of the object $L$
locally on $U$.
\begin{prop}\label{colimitsandstalks}
  Let $\me{C}$ be a stack on a topological space $X$ and $\mc{I}$ be a
  finite category.\\
  Suppose that for every open subset $U\subset X$ and every functor
  $\alpha:\,\mc{I}\ra\me{C}(U)$ there exists an open covering
  $U=\bigcup_{j\in J}{U_j}$, objects $L_j\in\Ob{\me{C}(U_j)}$ and
  morphisms $\sigma^j_i:\,\alpha(i)|_{_{U_i}}\ra \me{C}(U_j)$
  verifying condition (ii) of Lemma \ref{preplemma}.\\
  Then $\me{C}$ admits colimits indexed by $\mc{I}$.
\end{prop}
\begin{proof}
Consider an open subset $U\subset X$ and a functor
$\alpha:\,\mc{I}\ra\me{C}(U)$. By hypothesis and Lemma \ref{preplemma} there
exists an open covering $U=\bigcup_{j\in J}{U_j}$ such that
$\rho_{_{U_j}}\alpha$ is representable in $\me{C}(U_j)$ by an object $L_j$ and the 
restriction to any smaller open subset $W\subset U_j$ commutes to these colimits. 
The conditions clearly imply that we may patch together the colimits 
$L_j\in\Ob{\me{C}(U_j)}$ to an object $L\in\Ob{\me{C}(U)}$. Now applying again 
Lemma \ref{preplemma} we see that $L$ is a colimit of $\alpha$ and that all 
restrictions commute to this colimit.
\end{proof}
\begin{cor}\label{limitsandstalks}
  Let $\me{C}$ be a prestack, $\eta:\,\me{C}\ra\me{C}^{\ddagger}$ the natural 
  functor into the associated stack and $\mc{I}$ be a finite category.\\
  Suppose that the stalks of $\me{C}$ admit colimits indexed by $\mc{I}$.\\
  Moreover we assume that for any open subset $U\subset X$ and any functor 
  $\alpha:\,\mc{I}\ra\me{C}(U)$ the following statement holds:\\
  For any point $p\in U$ there exists an open neighborhood $U_p\subset U$, 
  an object $L^p\in\me{C}(U_p)$ and morphisms 
  $\sigma^p_i:\,\alpha(i)|_{_V}\ra L^p$ such that condition (ii) of Lemma
  \ref{preplemma} is verified.\\
  Then $\me{C}^\ddagger$ admits colimits indexed by $\mc{I}$. 
\end{cor}
\begin{proof}
Let $\alpha:\,\mc{I}\ra\me{C}^\ddagger(U)$ be a functor. Consider the functors
$\rho^{\ddagger U}_p\alpha$ for all $p\in X$. Since $\mc{I}$ is finite there 
exists an open neighborhood $U_p$ of $p$ such that $\rho^{\ddagger U}_p\alpha$ factors 
through $\me{C}(U_p)$. Hence we can apply Proposition \ref{colimitsandstalks}.
\end{proof}
In particular we get the much weaker statement that if a prestack $\me{C}$ admits 
colimits indexed by a finite category $\mc{I}$ then $\me{C}^\ddagger$ admits 
colimits indexed by $\mc{I}$.

\subsection{A criterion for abelian stacks}

We can apply the results of the last paragraph to additive prestacks with abelian 
stalks. First recall that if $\me{C}$ is an additive prestack then $\me{C}^\ddagger$ 
is additive. 
\begin{thm}\label{abstack}
  Let $\me{C}$ be an additive prestack with abelian stalks. 
  Suppose that for every $p\in X$ and every morphism $f:\,A\ra B$ in $\me{C}_p$ 
  there exists an open neighborhood $U$ of $p$ such that $f$ may be represented 
  by a morphism $\tilde f:\,\tilde A\ra\tilde B$ in $\me{C}(U)$ and there 
  are morphisms  $K\ra \tilde A$, $\tilde B\ra K'$ such that
  $K\ra\tilde A$ is a kernel in $\me{C}_q$ and $\tilde B\ra K'$ is a
  cokernel in $\me{C}_q$ for any $q\in U$.\\
  Then $\me{C}^{\ddagger}$ is an abelian stack.
\end{thm}
\begin{proof}
Clearly the conditions of Proposition \ref{colimitsandstalks} and 
Corollary \ref{limitsandstalks} are satisfied for cokernels and kernels.
Hence $\me{C}^\ddagger$ admits cokernels and kernels.\\
Let $f:\,A\ra B$ be a morphism of $\me{C}^\ddagger(U)$ and consider the 
natural morphism $\coim f\ra\im f$. Since the categories of germs are abelian 
this morphism is an isomorphism in the stalks. Since $\me{C}^\ddagger$ is 
separated it is also an isomorphism in $\me{C}^\ddagger(U)$.
\end{proof}
\begin{cor}\label{critab}
  Let $\me{C}$ be an additive stack on $X$ such that all stalks are abelian
  categories. Then $\me{C}$ is an abelian stack if and only if for
  every morphism $f:\,A\ra B$ in $\me{C}_p$ there is an open
  neighborhood $U$ of $p$ such that $f$ may be represented by a
  morphism $\tilde f:\,\tilde A\ra\tilde B$ in $\me{C}(U)$ and there 
  are morphisms  $K\ra \tilde A$, $\tilde B\ra K'$ such that
  $K\ra\tilde A$ is a kernel in $\me{C}_q$ and $\tilde B\ra K'$ is a
  cokernel in $\me{C}_q$ for any $q\in U$.
\end{cor}
\begin{remark}
  \emph{Hence in order to verify that an additive stack is abelian it is 
  enough to verify that its stalks are abelian categories and that ``kernels 
  and cokernels are constant in some neighborhood''.}
\end{remark}
   
\section{Microlocalization of sheaves}

\subsection{Notations}

Let $\R^+$ denote the group of strictly positive real numbers and $\C^\times$
the group of non-zero complex numbers.\\
We will mainly work on a fixed complex manifold\begin{footnote}{All manifolds 
(complex or real) in this paper are supposed to be finite dimensional with a countable 
base of open sets.}\end{footnote} $X$ of complex dimension 
$\dim_{_{\C}}\!X=d_X$. Let $T^*X$ be its cotangent bundle and $T^*_XX$
the zero section. Set $\dot T^*X=T^*X\setminus T^*_XX$ and let $P^*X=
\dot T^*X/\C^{\times}$ be the projective cotangent bundle. We denote by 
$$ \gamma:\, \dot T^*X \lra P^*X  $$
the natural map.\\
If $\Lambda\subset \dot T^*X$ is a subset, we define the antipodal set $\Lambda^a$ as
$$ \Lambda^a=\Big\{(x;\xi)\ |\ (x;-\xi)\in\Lambda\Big\}, $$
and we set
$$ \R^+\Lambda=\Big\{(x;\xi)\in\dot T^*X\ |\ \exists \alpha\in\R^+\ (x;\alpha\xi)\in
   \Lambda\Big\}. $$
We define similarly $\C^\times\Lambda$. Hence $\C^\times \Lambda=\gamma^{-1}
\gamma(\Lambda)$. If $\Lambda=\{p\}$ is a point, we will write $\C^\times p$ instead
of $\C^\times\{p\}$.\\
We say that a subset $\Lambda\subset \dot T^*X$ is $\R^+$-conic (resp. 
$\C^\times$-conic) if it is stable under the action of $\R^+$ (resp. $\C^\times$), i.e.
if $\R^+\Lambda=\Lambda$ (resp. $\C^\times\Lambda=\Lambda$).\\
In the sequel, we will often deal with $\R^+$-conic subsets that are
only locally $\C^\times$-conic.
More precisely, a subset $\Lambda\subset\dot T^*X$ is called $\C^\times$-conic at $p\in\dot T^*X$ if there exists an open neighborhood $U$ of $p$ such that 
$U\cap \C^\times \Lambda=U\cap\Lambda$. Note that this definition 
still makes sense if $\Lambda$ is a germ of a subset at $p$. An open subset is always
$\C^\times$-conic at each $p\in U$.\\
Let $S\subset \dot T^*X$ be another subset, and suppose that $\Lambda$ is defined on
a germ of a neighborhood of $S$. Then we say that $\Lambda$ is $\C^\times$-conic
on $S$ if it is $\C^\times$-conic at every point of $S$. Clearly this is equivalent to
the statement that there exists an open neighborhood $U$ of $S$ such that $U\cap 
\C^\times\Lambda=U\cap\Lambda$. In particular, $\Lambda$ is $\C^\times$-conic on 
$\dot T^*X$ if and only if it is $\C^\times$-conic.\\
Finally we call the following easy topological lemma to the reader's attention.
\begin{lemma}
  Let $S\subset \dot T^*X$ be a $\C^\times$-conic set and $U\supset S$ an $\R^+$-conic
  open neighborhood. Then there exists a $\C^\times$-conic open set $V$ such that
  $S\subset V\subset U$.  
\end{lemma}
Now let us fix the conventions for sheaves. All sheaves considered here are sheaves 
of vector spaces over a given field $k$.\\
We will consider the following categories: 
\begin{itemize}
    \item[] $\on{D^b}(k_X)$ is the derived category of bounded complexes
      of sheaves of $k$ vector spaces,
    \item[] $\on{D^b_{\text{$\R$-c}}}(k_X)$ is the full subcategory of
      $\on{D^b}(k_X)$ whose objects have $\R$-constructible cohmology,
    \item[] $\on{D^b_{\text{$\C$-c}}}(k_X)$ is the full subcategory of
      $\Db_{\text{$\R$-c}}(k_X)$ whose objects have $\C$-constructible 
      cohomology,
    \item[] $\Perv(k_X)$ is the full  abelian subcategory of $\Db_{\Cc}(k_X)$ whose
      objects are perverse sheaves.   
\end{itemize}
We will not recall the construction of these categories here, for more details 
see for instance \cite{KS3}.

\subsection{Microlocalization of sheaves}

In this section we recall the construction and some properties of the 
microlocalization of $\Db(k_X)$ on a subset $S\subset T^*X$ (\cite{KS3}, Chapter VI) 
which we will then discuss from the (pre)stack-theoretical point of view. Note that 
all definitions and statements below which do not involve $\C^\times$-conic subsets of 
$T^*X$ are valid on a real manifold.\\
Recall that if $\me{F}\in\Db(k_X)$ then one can associate to $\me{F}$ a closed 
$\R^+$-conic involutive subset $\SS(\me{F})$ of $T^*X$ called the micro-support of 
$\me{F}$. The theory of the micro-support can be found in \cite{KS3}. It is the set 
of codirections in which 
$\me{F}$ ``does not propagate''. More precisely, a point 
$p\in T^*X$ is not a point of the micro-support if and only if there exists 
an open neighborgood $U$ of $p$ such that for any $x\in X$ and any real map $\psi$ 
of class $\on{C}^1$ with $\psi(x)=0$ and $(\on{d}\psi)_{x}\in U$ we have
$$\Dr\Gamma_{\{x\ |\ \psi(x)\geqs 0\}}(\me{F})_x\simeq 0.$$
Let $S\subset T^*X$ be an arbitrary subset. Set
$$ \mc{N}_S=\Big\{ \me{F}\in \on{D^b}(k_X)\ |
   \ \on{SS}(\me{F})\cap S=\varnothing\Big\}. $$
The micro-support is invariant under the shift functor. Further, if
$$ \xymatrix{
   {\me{F}'} \ar[r] & {\me{F}} \ar[r] & {\me{F}''} \ar[r] & {\me{F}'[1]} }$$
is a distinguished triangle in which $\me{F}'$ and $\me{F}$ are objects
of $\mc{N}_S$,  then it follows from the inclusion 
$$\on{SS}(\me{F''})\subset \on{SS}(\me{F})\cup \on{SS}(\me{F'}) $$
that $\me{F}''$ is also an object of $\mc{N}_S$. Hence $\mc{N}_S$ is a full 
triangulated subcategory of $\Db(k_X)$.\\
Note that if $x\in S\cap T^*_XX$ then $\me{F}\in\mc{N}_S$ implies $\me{F}\simeq 0$
in a neighborhood of $x$. 
\begin{defi}
  The microlocalization of $\Db(k_X)$ on $S$ is the localization of the triangulated
  category $\Db(k_X)$ by the full triangulated subcategory $\mc{N}_S$, i.e. 
  $$ \Db(k_X,S)=\Db(k_X)/\mc{N}_S. $$
  If $S=\{p\}$, we will write $\Db(k_X,p)$ for $\Db(k_X,\{p\})$. 
\end{defi}
Note that an object $\me{F}$ in $\Db(k_X,S)$ is isomorphic to zero if and
only if $\me{F}\oplus\me{F}[1]\in\mc{N}_{S}$ and since
$\on{SS}(\me{F}\oplus\me{F}[1])=\on{SS}(\me{F})$ this is equivalent to 
$\on{SS}(\me{F})\cap S=\varnothing$.\\
A morphism $\me{F}\ra \me{G}$ of $\on{D^b}(k_X)$ is called an isomorphism on 
$S$ if it is an isomorphism in $\on{D^b}(k_X,S)$. This is equivalent 
to the existence of a distinguished triangle in $\Db(k_X)$
$$ \xymatrix{
   {\me{F}} \ar[r] & {\me{G}} \ar[r] & {\me{H}} \ar[r]^+ &  }$$ 
with $\on{SS}(\me{H})\cap S=\varnothing$.\\
From this it follows easily that if $\me{F}\overset{\sim}{\ra}\me{G}$ is an isomorphism 
in $\Db(k_X,S)$, then $\SS(\me{F})\cap S=\SS(\me{G})\cap S$. Hence the micro-support of 
$\me{F}\in\Db(k_X,S)$ is well-defined in a germ of a neighborhood of $S$.\\ 
Let $\me{F},\me{G}\in\Db(k_X,S)$. By definition we have
$$ \on{Hom}_{\Db(k_X,S)}(\me{F},\me{G}) \simeq
    \varinjlim_{\me{F}'\overset{\sim}{\ra}\me{F} \above 0pt \text{on $S$}}
     \on{Hom}_{\Db(k_X)}(\me{F'},\me{G}). $$
We will have constant recourse to the following easy lemma.
\begin{lemma}\label{micro1}
  Let $S\subset \dot T^*X$ be any subset. Consider a morphism $\me{F}\ra \me{G}$ 
  of $\Db(k_X)$ that is an isomorphism on $S$. Then there 
  exists an $\R^+$-conic open neighborhood $U$ of $S$ such that $\me{F}
  \ra \me{G}$ is an isomorphism on $U$. In particular $\SS(\me{F})\cap U=
  \SS(\me{G})\cap U$.\\
  If moreover $S$ is $\C^\times$-conic, then we can choose $U$ to be 
  $\C^\times$-conic. 
\end{lemma}
\begin{proof}
By hypothesis there exists a distinguished triangle in $\Db(k_X)$
$$  \xymatrix{
    {\me{F}} \ar[r] & {\me{G}} \ar[r] & {\me{H}} \ar[r]^+ &  }$$
such that $\SS(\me{H})\cap S=\varnothing$. Since 
$\dot{\SS}(\me{H})$ is a closed $\R^+$-conic subset of $\dot T^*X$, the 
set $U=\complement \dot{\SS}(\me{H})$ is an open and $\R^+$-conic 
neighborhood of $S$ such that $\SS(\me{H})\cap U=\varnothing$.\\
Now suppose that $S$ is $\C^\times$-conic. To prove the last statement we use the 
fact that every $\R^+$-conic open neighborhood $V$ of $S$ contains a 
$\C^{\times}$-conic open neighborhood of $S$.
\end{proof}
  Recall that to any $\me{F},\me{G}\in\Db(k_X)$ we can associate the object 
  $\mu hom(\me{F},\me{G})\in\Db(k_{T^*X})$ (see \cite{KS3}, Chapter IV). This complex
  satisfies 
  $$ \supp (\mu hom(\me{F},\me{G}))\subset \SS(\me{F})\cap \SS(\me{G}). $$
  Therefore $\mu hom(\me{F},\me{G})|_{_S}$ is well-defined for 
  $\me{F},\me{G}\in\Db(k_X,S)$.\\
  For an arbitrary subset $S$ there is a natural morphism
  \begin{equation}\label{mor1} 
   \on{Hom}_{\Db(k_X,S)}(\me{F},\me{G})\lra
   \on{H}^0\!(S,\mu hom(\me{F},\me{G})). 
  \end{equation}
  Let us recall its construction. For any two objects $\me{F}_1,\me{F}_2\in\Db(k_X)$ we 
  have a canonical isomorphism 
  $$ \on{Hom}_{\Db(k_X)}(\me{F}_1,\me{F}_2)\simeq 
     \on{H}^0(T^*X,\mu hom(\me{F}_1,\me{F}_2))$$
  which defines a morphism
   $$ \on{Hom}_{\Db(k_X)}(\me{F}_1,\me{F}_2)\lra 
      \on{H}^0(S,\mu hom(\me{F}_1,\me{F}_2)).$$
  Now if $\me{F}'\ra \me{F}$ is an isomorphism on $S$ we get an induced isomorphism
  $$  \on{H}^0(S,\mu hom(\me{F},\me{G})) \overset{\sim}{\lra}
      \on{H}^0(S,\mu hom(\me{F}',\me{G})). $$
  Thus we get morphisms 
   $$ \on{Hom}_{\Db(k_X)}(\me{F}',\me{G})\lra 
      \on{H}^0(S,\mu hom(\me{F},\me{G})).$$
  which induce the morphism (\ref{mor1}).\\
  There is a well-known situation in which this morphism is an isomorphism 
  (\cite{KS3}, Theorem 6.1.2).
\begin{prop}\label{prop51}
  Let $p\in T^*X$ and $\me{F},\me{G}\in\on{D}^b(X,p)$. Then the morphism (\ref{mor1})
$$ \on{Hom}_{\on{D}^b(X,p)}(\me{F},\me{G})\lra
   \on{H}^0\!\mu hom(\me{F},\me{G})_p $$
  is an isomorphism.
\end{prop}
The idea is to calculate both sides by using microlocal cut-off functors. We will
show that such a strategy works in the case of a closed cone in $\dot T^*_xX$, 
$x\in X$.

However, the morphism (\ref{mor1}) is not an isomorphism in general (cf. \cite{KS3},
Exercise VI.6 which gives a counter-example on an open subset). 
 
\begin{remark}
  \emph{The correspondance 
    $$ T^*X\supset U\mapsto \Db(k_X,U) $$ 
  defines a prestack on $T^*X$, which we will usually denote by $\Db(k_X,*)$.\\
  Therefore $\gamma_*\Db(k_X,*)|_{_{\dot T^*X}}$ defines a prestack on $P^*X$.} 
\end{remark}
\begin{prop}\label{germformula}
  Let $S$ be a subset of $\dot T^*X$. Then the natural functor
  $$ \twodilim{\text{$S\subset U\subset \dot T^*X$} \above 0pt
   \text{$U$ $\R^+$-conic}} \Db(k_X,U) \lra \Db(k_X,S) $$
  is an equivalence.
\end{prop}
\begin{proof}
The functor is obviously essentially surjective. Let us show that it is  
fully faithful. \\
Let $\me{F},\me{G}\in\Db(k_X)$. By Lemma \ref{micro1} we get 
\begin{align*}
   \Hom{\twodilim{\text{$S\subset U\subset \dot T^*X$} \above 0pt
   \text{$U$ $\R^+$-conic}}{\Db(k_X,U)}}(\me{F},\me{G})&\simeq
   \underset{\text{$S\subset U\subset \dot T^*X$} \above 0pt
   \text{$U$ $\R^+$-conic}}{\varinjlim}\Hom{\Db(k_X,U)}(\me{F},\me{G})\\
   &\simeq
   \underset{\text{$S\subset U\subset \dot T^*X$} \above 0pt
   \text{$U$ $\R^+$-conic}}{\varinjlim}
   \underset{\text{$\me{F}'\overset{\sim}{\ra} \me{F}$} \above 0pt
   \text{on $U$}}{\varinjlim}
   \Hom{\Db(k_X)}(\me{F}',\me{G})\\
   &\simeq
   \underset{S\subset U\subset \dot T^*X}{\varinjlim}
   \underset{\text{$\me{F}'\overset{\sim}{\ra}\me{F}$} \above 0pt
   \text{on $U$}}{\varinjlim}\Hom{\Db(k_X)}(\me{F}',\me{G})
   \simeq \Hom{\Db(k_X,S)}(\me{F},\me{G}).
\end{align*}
\end{proof}
\begin{cor}\label{cgermformula}
   Let $S$ be a $\C^{\times} $-conic subset of $\dot T^*X$ and $p\in T^*X$. 
   \begin{itemize}
      \item[(i)] The natural functor
  $$ \twodilim{\text{$S\subset U\subset \dot T^*X$} \above 0pt
     \text{$U$ $\C^\times$-conic}} \Db(k_X,U) \lra \Db(k_X,S) $$
      is an equivalence.
      \item[(ii)] The natural functor 
      $$ \Db(k_X,*)_p\lra \Db(k_X,p) $$
      is an equivalence. If moreover $p\in\dot T^*X$ then
  $$ \gamma_*(\Db(k_X,*)|_{\dot T^*X})_{\gamma(p)}
     \simeq \Db(k_X,\gamma^{-1}\gamma(p))=\Db(k_X,\C^{\times} p). $$ 
   \end{itemize}
\end{cor}
\begin{proof}
Any $\R^+$-conic neighborhood of a $\C^{\times} $-conic subset contains a 
$\C^{\times} $-conic neighborhood. Thus part (i) of the corollary follows from a 
cofinality argument applied to the result of Proposition \ref{germformula} and part 
(ii) follows from Proposition \ref{germformula} and (i).
\end{proof}

\subsection{Refined microlocal cut-off}

Let us recall the basic idea of a microlocal cut-off functor.\\
Let $X$ be a finite dimensional real vector space, $U\owns 0$ a relatively compact 
open neighborhood of $0$ and consider an open cone $\gamma\subset X^*$. 
\begin{defi}\label{deficutoff}
A microlocal cut-off functor on $U\times\gamma$ is a functor 
 $$\Phi_{U,\gamma}:\,\Db(k_X)\ra\Db(k_X)$$ 
such that 
\begin{itemize}
   \item[(i)] $\SS(\Phi_{U,\gamma}\me{F})\subset X\times \ol{\gamma}$ 
     (the micro-support has been cut off at $\ol{\gamma}$),
   \item[(ii)] $\SS(\Phi_{U,\gamma}(\me{F}))\cap U\times \gamma=
      \SS(\me{F})\cap U\times\gamma$,
   \item[(iii)] $\Phi_{U,\gamma}$ is equipped with a morphism of functors 
      $\alpha:\,\Phi_{U,\gamma}\ra \on{Id}$ such that $\alpha$ induces an 
      isomorphism in $\Db(k_X,U\times\gamma)$ which can be 
      visualized by
$$ \xymatrix@C=2cm{
      {\Db(k_X)} \ar[d] \ar[r]^{\Phi_{U,\gamma}} & {\Db(k_X)} \ar[d] \\
       {\Db(k_X,U\times\gamma)} \ar[r]_{\Phi_{U,\gamma}\simeq Id} &
       {\Db(k_X,U\times\gamma)}. }$$
\end{itemize}
\end{defi}
Note that condition (iii) implies (ii).\\
If the cut-off functor $\Phi_{U,\gamma}$ allows us to estimate the micro-support
of $\Phi_{U,\gamma}(\me{F})$ in the fiber $\{0\}\times X^*$, we usually call it
a refined microlocal cut-off.\\
A cut-off functor is easily constructed in the case of a convex open cone (see
Proposition 5.2.3 of \cite{KS3}). A generalization to non-convex cones is stated 
in Exercise V.8 of \cite{KS3} (for a proof see \cite{D'A}). These tools will allow us 
in Section 3.4 to calculate sections of $\mu hom$ along a complex line (or more 
generally along closed cones of $T^*_xX$ where $x\in X$). The result will
imply that the morphism \eqref{mor1} is an isomorphism in this case 
(see Section 3.4).\\ 
Later we will need to construct a functor $\Db_{\Cc,\Lambda}(k_X,\C^{\times} p)\ra 
\Db_{\Cc}(k_X,\pi(p))$ if $\Lambda$ is in generic position at $p$ (Section 6.1).
For this purpose we will need the refined microlocal cut-off of \cite{D'A}. 
It is an extension of the classical ``refined microlocal cut-off lemma'' 
(Proposition 6.1.4 of \cite{KS3}) to non-convex cones with a good estimate for 
the micro-support.\\
Let us recall the cut-off functor of \cite{KS3}, Exercise V.8. Let $X$ be a real, 
finite dimensional vector space, $\dot X=X\setminus\{0\}$, $U\subset X$ an open 
subset and $\gamma\subset X$ an open cone.\\
We have the following natural morphisms:
$$ \xymatrix@R=.5cm@C=.5cm{
        & X\times X \ar[dl]_{q_1} \ar@<.7ex>[d]^{s} 
         \ar@<-.7ex>[d]_{\tilde s} \ar[dr]^{q_2} & \\
       X & X & X } $$
where $q_1$ and $q_2$ are the natural projections and
\begin{eqnarray*}
   s:\ X\times X\lra X & ; & (x,y)\mapsto x+y, \\
   \tilde s:\ X\times X\lra X & ; & (x,y)\mapsto x-y. 
\end{eqnarray*}
We define the functor $\Phi_{U,\gamma}$ by setting for any $\me{F}\in\Db(k_X)$
$$ \Phi_{U,\gamma}(\me{F})=k_{\gamma^a}^{\sct{\wedge}}*\me{F}_U=
   \on{R}\!s_{!}(q_1^{-1}k_{\gamma^a}^{\sct{\wedge}}\otimes
    q_2^{-1}\me{F}_U). $$
It can be shown that $\Phi_{U,\gamma}$ is a microlocal cut-off functor in the sense of
Definition \ref{deficutoff}.\\
Let us add two easy lemmas which will be useful in the next section.
\begin{lemma}
  Let $\me{F},\me{G}\in\Db(k_X)$. Then we have a canonical isomorphism
  $$ \on{R}s_!(q_1^{-1}\me{F}\otimes q_2^{-1}\me{G})\simeq
     \on{R}q_{1!}(\tilde s^{-1}\me{F}\otimes q_2^{-1}\me{G}) $$
\end{lemma}
\begin{proof}
Consider the map
$$ \varphi:\ X\times X \lra X\times X \quad ; \quad (x,y)\mapsto (x+y,y) $$
We obviously have the formulas:
$$ \tilde s\circ \varphi=q_1 \qquad 
  q_1\circ \varphi=s \qquad
  q_2\circ \varphi=q_2 $$
Hence we get the chain of natural isomorphisms
\begin{align*}
   \on{R}s_!(q_1^{-1}\me{F}\otimes q_2^{-1}\me{G})&\simeq
   \on{R}q_{1!}\on{R}\varphi_!(\varphi^{-1}{\tilde{s}}^{-1}\me{F}
      \otimes \varphi^{-1}q_2^{-1}\me{G})\simeq
    \on{R}q_{1!}(\on{R}\varphi_!\varphi^{-1}{\tilde{s}}^{-1}\me{F}
      \otimes q_2^{-1}\me{G})\\
   &\simeq \on{R}q_{1!}({\tilde{s}}^{-1}\me{F}
      \otimes q_2^{-1}\me{G}).
\end{align*}
\end{proof}
\begin{lemma}\label{lemma7.2}
  Let $\gamma_1,\gamma_2\subset X$ be two open cones and $U\subset X$
  open. Then there is a natural distinguished triangle
  $$ \xymatrix{
    \Phi_{U,\gamma_1\cap\gamma_2}(\me{F}) \ar[r] &
    \Phi_{U,\gamma_1}(\me{F})\oplus \Phi_{U,\gamma_2}(\me{F}) \ar[r] &
    \Phi_{U,\gamma_1\cup \gamma_2}(\me{F}) \ar[r]^(.7){+} & }$$
\end{lemma}
\begin{proof}
Follows immediately from the distinguished triangle
$$ 
\xymatrix{
  k_{\gamma_1^a\cap \gamma_2^a} \ar[r]  &
  k_{\gamma_1^a} \oplus k_{\gamma_2^a} \ar[r] & 
  k_{\gamma_1^a\cup \gamma_2^a} \ar[r]^(.6){+} & }$$
by applying the triangulated functors $(\,\cdot\,)^\wedge$ and 
$(\,\cdot\,)*\me{F}_{U}$.
\end{proof}

Next, we recall D'Agnolo's condition under which the cut-off functor 
$\Phi_{U,\gamma}$ is refined. These results will not be needed until 
Section 6.\\

\begin{defi}
  If $\gamma\subset X^*$ is an open cone, set 
    $$\partial^\circ\gamma=
  \pi\chi(\SS(\C_{\gamma})\setminus\{0;0\})$$ 
  where $\chi:\,T^*X^*\ra T^*X$ is defined by $\chi(\xi;x)=(x;-\xi)$.\\
  One says that $(U,\gamma)$ is a refined cutting pair at $0$, if $U\subset
  X$ is a relatively compact open neighborhood of $0$ and for any
  $x\in\partial U\cap \partial^\circ\gamma$ there exists $\xi\in\dot
  X$ such that $\on{N}_x^*(U)=\R_{\geqs 0}\xi$ and $\chi(\SS(k_{_\gamma}))\cap 
  \pi^{-1}(x)=\R_{\leqs 0}\xi$.
\end{defi}
This allows one to give a good estimate of the micro-support of 
$\Phi_{U,\gamma}(\me{F})$ at $0$.
\begin{prop}
  Let $(U,\gamma)$ be a refined cutting pair at $0$. Then
$$ \SS(\Phi_{_{U,\gamma}}(\me{F}))\cap \pi^{-1}(0)\subset
   \Big\{\xi\in\gamma\ |\ (0;\xi)\in\SS(\me{F})\Big\}\cup
  \Big\{\xi\in\partial\gamma\ |\ \exists x\in \overline{U} :
  (x,\xi)\in\SS(\me{F})\Big\}.$$
\end{prop}
Let us add a useful corollary:
\begin{cor}\label{kleineergaenzung}
   Let $(U,\gamma)$ be a refined cutting pair at $0$ and suppose that
   $\SS(\me{F})\cap (\ol{U}\times \partial\gamma)=\varnothing$. Then there
   exists an open neighborhood $V$ of $0$ such that
   $$ \SS(\Phi_{U,\gamma}(\me{F}))\cap \pi^{-1}(V)=
      \SS(\me{F})\cap (V\times\gamma). $$
\end{cor}
\begin{proof}
Since $\SS(\Phi_{U,\gamma}(\me{F}))\subset X\times\ol{\gamma}$ and
$\SS(\Phi_{U,\gamma}(\me{F}))\cap (U\times\gamma)=\SS(\me{F})\cap (U\times
\gamma)$, it is enough to show that $\SS(\Phi_{U,\gamma}(\me{F}))\cap V
\times\partial \gamma=\varnothing$ for some neighborhood $V$ of $0$. 
D'Agnolo's estimate of the microsupport implies that this is at least 
true at $0$.\\
Now suppose that such a neighborhood $V$ does not exist. Then we can construct
a sequence $(x_n,\xi_n)$ such that $x_n\lra 0$ and $\xi_n\in
\SS_{x_n}(\Phi_{U,\gamma}(\me{F}))\cap\partial\gamma$.
Since both sets are invariant by $\R_{\geqs 0}$ we can assume that $|\xi_n|=1$,
hence by extracting a subsequence we can suppose that $\xi_n\lra \xi$. 
Since $\partial\gamma$ and $\SS(\Phi_{U,\gamma}(\me{F}))$ are 
closed we get by the estimate of the micro-support that there exists
$x\in\ol{U}$ such that $(x,\xi)\in\SS(\me{F})$ which is impossible by 
hypothesis.
\end{proof}
Finally, let us state an existence lemma for refined cutting pairs.
\begin{lemma}
  Let $X$ be a real vector space and $L\subset X$ a subspace of $X$. Then 
  there exists a fundamental system of open conic neighborhoods $\gamma$
  of $(\dot T^*_LX)_0$ such that for each $\gamma$ there exists a fundamental 
  system of open neighborhoods $U$ of $0$ in $X$ such
  that $(U,\gamma)$ is a refined cutting pair. 
\end{lemma}
\begin{proof}
This lemma is shown during the proof of Corollary 3.4 using Lemma 3.3 
in \cite{D'A}.
\end{proof}
We are now ready to give the refined version of the microlocal cut-off lemma along 
a linear subspace.
\begin{prop}[Refined microlocal cut-off]\emph{ }\\
  Let $X$ be a real vector space and $L\subset X$ a subspace. Let
  $\me{F}\in \Ob{\on{D}^{\on{b}}}(k_X)$.\\
  Then there exists a fundamental system of open conic neighborhoods
  $\gamma$ of $\dot (T^*_LX)_0$ in $T^*_0X=X^*$ and for each $\gamma$ there
  exists a fundamental system of open neighborhoods $U$ of $0$, such that 
  \begin{itemize}
    \item[(i)] the natural morphism $u:\,\Phi_{U,\gamma}(\me{F})\ra \me{F}$ 
      induces an isomorphism in $\on{D}^{\on{b}}(k_X,U\times\gamma)$.
    \item[(ii)] $\SS(\Phi_{U,\gamma}(\me{F}))\subset X\times \ol{\gamma}$ 
    \item[(iii)] $\dot{\on{SS}}(\Phi_{U,\gamma}(\me{F}))\cap 
      \partial\gamma\cap \pi^{-1}(0)\subset \Big\{\xi\in
      \partial\gamma\ |\ \exists x\in \overline{U} :(x,\xi)\in
      \SS(\me{F})\Big\}$ where we have set $\dot{\on{SS}}(\Phi_{U,\gamma}(\me{F}))= 
      \on{SS}(\Phi_{U,\gamma}(\me{F}))\cap \dot T^*X$.
  \end{itemize} 
\end{prop}

\subsection{Morphisms in $\Db(k_X,\{x\}\times\dot\delta)$}

In this section we will calculate the morphisms in the category
$\Db(k_X,\{x\}\times\dot\delta)$ where $x\in X$ and $\delta\subset T^*_xX$ is a closed
cone. More precisely we will show that the natural morphism
$$  \on{Hom}_{_{\Db(k_X,\{x\}\times\dot\delta)}}\big(\me{F},\me{G}\big) \lra
     \on{H}^0\!\big(\{x\}\times\dot\delta,\mu hom(\me{F},\me{G})\big) $$
is an isomorphism. In order to prove this, we will consider the composition 
$$  \inlim{U,\gamma}{\on{H}^0\!\Dr \on{Hom}\big(\Phi_{U,\gamma}(\me{F})_U,\me{G}\big)}
 \lra \on{Hom}_{_{\Db(k_X,\{x\}\times\dot\delta)}}\big(\me{F},\me{G}\big) \lra
     \on{H}^0\!\big(\{x\}\times\dot\delta,\mu hom(\me{F},\me{G})\big)
   \qquad \text{$(*)_\delta$}   $$
and show that it is an isomorphism. Here $U$ runs through the family of relatively
compact open neighborhoods of $0$ and $\gamma$ through the set of open
cones containing $\dot\delta$.
\begin{prop}\label{thecalculation}
  Let $X$ be a real vector space and consider a closed convex
  proper cone $\delta\subset X$. Then the natural morphism
  \begin{equation*} 
    \inlim{U,\gamma}{\on{H}^n\!\Dr \on{Hom}\big(\Phi_{U,\gamma}(\me{F})_U,\me{G}\big)}
   \lra \on{H}^n\!\big(\{0\}\times\dot\delta,\mu hom(\me{F},\me{G})\big) 
    \qquad \qquad \qquad \qquad \text{$(*)_\delta$}
  \end{equation*}
  is an isomorphism. Here $U$ runs through the family of relatively
  compact open neighborhoods of $0$ and $\gamma$ through the set of open
  cones containing $\dot\delta$.
\end{prop}
\begin{proof}
First note that since $\mu hom(\me{F},\me{G})$ is conic, we have
$$ \on{H}^n\!\big(\{0\}\times\dot\delta,\mu hom(\me{F},\me{G})\big)
  \simeq \inlim{U,\gamma}
  {\on{H}^n\!\big(U\times\gamma,\mu hom(\me{F},\me{G})\big).}$$
Since $\gamma$ is open, convex and proper, we have by \cite{KS3}, Theorem 4.3.2
\begin{align*}
  \on{H}^n\!&\big(U\times\gamma,\mu hom(\me{F},\me{G})\big)\simeq
   \inlim{V,Z}{\on{H}_{Z\cap V}^n\!\big(V,\on{R}
     \me{H}om(q_2^{-1}\me{F},q_1^{!}\me{G})\big)}
\end{align*}
where $V$ runs through the family of open subsets of
$T_{\Delta_{_X}}(k_X\times X)\simeq X\times X$
such that $V\cap \Delta_X=U$ (i.e. $q_1(V)=U$) and $Z$ through the 
family of closed subsets such that the inclusion
$\on{C}_{_{\Delta_{_X}}}\!Z\subset U\times\gamma^\circ$ holds
\begin{footnote}{Here $\on{C}_{\Delta_X}Z$ denotes the normal cone to $Z$ along the
diagonal $\Delta_X$ (see \cite{KS3}, Definiton 4.1.1) and $\gamma^\circ$ denotes
the polar cone of $\gamma$, i.e. $$ \gamma^\circ=\{x\in X^*\ |\ \langle y,x
\rangle\geqs 0\ \text{for all $y\in\gamma$}.\}$$}\end{footnote}.\\ 
We have the following chain of isomorphisms
\begin{align*} 
  \on{R}\Gamma_{Z\cap
    V}&\big(V,\on{R}\me{H}om(q_2^{-1}\me{F},q_1^{!}\me{G})\big)\simeq
   \on{R}\Gamma\big(V,\on{R}\Gamma_{Z\cap V}
    \on{R}\me{H}om(q_2^{-1}\me{F},q_1^{!}\me{G})\big) \\
  & \simeq \on{R}\Gamma\big(V,
      \on{R}\me{H}om((q_2^{-1}\me{F})_{Z\cap V},q_1^{!}\me{G})\big)
    \simeq \on{R\Gamma}\big(U\times X,
    \on{R}\me{H}om((q_2^{-1}\me{F})_{Z\cap V},q_1^{!}\me{G})\big)\\
  &\simeq \on{R\Gamma}\big(U,
  \on{R}\me{H}om(\on{R}q_{1!}(q_2^{-1}\me{F})_{Z\cap V},\me{G}\big)\simeq
  \on{RHom}\big((\on{R}q_{1!}(k_{Z\cap V}
   \otimes q_2^{-1}\me{F})_U,\me{G}\big). 
\end{align*}   
Hence we have
$$ \on{H}^n\!\big(\{0\}\times\dot\delta,\mu hom(\me{F},\me{G})\big)
  \simeq \inlim{U,\gamma}{\inlim{V,Z}{\on{H}^n\Dr\on{Hom}\big((\on{R}q_{1!}
    (k_{Z\cap V}\otimes q_2^{-1}\me{F}))_U,\me{G}\big).}} $$
Now fix $U,\gamma$ and $V,Z$. Then $V$ contains a small relatively
compact open neighborhood of $0$ of type $U'\times U'$. Moreover we
may assume by cofinality that $Z$ is of the form ${\tilde s}^{-1}\gamma^\circ$ 
in a neighborhood of $0$. Hence $Z\cap V$ contains 
${\tilde s}^{-1}\gamma^\circ\cap U'\times U'$. We can
therefore remove the second limit by replacing $V\cap Z$ with
${\tilde s}^{-1}\gamma^\circ\cap U\times U$. Then we get 
\begin{align*} 
\on{RHom}(&(\Dr q_{1!}(k_{Z\cap V}
    \otimes q_2^{-1}\me{F}))_U,\me{G})\simeq
 \on{RHom}((\Dr q_{1!}({\tilde s}^{-1}k_{\gamma^\circ}\otimes
    q_2^{-1}\me{F}_U))_U,\me{G}) \\
  &\simeq 
  \on{RHom}((\on{R}s_{!}(q_1^{-1}k_{\gamma^a}^{\sct{\wedge}}\otimes
    q_2^{-1}\me{F}_U))_U,\me{G}) \simeq
   \on{RHom}(\Phi_{U,\gamma}(\me{F})_U,\me{G}).
\end{align*}
Therefore
$$ \on{H}^n\!\big(\{0\}\times\dot\delta,\mu hom(\me{F},\me{G})\big)
  \simeq
  \inlim{U,\gamma}{\on{H}^n\Dr\on{Hom}
   \big(\Phi_{U,\gamma}(\me{F})_U,\me{G})\big)}. $$
\end{proof}
\begin{lemma}
  For any closed cone $\delta\subset X$ consider the morphism $(*)_\delta$. 
  Let $\delta_1,\delta_2\subset X$ be two closed 
  cones such that the morphisms $(*)_{\delta_1},(*)_{\delta_2}$ and 
  $(*)_{\delta_1\cap\delta_2}$ are isomorphisms.\\
  Then $(*)_{\delta_1\cup\delta_2}$ is an isomorphism.
\end{lemma}
\begin{proof}
By Lemma \ref{lemma7.2} we get a morphism of distinguished triangles such
that the vertical morphisms are given by $(*)_{\delta_1\cap\delta_2}$, 
$(*)_{\delta_1}\oplus (*)_{\delta_2}$ and $(*)_{\delta_1\cup\delta_2}$.
Then the lemma follows from the Five Lemma.
\end{proof}
\begin{prop}\label{newresult}
  Let $X$ be a real vector space and consider a closed 
  cone $\delta\subset X$. Then the natural morphism
  $$ \inlim{U,\gamma}{\on{H}^n\Dr\on{Hom}\big(\Phi_{U,\gamma}(\me{F})_U,
       \me{G}\big)}\lra
     \on{H}^n\!\big(\{0\}\times\dot\delta,\mu hom(\me{F},\me{G})\big)
   \qquad\qquad\qquad \text{$(*)_\delta$} $$
  is an isomorphism. Here $U$ runs through the family of relatively
  compact open subsets of $0$ and $\gamma$ through the set of open
  cones containing $\dot\delta$.
\end{prop}
\begin{proof}
First suppose that $\delta$ can be written as a finite union of closed convex 
proper cones. Note that the intersection of two proper, closed, convex cones is 
again proper, closed and convex. Therefore, if $\delta'$ is the union of $n$ closed, 
convex, proper cones then the intersection of a closed, convex, proper cone 
with $\delta'$ can be written as a union of $n$ closed convex proper cones.
Using Proposition \ref{thecalculation} and the previous lemma we can then easily show
the proposition by induction on the number of closed, convex, proper cones that cover
$\delta$.\\
Now let us consider the general case. Every closed cone is a decreasing intersection 
of closed cones $\delta_i$ that can be covered by a finite number of closed convex 
proper cones\begin{footnote}{The proof is done by a simple compacity argument in 
$\dot X/\R^+$. For every $p\in \delta$ choose a closed convex proper cone with 
angle $\varepsilon$ that contains $p$. Then a finite number of these cones, say 
$\gamma_1,\gamma_2,\ldots,\gamma_n$, will cover $\delta$ and we set 
$\delta_1=\underset{i=1,\ldots,n}{\bigcup}{\gamma_i}$. The next cone $\delta_2$ 
is constructed by choosing for each point a closed convex proper cone with angle 
$\varepsilon/2$. Again a finite number will cover $\delta$, say 
$\gamma_1',\ldots,\gamma_m'$. Then we define $\delta_2$ as the 
union of all intersections $\gamma_i\cap\gamma'_j$ and we proceed by induction. 
It is clear by construction that the intersection of the $\delta_i$ is the 
cone $\delta$.}\end{footnote}. Then
$(*)_\delta=\varinjlim {(*)_{\delta_i}}$ is an isomorphism.
\end{proof}
\begin{thm}\label{cutoffprop}
 Let $X$ be a real manifold, $x\in X$ and $\delta \subset T^*_xX$ a closed
 cone. Then the natural morphism
 $$  \on{Hom}_{_{\Db(k_X,\{x\}\times\dot\delta)}}\big(\me{F},\me{G}\big) \lra
     \on{H}^0\!\big(\{x\}\times\dot\delta,\mu hom(\me{F},\me{G})\big) $$
  is an isomorphism. 
\end{thm} 
\begin{proof}
By the last proposition we know that the morphism $(*)_\delta$ is an isomorphism. 
Therefore the morphism of the theorem is surjective.\\
Let us prove that it is injective.\\ 
Let $\me{F}\ra\me{G}$ be a morphism of $\Db(k_X,\{x\}\times\dot\delta)$ that is zero 
in $\on{H}^0(\{x\}\times\dot\delta,\mu hom(\me{F},\me{G}))$. Then we may represent
this morphism by a morphism $\me{F}'\ra\me{G}$ in $\Db(k_X)$ and a morphism 
$\me{F}'\ra\me{F}$ that is an isomorphism on $\{x\}\times \dot \delta$. We get a 
commutative diagram
$$ \xymatrix{  
    {\on{Hom}_{_{\Db(k_X)}}\big(\me{F}',\me{G}\big)} \ar[r] \ar[d]_{Id} &  
    {\on{Hom}_{_{\Db(k_X,\{x\}\times\dot\delta)}}\big(\me{F},\me{G}\big)} \ar[r] 
       \ar[d]^{\sim} & 
    {\inlim{U,\gamma}{\Dr^0\on{Hom}\big(\Phi_{U,\gamma}
     (\me{F})_U,\me{G}\big)}} \ar[d]^{\sim}  \\
     {\on{Hom}_{_{\Db(k_X)}}\big(\me{F}',\me{G}\big)} \ar[r] &  
    {\on{Hom}_{_{\Db(k_X,\{x\}\times\dot\delta)}}\big(\me{F}',\me{G}\big)} \ar[r] & 
    {\inlim{U,\gamma}{\Dr^0\on{Hom}\big(\Phi_{U,\gamma}
     (\me{F}')_U,\me{G}\big).}} 
     }$$
Using the diagram, we see that there exists $(U,\gamma)$ such that 
$\Phi_{U,\gamma}(\me{F}')_U\ra \me{F}'\ra\me{G}$ is the zero map in $\Db(k_X)$. 
But $\Phi_{U,\gamma}(\me{F}')_U\ra\me{F}'$ is an isomorphism on 
$\{x\}\times\dot\delta$, and therefore $\me{F}'\ra\me{G}$ represents the zero 
morphism in $\Db(k_X,\{x\}\times\dot\delta)$.
\end{proof}

\section{Microlocalization of constructible sheaves}

\subsection{Microlocalization of $\R$-constructible sheaves}

Consider the full triangulated subcategory  $\Db_{\Rc}(k_X)\subset\Db(k_X)$
and a subset $S\subset T^*X$. There are two obvious ways to define the 
microlocalization of the derived category of $\R$-constructible sheaves on $S$ 
that we recall now. Set
$$ \mc{N}_{\Rc,S}=\mc{N}_S\cap\Db_{\Rc}(k_X). $$
Then the inclusion $\Db_{\Rc}(k_X)\subset \Db(k_X)$ induces a functor
$$ \Db_{\Rc}(k_X)/\mc{N}_{\Rc,S}\lra \Db(k_X)/\mc{N}_S=\Db(k_X,S). $$
Clearly the objects of the image of this functors are complexes in $\Db(k_X,S)$
with $\R$-constructible cohomology. But the functor is not fully faithful
in general. For our purpose it will be convenient to work with the category
$\Db_{\Rc}(k_X)/\mc{N}_{\Rc,S}$ as does Andronikof in \cite{An1},\cite{An2}.
\begin{defi}We set
  $$ \Db_{\Rc}(k_X,S)=\Db_{\Rc}(k_X)/\mc{N}_{\Rc,S}. $$
\end{defi}
\begin{remark}
  \em   Let us emphasize again that the natural functor
  \begin{equation}  \Db_{\Rc}(k_X,S) \lra \Db(k_X,S) \label{comparison} 
  \end{equation}
  is not fully faithful. Hence although the objects of $\Db_{\Rc}(k_X,S)$
  have $\R$-constructible cohomology, the category can not be identified to
  the full subcategory of $\Db(k_X,S)$ whose objects have $\R$-constructible 
  cohomology. More precisely, let $\me{F},\me{G}\in\Db_{\Rc}(k_X)$.
  \begin{itemize}
     \item[(i)] A morphism $\me{F}\ra\me{G}$ in $\Db_{\Rc}(k_X,S)$ is represented by a 
       morphism $\me{F}'\ra\me{G}$ where $\me{F}'$ has $\R$-constructible cohomology
       sheaves and is isomorphic to $\me{F}$ on $S$.
     \item[(ii)] A morphism $\me{F}\ra\me{G}$ in $\Db(k_X,S)$ is represented by a 
       morphism $\me{F}'\ra\me{G}$ where $\me{F}'\in\Db(k_X)$ is isomorphic to
       $\me{F}$ on $S$.
  \end{itemize} 
  Therefore it is not obvious to compare morphisms in these two categories.\\ 
  Nevertheless note that for any $\me{F}\in\Db_{\Rc}(k_X,S)$ we have by definition 
  $\me{F}\simeq 0$ if and only if $\SS(\me{F})\cap S=\varnothing.$ Hence 
  if $\me{F}\ra\me{G}$ is a morphism in $\Db_{\Rc}(k_X)$ we get 
  that $\me{F}\ra\me{G}$ is an isomorphism in $\Db_{\Rc}(k_X,S)$ if and 
  only if it is an isomorphism in $\Db(k_X,S)$, hence if and only if there is a 
  distinguished triangle in $\Db_{\Rc}(k_X)$ 
  $$ \xymatrix{
       {\me{F}}\ar[r] & {\me{G}}\ar[r] & {\me{H}} \ar[r]^+ & }$$ 
  such that $\SS(\me{H})\cap S=\varnothing$. More generally we get
\end{remark}
\begin{prop}
  The natural functor  
  $$ \Db_{\Rc}(k_X,S) \lra \Db(k_X,S) $$ 
  is conservative, i.e. a morphism $\me{F}\ra\me{G}$ of $\Db_{\Rc}(k_X,S)$ is an
  isomorphism in $\Db_{\Rc}(k_X,S)$ if and only if it is an isomorphism in
  $\Db(k_X,S)$.
\end{prop}
\begin{proof}
We embed $\me{F}\ra\me{G}$ in a distinguished triangle
$$ \xymatrix{
    {\me{F}} \ar[r] & {\me{G}} \ar[r] & {\me{H}} \ar[r]^+ & } $$
in $\Db_{\Rc}(k_X,S)$. If $\me{F}\ra\me{G}$ is an isomorphism in $\Db(k_X,S)$ then
$\me{H}\simeq 0$ in $\Db(k_X,S)$. Hence $\SS(\me{H})\cap S=\varnothing$. Therefore
$\me{H}\simeq 0$ in $\Db_{\Rc}(k_X,S)$ and $\me{F}\ra\me{G}$ is an isomorphism in
$\Db_{\Rc}(k_X,S)$.
\end{proof}
However there are some situations when the functor \eqref{comparison} is fully
faithful. For instance, Andronikof remarked that the proof of Proposition
\ref{prop51} holds in the constructible case, hence
\begin{prop}
  Let $\me{F},\me{G}\in\on{D}^b_{\text{$\R$-c}}(k_X,p)$. 
  Then there is a canonical isomophism
$$ \on{Hom}_{\on{D}^b_{\text{$\R$-c}}(k_X,p)}(\me{F},\me{G})
   \overset{\sim}{\lra}
   \on{H}^0\!\mu hom(\me{F},\me{G})_p $$
  and the natural functor
  $$  \on{D}^b_{\text{$\R$-c}}(k_X,p)\lra  \on{D}^b(k_X,p) $$
  is fully faithful. 
\end{prop}

\subsection{The category $\Db_{\Rc}(k_X,\{x\}\times\dot\delta)$}

We will see that the natural morphism (\ref{comparison}) is a an isomorphism
for orbits of $\C^{\times}$ and more generally for closed cones $\dot\delta$ in 
$\dot T^*_xX$ for some $x\in X$. 
\begin{prop}\label{prop45}
  Let $\me{F},\me{G}\in\Db_{\Rc}(k_X)$. Then the natural morphism
  $$ \on{Hom}_{\Db_{\Rc}(k_X,\{x\}\times\dot\delta)}(\me{F},\me{G}) \lra
     \on{Hom}_{\Db(k_X,\{x\}\times\dot\delta)}(\me{F},\me{G}) $$
  is an isomorphism.
\end{prop}
\begin{proof}
We may assume that $X$ is a vector space.
Note that if $\me{F}\in\Db_{\Rc}(k_X)$, then $\Phi_{U,\gamma}(\me{F})_U$ is 
$\R$-constructible for any relatively compact subanalytic open subset $U\owns 0$ 
and any subanalytic open cone $\gamma\subset X$.\\ 
Then the results of Section 3.4 hold in the $\R$-constructible case. More precisely,
we see first (as in Proposition \ref{newresult}) that the composition
\begin{align*}  
  \inlim{U,\gamma}{\on{H}^0\Dr\on{Hom}\big(\Phi_{U,\gamma}(\me{F})_U,\me{G}\big)}\lra
    \on{Hom}_{_{\Db_{\Rc}(k_X,\{x\}\times\dot\delta)}}\big(\me{F},\me{G}\big) &\lra
    \on{Hom}_{_{\Db(k_X,\{x\}\times\dot\delta)}}\big(\me{F},\me{G}\big)  \\
    &\lra   \on{H}^0\big(\{x\}\times\dot\delta,\mu hom(\me{F},\me{G})\big)
\end{align*}
is an isomorphism and then (as in Theorem \ref{cutoffprop}) that the composition
$$  \on{Hom}_{_{\Db_{\Rc}(k_X,\{x\}\times\dot\delta)}}\big(\me{F},\me{G}\big) \lra
    \on{Hom}_{_{\Db(k_X,\{x\}\times\dot\delta)}}\big(\me{F},\me{G}\big)  
    \lra   \on{H}^0\!\big(\{x\}\times\dot\delta,\mu hom(\me{F},\me{G})\big) $$
is an isomorphism. Since the second morphism of the last composition is an 
isomorphism, we get the result.
\end{proof}
Combining Proposition \ref{prop45} and Theorem \ref{cutoffprop} we get the following
theorem.
\begin{thm}\label{cor46}
  The natural functor
  $$ \Db_{\Rc}(k_X,\{x\}\times\dot\delta) \lra
     \Db(k_X,\{x\}\times\dot\delta) $$
  is fully faithful. Moreover for every $\me{F},\me{G}\in
  \Db_{\Rc}(k_X,\{x\}\times\dot\delta)$ we have
  $$ \on{Hom}_{\Db_{\Rc}(k_X,\{x\}\times\dot\delta)}(\me{F},\me{G})\overset{\sim}{\lra}
      \on{H}^0\!(\{x\}\times\dot\delta,\mu hom(\me{F},\me{G})). $$
\end{thm}

\subsection{Microlocally $\C$-constructible sheaves}

Recall the microlocal characterization of complexes with $\C$-constructible
cohomology sheaves given in \cite{KS3} (Theorem 8.5.5).
\begin{prop}
   A complex $\me{F}$ in $\Db(k_X)$ has $\C$-constructible cohomology if and only if 
   $\me{F}\in \Db_{\text{$\R$-c}}(k_X)$ and $\on{SS}(\me{F})$ is a 
   $\C^\times$-conic subset of $T^*X$.
\end{prop} 
If $S\subset T^*X$ is a not necessarily $\C^\times$-conic subset, then this suggests 
the definition of a microlocally $\C$-constructible sheaf on $S$ in as follows:
\begin{defi}\label{deficconst}
  \begin{itemize}
     \item[(i)] An object $\me{F}\in\on{D}^b_{\text{$\R$-c}}(k_X)$ (or 
          $\on{D}^b_{\text{$\R$-c}}(k_X,S')$ for $S'\supset S$) is called 
         microlocally 
          $\C$-constructible on $S$ if $\SS(\me{F})$ is $\C^\times$-conic on $S$.
     \item[(ii)] We denote by $\Db_{\Cc}(k_X,S)$ be the full subcategory 
          of $\Db_{\Rc}(k_X,S)$ consisting of microlocally $\C$-constructible 
          sheaves (on $S$).
  \end{itemize}
\end{defi}
\begin{remark}
 \em Note that the category of microlocally $\C$-constructible sheaves (on $S$) 
  is different from the category $\Db_{\Cc}(k_X)/(\mc{N}_{_S}\cap
  \Db_{\Cc}(k_X))$ (i.e. the microlocalization of $\C$-constructible sheaves). There is 
  a natural functor
  $$  \Db_{\Cc}(k_X)/(\mc{N}_{_S}\cap\Db_{\Cc}(k_X)) \lra \Db_{\Cc}(k_X,S), $$
  but in general, an object in $\Db_{\Cc}(k_X,S)$ cannot be represented by a 
  complex with $\C$-constructible cohomology sheaves. One shall keep in mind that 
  by definition an object in $\Db_{\Rc}(k_X,S)$ (the microlocalization of 
  $\R$-constructible sheaves) is represented by an $\R$-constructible sheaf on $X$.
  Since $\R$-constructible sheaves can be defined by a microlocal 
  property\begin{footnote}{An object $\me{F}\in\Db(k_X)$ is weakly $\R$-constructible
  if and only if its micro-support is Lagrangian.}\end{footnote}, there also exists a
  natural definition of a microlocally $\R$-constructible sheaf which we do not 
  consider in this paper.
\end{remark}
\begin{remark}
  \em Of course Definition \ref{deficconst}, (i) is equivalent to the statement
  \begin{itemize}
   \item[(i)] A sheaf $\me{F}\in\Db_{\Rc}(k_X)$ (resp. $\me{F}\in
      \on{D}^b_{\text{$\R$-c}}(k_X,S')$ for $S'\supset S$) is microlocally 
      $\C$-constructible on $S$ if for every point $p$ of $S$ there exists an 
      open neighborhood $U$ of $p$ such that $U\cap \SS(\me{F})=U\cap 
      \C^{\times}\SS(\me{F})$.
  \end{itemize}
  Obviously $\Db_{\Cc}(k_X,T^*X)=\Db_{\Cc}(k_X)$ and if $x\in X$, then 
  $\me{F}\in\Db_{\Rc}(k_X)$ defines an object of $\Db_{\Cc}(k_X,x)$ if and only if 
  $\me{F}|_{_V}$ is $\C$-constrcutible for some neighborhood $V$ of $x$.\\
  However the category $\Db_{\Cc}(k_X,S)$ is not very easy to understand in 
  general, especially if $S$ is not $\C^{\times} $-conic.\em
\end{remark}
\begin{lemma}
   Let $\me{F}\in\Db(k_X)$ and $S\subset \dot T^*X$
  \begin{itemize}
  \item[(i)] The object $\me{F}$ is microlocally $\C$-constructible on $S$ if and only
      if $\me{F}$ is micolocally $\C$-constructible on $\R^+S$.
  \item[(ii)] Suppose that $S\subset\dot T^*X$ is $\C^{\times} $-conic.\\ 
     Then $\me{F}$ is microlocally $\C$-constructible on $S$ if and only
     if $\me{F}$ is micolocally $\C$-constructible in $\gamma^{-1}(U)$ where
     $U$ is a germ of a neighborhood of $\gamma(S)$ in $P^*X$.
  \end{itemize}
\end{lemma}
\begin{proof}
Statemen (i) is a consequence of Lemma \ref{micro1} and (ii) follows from (i) 
and the fact that any $\R^+$-conic neighborhood of $S$ contains a $\C^{\times}$-conic 
neighborhood.
\end{proof}
\begin{remark}
  \em There is an obvious functor $\Db_{\Cc}(k_X,\C^\times S)\ra 
  \Db_{\Cc}(k_X,S)$. One might ask the question whether or not a sheaf $\me{F}$ 
  of $\Db_{\Cc}(k_X,S)$ can be lifted to $\Db_{\Cc}(k_X,\C^\times S)$ and if there
  is a tool to produce an object of $\Db_{\Cc}(k_X,\C^\times S)$ that is 
  isomorphic to $\me{F}$ on $S$. There does not seem to be an obvious answer 
  as the following example shows:\\
  Consider $X=\C^2$ and the sheaf $\me{F}=\C_{\C\times\{(0,0)\}}\oplus 
  \C_{\{0\}\times\R\times\R\times\{0\}}$. Then 
$$\SS(\me{F})=T^*_{\C\times \{(0,0)\}}X\cup T^*_{\{0\}\times\R\times\R\times 
    \{0\}}X.$$
  Take $p=((0,0,0,0);(0,0,1,1))\in\dot T^*X$. Then $\SS(\me{F})=\C^{\times} p$ in 
  a neighborhood of $p$.\\ 
  But if $U\supset \C^{\times} p$ is an arbitrary neighborhood of $\C^{\times} p$, 
  $\SS(\me{F})$ is not $\C^{\times} $-conic on $U$, hence $\me{F}$ is 
  microlocally $\C$-constructible at $p$ but not on $\C^\times p$. However $\me{F}$  
  is isomorphic in $\Db_{\Rc}(k_X,p)$ to the sheaf $\C_{\C\times\{(0,0)\}}$ 
  which is globally $\C$-constructible. The problem is how to construct 
  $\C_{\C\times\{(0,0)\}}$ functorially from $\me{F}$. It cannot be done by a cut-off 
  functor which will always preserve the micro-support in a neighborhood of 
  $\C^\times p$.\\
  Hence microlocally $\C$-constructible sheaves should be defined on $P^*X$ rather then
  $T^*X$ and Definition \ref{deficconst} will mostly be used for a 
  $\C^\times$-conic subset $X$.
\end{remark}

\subsection{The category $\Db_{\Cc}(k_X,\{x\}\times\dot\delta)$}

\begin{prop}\label{prop414}
  The natural functors
  $$ \Db_{\Cc}(k_X,\{x\}\times\dot\delta) \lra
      \Db_{\Rc}(k_X,\{x\}\times\dot\delta) \lra
       \Db(k_X,\{x\}\times\dot\delta) $$
  are fully faithful. Moreover for every $\me{F},\me{G}\in
  \Db_{\Rc}(k_X,\{x\}\times\dot\delta)$ we have
  $$ \on{Hom}_{\Db_{\Cc}(k_X,\{x\}\times\dot\delta)}(\me{F},\me{G})\overset{\sim}{\lra}
      \on{H}^0(\C^\times p,\mu hom(\me{F},\me{G})). $$
\end{prop}
\begin{proof}
The first functor is fully faithful by definition, the second by Theorem \ref{cor46}.
The second part follows again from Theorem \ref{cor46}.  
\end{proof}

\section{Invariance by quantized contact transformations}

Let $\Omega_X\subset T^*X$ be an open subset of a real manifold $X$.
In \cite{KS3}, Kashiwara-Schapira showed that the category $\Db(k_X,\Omega_X)$ (or more
generally the prestack $\Db(k_X,\,*\,)|_{\Omega_X}$) is invariant under ``quantized 
contact transformations''. Let us briefly explain this statement.\\
Consider real manifolds $X,Y$ of the same dimension and open 
subsets $\Omega_X\subset \dot T^*X$, $\Omega_Y\subset\dot T^*Y$. An $\R^+$-homogenous
symplectic isomorphism
$$ \chi:\ \Omega_X \overset{\sim}{\lra} \Omega_Y $$
is often called a contact transformation (although strictly speaking, the contact 
structures are defined on the projective bundles). Invariance under ``quantized 
contact transformations'' means that locally we can construct from $\chi$ an 
equivalence of categories 
$$ \Phi_{\me{K}}:\ \Db(k_X,\Omega_X) \overset{\sim}{\lra} \Db(k_Y,\Omega_Y). $$
The equivalence $\Phi_{\me{K}}$ is explicitly given by an integral transform and 
depends on the choice of a kernel $\me{K}\in\Db(k_{Y\times X})$. The main 
result (\cite{KS3}, Corollary 7.2.2) is:
\begin{thm}\label{kerneltheorem}
   Let $X,Y$ be two real manifolds, $\Omega_X\subset T^*X$, $\Omega_Y\subset T^*Y$
   open subsets and
   $$ \chi:\ \Omega_X\lra \Omega_Y $$ 
   a real contact transformation. Set
$$ \Lambda=\Big\{((y;\eta),(x;\xi))\in\Omega_Y\times\Omega_X^a\ |\ (y,\eta)=
   \chi(x,-\xi)\Big\}. $$
   Let $p_X\in\Omega_X$ and $p_Y=\chi(p_X)$.\\
   There exist open neighborhoods $X'$ of $\pi(p_X)$, $Y'$ of $\pi(p_Y)$, $\Omega_X'$
   of $p_X$, $\Omega_Y'$ of $p_Y$ with $\Omega_X'\subset T^*X'\cap \Omega_X$,
   $\Omega_Y'\subset T^*Y'\cap \Omega_Y$ and a kernel $\me{K}\in
   \Db(k_{Y'\times X'})$ such that: 
   \begin{itemize}
     \item[(1)] $\chi$ induces a contact transformation $\Omega_X'\overset{\sim}{\ra}
       \Omega_Y'$,
     \item[(2)] $$ \Big((\Omega_Y''\times T^*X')\cup (T^*Y'\times {\Omega_X''}^a)\Big)
       \cap\SS(\me{K})\subset \Lambda \cap (\Omega_Y''\times {\Omega_X''}^a)$$ for 
       every open subsets $\Omega_{X}''\subset\Omega_{X}'$ and $\Omega_{Y}''=
        \chi(\Omega_{X}'')$,
     \item[(3)] composition with $\me{K}$ induces an equivalence of prestacks
       $$\Phi_{\me{K}}=\me{K}\circ :\, \chi_*\Db(k_X',\,*\,)|_{\Omega_{X}'} \lra 
             \Db(k_{Y'},\,*\,)|_{\Omega_{Y}'},$$
       a quasi-inverse being given by $\Phi_{\me{K}^*}$ whith
       $\me{K}^*=r_*\Dr \me{H}om(\me{K},\omega_{Y\times X|X})$ where $r:\,Y\times X\ra
       X\times Y$ switches the factors.
     \item[(4)] $\SS(\Phi_{\me{K}}(\me{F}))\cap \Omega_{Y}''=
                   \chi(\SS(\me{F})\cap \Omega_{X}''),$
     \item[(5)] $\chi_*\mu hom(\me{F},\me{G})|_{\Omega_{X'}}\simeq 
                \mu hom(\Phi_{\me{K}}(\me{F}),\Phi_{\me{K}}(\me{G}))|_{\Omega_{Y'}}.$
   \end{itemize}
\end{thm}
\begin{remark}
  \em There is a slight difference to \cite{KS3}, Corollary 7.2.2 where the smaller
  open subset $\Omega_X''$ is not introduced. For instance, statement (3) is 
  stated as an equivalence of categories 
 $$\Phi_{\me{K}}:\  \Db(k_{X'},\Omega_{X}') \lra \Db(k_{Y'},\Omega_{Y}').$$
  However, it is obvious that if we have the theorem for fixed $X',Y',\Omega_X',
  \Omega_Y'$, then (1), (2), (4) and (5) are valid for any open subset 
  $\Omega_X''\subset\Omega_X'$ and $\Omega_Y''=\chi(\Omega_X'')$. The fact that (3) 
  is valid on $\Omega_X''$ is based on the observation that (3) is essentially a 
  consequence of (2) (plus some independent conditions on $\me{K}$, cf. \cite{KS3}, 
  Theorem 7.2.1).
\end{remark}
Now let us consider constructible sheaves. It is not immediately obvious that the 
equivalence (3) of Theorem \ref{homogenouskerneltheorem} should induce an 
equivalence on the microlization of $\R$-constructible (resp. on microlocally 
$\C$-constructible or later on microlocally perverse) sheaves 
$$\chi_*\Db_{\Rc}(k_X,\,*\,)|_{\Omega_X} \overset{?}{\lra} 
   \Db_{\Rc}(k_Y,\,*\,)|_{\Omega_Y} $$
since the functor $\Phi_{\me{K}}$ is not well defined on $\R$-constructible sheaves.
This problem can be solved at a point $p\in T^*X$ by using the microlocal composition 
of \cite{KS3}. In \cite{An2}, Andronikof uses this tool to construct the functor 
$\Phi^\mu_{\me{K}}$. Then one can treat a variety of kernels $\me{K}$ but a priori,
one can no longer work in an open neighborhood.\\
However, under the hypothesis of Theorem \ref{kerneltheorem}, we do not need 
to use microlocal composition of kernels. We can always define the functor
$$ \chi_*\Db_{\Rc}(k_X,\,*\,)|_{\Omega_X} \lra  \chi_*\Db(k_X,\,*\,)|_{\Omega_X} 
   \overset{\Phi_{\me{K}}}\lra \Db(k_Y,\,*\,)|_{\Omega_Y} $$
and hope that it factors through $\Db_{\Rc}(k_Y,\,*\,)|_{\Omega_Y}$. However, one 
has to beware that  $\Db_{\Rc}(k_Y,\,*\,)|_{\Omega_Y}$ is not a full subprestack of 
$\Db(k_Y,\,*\,)|_{\Omega_Y}$, hence it is not sufficent to show that 
$\Phi_{\me{K}}(\me{F})$ is isomorphic to an $\R$-constructible sheaf in 
$\Db(k_Y,\,*\,)|_{\Omega_Y}$. We encounter this problem in Section 5.1 (see Theorem 
\ref{kerneltheoremrcons}).

\subsection{Quantization of $\R$-constructible sheaves}

Recall the definition of the full subcategory $\on{N}(Y,X,\Omega_Y,\Omega_X)$ of
$\Db(k_{Y\times X},\Omega_Y\times T^*X)$. Its objects are kernels $\me{K}$ on 
$Y\times X$ such that
\begin{itemize}
  \item[(i)] $\SS(\me{K})\cap (\Omega_Y\times T^*X)\subset \Omega_Y\times \Omega_X^a$
  \item[(ii)] The projection $p_1:\,\SS(\me{K})\cap (\Omega_Y\times T^*X)\ra \Omega_Y$
       is proper.
\end{itemize}
If $V=\pi_X(\Omega_X)$ is a subanalytic relatively compact open subset of $X$ we set
 $$ \on{N}_{\Rc}(Y,X,\Omega_Y,\Omega_X)=\on{N}(Y,X,\Omega_Y,\Omega_X)\cap 
\Db_{\Rc}(k_{Y\times X},\Omega_Y\times T^*X).$$ 
\begin{defi}
  Let $\me{K}\in \on{N}_{\Rc}(Y,X,\Omega_Y,\Omega_X)$. We define the functor 
  $\Phi^{\Rc}_{\me{K}}$ as 
 $$ \Db_{\Rc}(k_X,\Omega_X) \lra \Db_{\Rc}(k_Y,\Omega_Y) \quad ; \quad
     \me{F}\mapsto \me{K}\circ \me{F}_{\ol{V}}. $$
\end{defi}
Note that $\me{K}\circ\me{F}_{\ol{V}}$ 
is $\R$-constructible since $V$ is subanalytic and relatively compact. We may 
visualize the situation by the following diagram
$$ \xymatrix{
       \Db_{\Rc}(k_X,\Omega_X)  \ar[r]^{\Phi^{\Rc}_{\me{K}}} \ar[d] &   
       \Db_{\Rc}(k_Y,\Omega_Y) \ar[d] \\
       \Db(k_X,\Omega_X)   \ar[r]_{\me{K}\circ}  &   \Db(k_Y,\Omega_Y). } $$
The square is commutative up to natural isomorphism since the natural morphism
$\me{F}\ra\me{F}_{\ol{V}}$ induces an isomorphism $\me{K}\circ\me{F}\overset{\sim}{\ra}
\me{K}\circ\me{F}_{\ol{V}}$ in $\Db(k_Y,\Omega_Y)$.\\
Now suppose that $\me{K}\in\on{N}_{\Rc}(Z,Y,\Omega_Z,\Omega_Y)$ and 
$\me{L}\in\on{N}_{\Rc}(Y,X,\Omega_Y,\Omega_X)$. Then their composition 
$\me{K}\circ\me{L}$ is well defined in $\on{N}(Z,X,\Omega_Z,\Omega_X)$ but not 
necessarily in $\on{N}_{\Rc}(Z,X,\Omega_Z,\Omega_X)$. Note however that by definition 
$\me{L}\in\Db(k_{Y\times X},\Omega_Y\times T^*X)$. Hence, if we set 
$W=\pi_Y(\Omega_Y)$, we do not distinguish between $\me{L}$ and 
$\me{L}_{\ol{W}\times X}$ in $\me{L}\in\on{N}(Y,X,\Omega_Y,\Omega_X)$. 
In other words we get a natural isomorphism
$$ \me{K}\circ\me{L} \overset{\sim}{\lra} \me{K}\circ\me{L}_{\ol{W}\times X} $$
in $\Db(k_{Z\times X},\Omega_Z\times T^*X)$. Then we get natural isomorphisms
$$ \Phi_{\me{K}}^{\Rc}\Phi_{\me{L}}^{\Rc}(\me{F})\simeq
     \me{K}\circ (\me{L}\circ \me{F}_{\ol{V}})_{\ol{W}} \simeq
     (\me{K} \circ \me{L}_{\ol{W}\times X}) \circ \me{F}_{\ol{V}}\simeq\Phi_{\me{K}
      \circ \me{L}_{\ol{W}\times X}}^{\Rc}(\me{F}). $$
Hence, the theory of microlocal kernels (Section 7.1. of \cite{KS3})
works well in the $\R$-constructible case if we restrict ourselves to relatively 
compact subanalytic open sets.  

Finally suppose that $\Omega_X'\subset \Omega_X$, $\Omega_Y'\subset \Omega_Y$ and
that $\me{K}\in\on{N}_{\Rc}(Y,X,\Omega_Y,\Omega_X)\cap
\on{N}_{\Rc}(Y,X,\Omega_Y',\Omega_X')$. Then we get a diagram
$$ \xymatrix{
   {\Db_{\Rc}(k_X,\Omega_{X})} \ar[r]^{\Phi_{\me{K}}^{\Rc}} \ar[d] & 
   {\Db_{\Rc}(k_Y,\Omega_Y)} \ar[d]  \\
   {\Db_{\Rc}(k_X,\Omega_{X}')} \ar[r]_{\Phi_{\me{K}}^{\Rc}}  & 
   {\Db_{\Rc}(k_Y,\Omega_Y')} } $$
that is commutative up to natural isomorphism induced by $(\me{F}_{\ol{V}})_{\ol{V'}}
\simeq \me{F}_{\ol{V'}}$ where $V'=\pi(\Omega_Y')$. 

Now suppose that we are given a contact transformation
$$ \chi:\ \Omega_X\overset{\sim}{\lra} \Omega_Y $$
where we assume that $\pi(\Omega_X),\pi(\Omega_Y)$ are relatively compact subanalytic 
open sets. If there exists an object $\me{K}$ such that 
$\me{K}\in\on{N}_{\Rc}(Y,X,\Omega_Y',\Omega_X')$ for all open subsets 
$\Omega_X'\subset \Omega_X$ and $\Omega_Y'=\chi(\Omega_X')$ then 
we get a commutative diagram of functors of prestacks
$$ \xymatrix{
      {\Db_{\Rc}(k_X,\,*\,)|_{\Omega_X}} \ar[r]^{\Phi^{\Rc}_{\me{K}}} \ar[d]  & 
      {\Db_{\Rc}(k_Y,\,*\,)|_{\Omega_Y}} \ar[d] \\
      {\Db(k_X,\,*\,)|_{\Omega_X}} \ar[r]_{\Phi_{\me{K}}} & 
      {\Db(k_Y,\,*\,)|_{\Omega_Y}} }$$
All compatibility conditions are easily verified by diagram chases.

We are now ready to quantize $\R$-constructible sheaves. All we need to know is that
the kernel produced in Theorem \ref{kerneltheorem} can be taken $\R$-constructible.
\begin{thm}\label{kerneltheoremrcons}
   Let $X,Y$ be two real manifolds, $\Omega_X\subset T^*X$, $\Omega_Y\subset T^*Y$
   open subsets and
   $$ \chi:\ \Omega_X\overset{\sim}{\lra} \Omega_Y $$ 
   a contact transformation. Set
$$ \Lambda=\Big\{((y;\eta),(x;\xi))\in\Omega_Y\times\Omega_X^a\ |\ (y,\eta)=
   \chi(x,-\xi)\Big\}. $$
   Let $p_X\in\Omega_X$ and $p_Y=\chi(p_X)$.\\
   There exist open neighborhoods $X'$ of $\pi(p_X)$, $Y'$ of $\pi(p_Y)$, $\Omega_X'$
   of $p_X$, $\Omega_Y'$ of $p_Y$ with $\Omega_X'\subset T^*X'\cap \Omega_X$,
   $\Omega_Y'\subset T^*Y'\cap \Omega_Y$, and there exists $\me{K}\in
   \Db_{\Rc}(k_{Y'\times X'})$ such 
   \begin{itemize}
     \item[(1)] $\chi$ induces a contact transformation $\Omega_X'
        \overset{\sim}{\ra} \Omega_Y'$,
     \item[(2)] $$\Big((\Omega_Y''\times T^*X')\cup (T^*Y'\times {\Omega_X''}^a)\Big)
      \cap\SS(\me{K})\subset \Lambda \cap (\Omega_Y''\times {\Omega_X''}^a)$$ for every
        open subsets $\Omega_{X}''\subset\Omega_{X}'$ and $\Omega_{Y}''=
        \chi(\Omega_{X}'')$,
     \item[(3)] composition with $\me{K}$ induces an equivalence of prestacks
       $$\Phi_{\me{K}}^{\Rc}:\ \chi_*\Db_{\Rc}(k_{X'},\,*\,)|_{\Omega_{X}'} \lra 
             \Db_{\Rc}(k_{Y'},\,*\,)|_{\Omega_{Y}'},$$
       and $\Phi_{\me{K}^*}^{\Rc}$ defines a quasi-inverse where 
       $\me{K}^*=r_*\Dr \me{H}om(\me{K},\omega_{Y\times X|X})$ and $r:\,Y\times X\ra
       X\times Y$ switches the factors.
     \item[(4)] $\SS(\Phi_{\me{K}}^{\Rc}(\me{F}))\cap \Omega_{Y}''=
                   \chi(\SS(\me{F})\cap \Omega_{X}''),$
     \item[(5)] $\chi_*\mu hom(\me{F},\me{G})|_{\Omega_{X}'}\simeq 
                \mu hom(\Phi_{\me{K}}(\me{F}),\Phi_{\me{K}}(\me{G}))|_{\Omega_{Y}'}.$
   \end{itemize}
\end{thm}
\begin{proof}
Let us show that in the situation of Theorem \ref{kerneltheorem} we can choose the 
kernel to be $\me{K}$ $\R$-constructible.\\
First assume that $\Lambda$ is the conormal bundle of a smooth hypersurface 
$S\subset Y\times X$. Then one can take the kernel $\me{K}=k_S$ (cf. \cite{KS3}, 
Corollary 7.2.2). Recall that locally $\chi$ may be decomposed as 
$$\chi=\chi_1\circ\chi_2:\ \Omega_X\lra\Omega_Z\lra\Omega_Y$$ 
where the Lagrangian manifold $\Lambda_i$ associated to the contact 
transformation $\chi_i$ is the conormal bundle to a smooth hypersurface $S_i$.
By shrinking $\Omega_X$ and $\Omega_Y$ we may assume that $\pi_Z(\Omega_Z)$ is 
subanalytic and relatively compact. Then $k_{S_1}\circ 
k_{S_2\cap (\pi_Z(\Omega_Z)\times X)}$ is $\R$-constructible and satifies (1), (2),
(4) and (5) of the theorem (cf. \cite{KS3}, Corollary 7.2.2).\\
Moreover we know that $\Phi^{\Rc}_{\me{K}}$ and $\Phi^{\Rc}_{\me{K}^*}$ are 
well-defined and that they are quasi-inverse functors in the non-constructible case.\\
By definition, we have
$$\Phi^{\Rc}_{\me{K}}\circ\Phi^{\Rc}_{\me{K}^*}\simeq \Phi^{\Rc}_{\me{K}\circ
  \me{K}^*_{\ol{V}\times Y}}. $$
Recall that there is a natural isomorphism
$$ k_{\Delta_X} \lra  \me{K}\circ \me{K}^* $$
in $\Db(k_X,\Omega_X)$. Hence we get an isomorphism
\begin{equation} \label{equivb}
   k_{\Delta_X} \lra  \me{K}\circ \me{K}^* \lra \me{K}\circ \me{K}^*_{\ol{V}\times Y} 
\end{equation}
in $\Db(k_X,\Omega_X)$. It is sufficent to prove that this morphism is well defined in 
$\Db_{\Rc}(k_X,\Omega_X)$.\\
Denote by $q_{12},q_{13},q_{23}$ the obvious projections from $X\times X\times Y$. 
We get a commutative diagram
$$ 
 \xymatrix{
   {k_{\Delta_X}} \ar[r] \ar[dr] & 
    {\Dr q_{12*}\Dr\me{H}om(q_{13}^{-1}\me{K},q_{23}^!\me{K})} 
         \ar[d] & {\me{K}\circ\me{K}^*} \ar[l] \ar[d] & \\
    &  {\Dr q_{12*}\Dr\me{H}om(q_{13}^{-1}\me{K},q_{23}^!\me{K}_{\ol{V}\times Y})} &
           {\me{K}\circ\me{K}^*_{\ol{V}\times Y}} \ar[l]  & } $$
that is defined in $\Db(k_X)$. Note that the lower part is well-defined in 
$\Db_{\Rc}(k_X)$. All morphisms become isomorphisms in $\Db(k_X,\Omega_X)$
(cf. \cite{KS3}, Theorem 7.2.1). Since the natural functor $\Db_{\Rc}(k_X,\Omega_X)\ra
\Db(k_X,\Omega_X)$ is conservative, this shows that \eqref{equivb} is well-defined in 
$\Db_{\Rc}(k_X,\Omega_X)$.\\
Similarly one shows that the kernel $r^{-1}\Dr\me{H}om(\me{K},\omega_{X\times Y|X})$ 
defines a right inverse of $\me{K}$ which proves that $\Phi^{\Rc}_{\me{K}}$ is an 
equivalence. Then $\Phi^{\Rc}_{\me{K}^*}$ is actually a quasi-inverse since it is
a left inverse of an equivalence.
\end{proof}

\subsection{Quantization of $\C$-constructible sheaves}

It is now easy to transfer the results of the last section to  microlocally 
$\C$-constructible sheaves.\\
Consider complex manifolds $X,Y$ of the same dimension and open $\C^\times$-conic
subsets $\Omega_X\subset \dot T^*X$, $\Omega_Y\subset\dot T^*Y$. We will call
a $\C^\times$-homogenous symplectic isomorphism
$$ \chi:\ \Omega_X \overset{\sim}{\lra} \Omega_Y $$
a contact transformation (omitting ``complex'' since we will never consider real
contact transformations when dealing with microlocally $\C$-constructible sheaves). 
Then we get the analogous statements of Theorems \ref{kerneltheorem} and 
\ref{kerneltheoremrcons} by replacing open sets with $\C^\times$-conic open sets. 
Let us give the precise formulation only in the case of microlocally 
$\C$-constructible sheaves.  
\begin{thm}\label{homogenouskerneltheorem}
   Let $X,Y$ be two complex manifolds, $\Omega_X\subset \dot T^*X$, $\Omega_Y\subset 
   \dot T^*Y$ open subsets and
   $$ \chi:\ \Omega_X\lra \Omega_Y $$ 
   a homogenous complex contact transformation. Set
$$ \Lambda=\Big\{((y;\eta),(x;\xi))\in\Omega_Y\times\Omega_X^a\ |\ (y,\eta)=
  \chi(x,-\xi)\Big\}. $$
   Then $\Lambda$ is $\C^\times$-conic. Let $p_X\in\Omega_X$ and $p_Y=\chi(p_X)$.\\
   There exist open neighborhoods $X'$ of $\pi(p_X)$, $Y'$ of $\pi(p_Y)$,
   $\C^\times$-conic open neighborhoods $\Omega_X'$ of $\C^\times p_X$, $\Omega_Y'$ of 
   $\C^\times p_Y$ with $\Omega_X'\subset T^*X'\cap \Omega_X$, $\Omega_Y'\subset T^*Y'
   \cap \Omega_Y$ and a kernel $\me{K}\in\Db_{\Rc}(k_{Y'\times X'})$ satisfying 
   \begin{itemize}
     \item[(1)] $\chi$ induces a homogenous contact transformation 
         $\chi:\,\Omega_X'\overset{\sim}{\ra} \Omega_Y'$,
     \item[(2)] $$\Big((\Omega_Y''\times T^*X')\cup (T^*Y'\times {\Omega_X''}^a)\Big)
   \cap\SS(\me{K})\subset \Lambda \cap (\Omega_Y''\times {\Omega_X''}^a)$$ for every
        $\C^\times$-conic open subset $\Omega_{X}''\subset\Omega_{X}'$ and 
        $\Omega_{Y}''=\chi(\Omega_{X}'')$,
     \item[(3)] the functor $\Phi_{\me{K}}^{\Rc}$ induces an equivalence of prestacks
       $$\Phi_{\me{K}}^{\Cc}:\ \chi_*\Db_{\Cc}(k_{X'},\,*\,)|_{\Omega_{X}'} \lra 
             \Db_{\Cc}(k_{Y'},\,*\,)|_{\Omega_{Y}'}.$$
     \item[(4)] $\SS(\Phi_{\me{K}}^{\Cc}(\me{F}))\cap \Omega_{Y}''=
                   \chi(\SS(\me{F})\cap \Omega_{X}'')$
     \item[(5)] $\chi_*\mu hom(\me{F},\me{G})|_{\Omega_X'}\simeq \
        \mu hom(\Phi^{\Cc}_{\me{K}}(\me{F}),
        \Phi_{\me{K}}^{\Cc}(\me{G}))|_{\Omega_{Y'}}$.
   \end{itemize}
\end{thm}
\begin{proof}
Since by definition $\Db_{\Cc}(k_X,\Omega_X'')$ is a full subcategory of 
$\Db_{\Rc}(k_X,\Omega_X'')$, it is enough to show that for any 
$\me{F}\in\Db_{\Cc}(k_X,\Omega_X)$ the object $\Phi^{\Rc}_{\me{K}}(\me{F})$ is 
an object of $\Db_{\Cc}(k_Y,\Omega_Y)$. Hence we have to show that 
$\SS(\Phi(\me{F}))$ is $\C^\times$-conic on $\Omega_Y$. Since $\chi$ is 
$\C^\times$-homogenous this is easily verified by the formula
$$ \SS(\Phi_{\me{K}}^{\Rc}(\me{F}))\cap\Omega_X''=\chi(\SS(\me{F})\cap\Omega_Y''). $$
\end{proof}

\section{Microlocally complex constructible sheaves on $\C^\times p$}

Let $p\in\dot T^*X$. As a special case of Proposition \ref{prop414} we get that
the natural functor 
$$ \Db_{\Cc}(k_X,\C^{\times} p)\lra \Db(k_X,\C^{\times} p)$$ 
is fully faithful. Moreover, morphisms in $\Db(k_X,\C^{\times} p)$ between microlocally 
$\C$-constructible sheaves $\me{F},\me{G}\in\Db_{\Cc}(k_X,\C^\times p)$ are given 
by sections of $\mu hom(\me{F},\me{G})$ on $\C^{\times} p$.\\
In this section we will give a description of the objects of 
$\Db_{\Cc}(k_X,\C^\times p)$ using quantized contact transformation and the generic 
position theorem. More precisely we will show in Section 6.1 that if $\me{F}$ is 
microlocally $\C$-constructible on $\C^\times p$ and $\SS(\me{F})$ is in generic 
position (i.e. $\SS(\me{F})\cap \C^\times p$ is isolated
in $\pi^{-1}\pi(p)$) then $\me{F}$ is isomorphic in $\Db(k_X,\C^\times p)$ to an object 
of $\Db_{\Cc}(k_X)$. 

It will often be convenient to fix an $\R^+$-conic Lagrangian variety $\Lambda$ in 
$T^*X$ and to consider only sheaves whose micro-support is contained in $\Lambda$. We 
introduce the following categories:
\begin{defi}
  Let $\Lambda\subset T^*X$ be an $\R^+$-conic Lagrangian variety defined in a 
  neighborhood of a subset $S\subset T^*X$. Then we define the following two 
  categories:
  \begin{itemize}
     \item[(1)] $\Db_{\Rc,\Lambda}(k_X,S)\subset \Db_{\Rc}(k_X,S)$ is the
        full subcategory of $\Db_{\Rc}(k_X,S)$ whose objects $\me{F}$ satisfy 
        $\SS(\me{F})\subset\Lambda$ in a neighborhood of $S$.
     \item[(2)] $\Db_{\Cc,\Lambda}(k_X,S)=\Db_{\Cc}(k_X,S)\cap
         \Db_{\Rc,\Lambda}(k_X,S)$.
  \end{itemize}
\end{defi}
Of course, the second definition is only of interest when $S$ is $\C^\times$-conic
and $\Lambda$ is $\C^\times$-conic on the subset $S$.

\subsection{Microlocal complex constructible sheaves in generic position}
\begin{defi}
  Let $p\in \dot T^*X$ and $\Lambda\subset T^*X$ be an $\R^+$-conic Lagrangian
  subset such that $\Lambda$ is $\C^{\times}$-conic in a neighborhood of
  $\C^{\times} p$. We say that $\Lambda$ is in generic position at $p$ if 
  $$ \Lambda \cap \pi^{-1}\big(\pi(p)\big)\subset \C^{\times} p $$
  in a neighborhood of $\C^{\times} p$.
\end{defi}
\begin{remark}
  \em If $\Lambda$ is in generic position at $p$, then $\Lambda$ is  
  $\C^\times$-conic in a neighborhood of $\C^\times p$ by definition. Hence either
  $\Lambda \cap \pi^{-1}\big(\pi(p)\big)\subset \C^{\times} p$ or 
  $\Lambda \cap \pi^{-1}\big(\pi(p)\big)=\varnothing$. Moreover, being in generic 
  position is an open property. Also note that if $\Lambda$ 
  is in generic position then it is in generic position in the sense of \cite{KK}, 
  Chapter I.6 where $\Lambda$ is supposed to be only locally $\C^\times$-conic 
  in a neighborhood of $p$. 
\end{remark}
\begin{prop}\label{lambdafix}
  Let $\Lambda\subset T^*X$ be an $\R^+$-conic Lagrangian variety 
  that is $\C^\times$-conic in a neighborhood of $\C^{\times} p$. Suppose that 
  $\Lambda$ is in generic position at $p$. Then there exists a fundamental system of 
  conic open subanalytic neighborhoods $\gamma$ of $\C^{\times} p$ in 
  $\dot T^*_{\pi(p)}X$ and for each $\gamma$ a fundamental system of open relatively 
  compact subanalytic neighborhoods $U$ of $\pi(p)$ such that  
  \begin{itemize}
    \item[(1)] the microlocal cut-off functor $\Phi_{U,\gamma}:\,\Db(k_X)\ra\Db(k_X)$
       induces a functor
      $$\Phi_{U,\gamma}:\, \Db_{\Cc,\Lambda}(k_X,U\times\gamma)\lra
        \Db_{\Cc}(k_X,\C^\times p), $$
        and this functor factors as 
       $$\Phi_{U,\gamma}:\, \Db_{\Cc,\Lambda}(k_X,U\times\gamma)
          \lra \Db_{\Cc}(k_X,\pi(p))\lra \Db_{\Cc}(k_X,\C^\times p). $$
    \item[(2)] There is an isomorphism of functors $\Phi_{U,\gamma}
       \overset{\sim}{\lra} \iota$ where 
       $\iota:\,\Db_{\Cc,\Lambda}(k_X,U\times\gamma)
       \ra \Db_{\Cc}(k_X,\C^{\times} p)$ is the natural functor.
    \item[(3)] $\SS\big(\Phi_{U,\gamma}(\me{F})\big) \cap \dot\pi^{-1}(V)=
       \SS(\me{F})\cap (V\times\gamma)$ for sufficently small open 
       neighborhoods $V$ of $\pi(p)$.
    \item[(4)] $\SS\big(\Phi_{U,\gamma}(\me{F})\big)\cap
       \dot\pi^{-1}\pi(p)\subset \C^{\times} p$.
  \end{itemize}
\end{prop}
\begin{proof}
Since the functor $\Phi_{U,\gamma}$ sends $\Db_{\Rc}(k_X)$ to $\Db_{\Rc}(k_X)$ 
(cf. Proposition \ref{prop45}), it induces 
$$ \Phi_{U,\gamma}:\, \Db_{\Rc}(k_X,U\times \gamma)\lra \Db_{\Rc}(k_X,U\times\gamma). $$
Since $\SS(\Phi_{U,\gamma}(\me{F}))\cap U\times \gamma=\SS(\me{F})\cap U\times \gamma$
the functor $\Phi_{U,\gamma}$ preserves microlocally $\C$-constructible sheaves and
we get the functor of (1) as
$$ \Phi_{U,\gamma}:\, \Db_{\Cc}(k_X,U\times \gamma)\lra \Db_{\Cc}(k_X,U\times\gamma)
   \lra \Db_{\Cc}(k_X,\C^\times p). $$
Recall that there exist a fundamental system of conic subanalytic open neighborhoods 
$\gamma$ of $\C^{\times} p$ and for each $\gamma$  a fundamental system of relatively
 compact subanalytic open neighborhoods $U$ of $\pi(p)$ such that $(U,\gamma)$ is a 
refined cutting pair. If $\gamma$ is sufficently small then $\gamma\cap 
\pi^{-1}\pi(p)\cap\Lambda=\C^{\times} p$. Next we choose $U$ sufficently small such 
that $U\cap \partial\gamma\cap \Lambda=\varnothing$. Then (2) follows from the 
refined micolocal cut-off lemma, (3) from Corollary \ref{kleineergaenzung} and (4) 
from (3).\\
Finally let us prove the factorization of (1). For this purpose let us first write 
$\Db_{\Cc,\Lambda}(k_X,U\times\gamma)$ as a localization of a full subcategory of 
$\Db_{\Rc}(k_X)$.\\ 
Denote by $\Db_{\Rc,\Lambda,U\times\gamma}(k_X)$ the full subcategory 
of $\Db_{\Rc}(k_X)$ such that 
$$\SS(\me{F})\cap U\times\gamma\subset \Lambda. $$ 
The category $\Db_{\Rc,\Lambda,U\times\gamma}(k_X)$ is obviously a full triangulated 
subcategory. Now we localize $\Db_{\Rc,\Lambda,U\times\gamma}(k_X)$ by 
complexes whose micro-support is disjoint from $U\times\gamma$ (hence by 
$\mc{N}_{\Rc,U\times\gamma}$). Then we get a natural functor
\begin{equation} \label{dafleche} 
   \Db_{\Rc,\Lambda,U\times\gamma}(k_X)/\mc{N}_{\Rc,U\times\gamma} \lra
   \Db_{\Rc,\Lambda}(k_X,U\times\gamma).
\end{equation}
If $\me{F}\in\Db_{\Rc,\Lambda,U\times\gamma}(k_X)$ and  
 $$ (\SS(\me{F})\cap U\times\gamma)\subset (\Lambda\cap U\times\gamma),$$
then any object $\me{F}'$ that is isomorphic to $\me{F}$ on $U\times\gamma$ is also 
an object of $\Db_{\Rc,\Lambda,U\times\gamma}(k_X)$. Hence we get 
that \eqref{dafleche} is an equivalence.\\
By assumption $\Lambda$ is $\C^{\times}$-conic in a neighborhood of $\C^\times p$. 
Hence we may assume that $\Lambda$ is $\C^\times$-conic on $U\times\gamma$. One can 
show that if $\SS(\me{F})\subset\Lambda$ then $\SS(\me{F})$ is $\C^\times$-conic on 
$U\times\gamma$ (cf. Theorem 8.5.5 of \cite{KS3}). Let us recall the idea of the proof.
First one shows that $\SS(\me{F})$ is open in $\Lambda$ (on $U\times\gamma$) and 
therefore locally $\C^\times$-conic, i.e. for every $\C^\times$-orbit $S$ the set
$\SS(\me{F})\cap S\cap U\times\gamma$ is open in $S\cap U\times\gamma$ (Lemma 8.3.14
of \cite{KS3}). Then, by Proposition 8.5.2 of \cite{KS3}, one gets that
$\Lambda$ is $\C$-analytic on $U\times\gamma$. Hence $\Lambda$ is $\C$-analytic and
$\R^+$-conic and therefore $\C^\times$-conic.\\ 
Thus, we get the equivalence 
$$  \Db_{\Rc,\Lambda,U\times\gamma}(k_X)/\mc{N}_{\Rc,U\times\gamma}\simeq
   \Db_{\Rc,\Lambda}(k_X,U\times\gamma)\simeq \Db_{\Cc,\Lambda}(k_X,U\times\gamma).$$
By (3) and the assumption that $\Lambda$ is $\C^{\times}$-conic in a neighborhood of
$\C^{\times} p$ we get the functor
$$ \Db_{\Rc,\Lambda,U\times\gamma}(k_X)\overset{\Phi_{U,\gamma}}{\lra} 
   \Db_{\Cc}(k_X,\pi(p)) $$
for sufficently small $(U,\gamma)$. 
If $\me{F}\in \Db_{\Rc,\Lambda,U\times\gamma}(k_X)\cap
\mc{N}_{\Rc,U\times\gamma}$ then again by (3) $\me{F}$ is zero in  
$\Db_{\Cc}(k_X,\pi(p))$. Hence we get the factorization of (1).
\end{proof}
\begin{lemma}\label{compatibilitycondition}
  In the situation of Proposition \ref{lambdafix}, suppose that
  we have two pairs $(V,\delta)\subset (U,\gamma)$ such that (1),(2),(3),(4) are 
  satisfied. Then the natural morphism 
  $$ \Phi_{V,\delta}(\me{F}) \lra \Phi_{U,\gamma}(\me{F}) $$
  is an isomorphism in $\Db_{\Cc}(k_X,\pi(p))$.   
\end{lemma}
\begin{proof}
Embed the morphism in a distinguished triangle
$$ \xymatrix{
    {\Phi_{V,\delta}(\me{F})} \ar[r] &  {\Phi_{U,\gamma}(\me{F})} 
   \ar[r] & {\me{H}} \ar[r] &}$$  
in $\Db_{\Rc}(k_X)$. We want to show that $\me{H}$ is isomorphic to zero in a 
neighborhood of $\pi(p)$.\\
First note that we have $\SS(\me{H})\cap V\times\delta=\varnothing$ because 
$\Phi_{V,\delta}(\me{F}) \lra \Phi_{U,\gamma}(\me{F})$ is an isomorphism on 
$V\times \delta$. By (3), there exists an open neighborhood $V'\subset V$ of $\pi(p)$ 
such that
$$ \SS(\Phi_{V,\delta}(\me{F}))\cap \dot\pi^{-1}(V')=\SS(\me{F})\cap (V'\times\delta) $$
and 
$$ \SS(\Phi_{U,\gamma}(\me{F}))\cap \dot\pi^{-1}(V')=\SS(\me{F})\cap (V'\times\gamma). $$
However if we take $V'$ sufficently small then we will show that  
$$ \SS(\Phi_{U,\gamma}(\me{F}))\cap \dot\pi^{-1}(V')=\SS(\me{F})\cap (V'\times\delta). $$
Indeed suppose that this is not possible. Then we can construct a sequence 
$(x_n,\xi_n)\in \dot T^*X$ such that $x_n\ra \pi(p)$ and $\xi_n\in\gamma\setminus\delta$
and $|\xi_n|=1$. By extracting a convergent subsequence, we get a point
$(\pi(p),\xi)\in\SS(\Phi_{U,\gamma}(\me{F}))\cap \{\pi(p)\}\times 
(\gamma\setminus\delta)$. But this is impossible because by (4) we have
$\SS(\Phi_{U,\gamma}(\me{F}))\cap \pi^{-1}\pi(p)=\C^\times p$.\\
Hence 
$$ \SS(\me{H}) \cap \dot\pi^{-1}(V')\subset 
   (\SS(\Phi_{V,\delta}(\me{F})\cup \SS(\Phi_{U,\gamma}(\me{F}))\cap
   \dot\pi^{-1}\pi(V') \subset V'\times\delta $$
and therefore  $\SS(\me{H}) \cap \pi^{-1}(V')\subset T^*_XX$. 
But if $V'$ is sufficently small then
$$ \SS(\me{H}) \cap  V'  \subset 
   (\SS(\Phi_{V,\delta}(\me{F})\cup \SS(\Phi_{U,\gamma}(\me{F}))\cap
     V' \subset \pi(\Lambda)\cap V'  $$
where $\pi(\Lambda)$ is a closed hypersurface and we can conclude that
$\me{H}\simeq 0$ by involutivity of the micro-support.
\end{proof}
\begin{thm}\label{dafunctorphi}
  Suppose that $\Lambda$ is in generic position at $p$. Then
  the microlocal cut-off functors induce a fully faithful functor
  $$ \Phi:\ \Db_{\Cc,\Lambda}(k_X,\C^{\times} p) \lra \Db_{\Cc}(k_X,\pi(p)) $$
\end{thm}
\begin{proof}
 We obviously have the equivalence
  $$ \twodilim{\C^{\times} p\in U\times\gamma}
      \Db_{\Cc,\Lambda}(k_X,U\times\gamma) 
      \overset{\sim}{\lra}\Db_{\Cc,\Lambda}(k_X,\C^{\times} p). $$
We may assume by cofinality that $(U,\gamma)$ is a sufficently small refined cutting 
pair such that Proposition \ref{lambdafix} holds. By Lemma \ref{compatibilitycondition}
the functors $\Phi_{U,\gamma}$ of Proposition \ref{lambdafix} induce a functor
$$ \Phi:\ \Db_{\Cc,\Lambda}(k_X,\C^{\times} p) \lra \Db_{\Cc}(k_X,\pi(p)). $$
Let us show that $\Phi$ is fully faithful.\\
Let $\me{F},\me{G}\in\Db_{\Cc,\Lambda}(k_X,\C^{\times} p)$. Note that by 
Proposition \ref{lambdafix} if $\me{H}\ra\me{F}$ is an isomorphism on $\C^\times p$
then $\Phi_{U,\gamma}(\me{H})\ra\Phi_{U,\gamma}(\me{F})$ is an isomorphism at 
$\pi(p)$ for sufficently small $(U,\gamma)$. Consider the following chain of 
morphisms:
\begin{align*}
  \Hom{\Db_{\Cc}(k_X,\C^{\times} p)}&(\me{F},\me{G})=
  \underset{\me{H}\ra\me{F}\above 0pt iso\ on\ \C^\times p}
       {\varinjlim}\Hom{\Db(k_X)}(\me{H},\me{G}) \\
    & \lra \underset{\Phi_{U,\gamma}(\me{H})\ra\Phi_{U,\gamma}(\me{F}) 
     \above 0 pt iso\ on\ \pi(p)}
      {\varinjlim}\Hom{\Db(k_X)}(\Phi_{U,\gamma}(\me{H}),
      \Phi_{U,\gamma}(\me{G}))\\
   &\lra \underset{\me{H}'\ra\Phi_{U,\gamma}(\me{F})\above 0 pt iso\ on
      \ \pi(p)}{\varinjlim}
    \Hom{\Db(k_X)}(\me{H}',\Phi_{U,\gamma}\me{G})
   =\Hom{\Db_{\Cc}(k_X,\pi(p))}(\Phi_{U,\gamma}(\me{F}),
    \Phi_{U,\gamma}(\me{G})). 
\end{align*}
We have to check that the composition is an isomorphism for a sufficently small 
cutting pair $(U,\gamma)$.\\
Consider a morphism $\me{H}'\ra\Phi_{U,\gamma}(\me{G})$ and an isomorphism 
$\me{H}'\ra\Phi_{U,\gamma}(\me{F})$ on $\pi(p)$. Then $\SS(\me{H}')=
\SS(\Phi_{U,\gamma}(\me{F}))$ in $\pi^{-1}(V)$  for some neighborhood
$V$ of $\pi(p)$ and in particular $\me{H}'\in\Db_{\Cc,\Lambda}(k_X,\C^{\times} p)$, 
hence $\Phi_{U,\gamma}(\me{H}')\ra \me{H}'$ is an isomorphism in $\pi(p)$. 
Then $\me{H}'\ra\Phi_{U,\gamma}(\me{G})\ra\me{G}$ is sent to 
$\Phi(\me{H}')\ra \Phi(\me{G})$ which is the same morphism in 
$\Db_{\Cc}(k_X,\pi(p))$ as $\me{H}'\ra\Phi_{U,\gamma}(\me{G})$.
Hence the map is surjective.\\
Let us show that it is injective. If $\Phi_{U,\gamma}(\me{H})\ra
\Phi_{U,\gamma}(\me{G})$ is zero in $\Db_{\Cc}(k_X,\C^{\times} p)$ then there 
exists
$\me{K}$ and an isomorphism $\me{K}\ra\Phi_{U,\gamma}(\me{H})\ra
\Phi_{U,\gamma}(\me{F})$ on $\pi(p)$ such that $\me{K}\ra\Phi(\me{H})\ra
\Phi(\me{G})$ is the zero morphism. Then $\me{K}\ra\Phi(\me{H})\ra\me{H}\ra
\me{F}$ is an isomorphism on $\C^{\times} p$ and $\me{K}\ra\Phi(\me{H})\ra\me{H}\ra
\me{G}$ is zero hence $\me{H}\ra\me{G}$ is zero in 
$\Db_{\Cc}(k_X,\C^{\times} p)$.
\end{proof}
\begin{remark}
  \em Note that
$$ \Db_{\Cc}(k_{\,*\,})_{\pi(p)} \overset{\sim}{\lra}  \Db_{\Cc}(k_X,\pi(p)), $$
  where $\Db_{\Cc}(k_{\,*\,})$ denotes the prestack $U\mapsto \Db_{\Cc}(k_U)$. 
  In particular we get the fully faithful functor
  $$ \Db_{\Cc,\Lambda}(k_X,\C^{\times} p) \lra \Db_{\Cc}(k_{\,*\,})_{\pi(p)}. $$
  Also recall that $\Db_{\Cc,\Lambda}(k_X,\C^{\times} p)$ is  additive (because of 
  the triangular inequality for the micro-support). Hence 
  $\Db_{\Cc,\Lambda}(k_X,\C^{\times} p)$ is a full additive subcategory of 
  $\Db_{\Cc}(k_{\,*\,})_{\pi(p)}$.\\
  Therefore we have shown that if we are given $\me{F}\in\Db_{\Cc}(k_X,\C^\times p)$ 
  such that  $\SS(\me{F})$ is in generic position at $p$ then there exists an object
  $\widetilde{\me{F}}\in\Db_{\Rc}(k_X)$ that is $\C$-constructible in a neighborhood
  of $\pi(p)$. Moreover we can construct $\widetilde{\me{F}}$ in a functorial way. 
\end{remark}

\subsection{The generic position theorem} 

Let us recall Kashiwara-Kawai's generic position theorem (cf. \cite{KK}, Chapter I.6). 
\begin{prop}
  Let $p\in\dot T^*X$ and $\Lambda\subset T^*X$ be an $\R^+$-conic Lagrangian
  subset such that $\Lambda$ is $\C^{\times}$-conic in a neighborhood of
  $\C^{\times} p$. Then there exists a complex contact transformation
  $$\chi:\, T^*X\ra T^*X, $$
  defined in a neighborhood of $\C^\times p$, such that $\chi(\Lambda)$ is in 
  generic position at $q=\chi(p)$.
\end{prop}
\begin{proof}
Kashiwara-Kawai show this result on a neighborhood of 
$p$ for locally $\C^\times$-conic Lagrangian varieties assuming only that 
$\Lambda\cap\pi^{-1}\pi(p)=\C^\times p$ in a neighborhood of $p$. Hence if we suppose
that $\Lambda$ is $\C^{\times}$-conic in a neighborhood of $\C^{\times} p$, we get
that $\chi(\Lambda)$ is in generic position at $q$ since $\chi$ is 
$\C^\times$-homogenous.
\end{proof}
\begin{remark}
  \em Consider a microlocal complex constructible sheaf $\me{F}\in\Db_{\Cc}(k_X,
  \C^\times p)$. By the generic position theorem and invariance under quantized 
  contact transformations we can find a contact transformation 
  $\chi:\,T^*X \ra T^*X$ defined in a neighborhood of $\C^\times p$ and 
  an equivalence of categories
  $$ \Phi^{\Cc}_{\me{K}}:\ \Db_{\Cc}(k_X,\C^\times p) \overset{\sim}{\lra} 
       \Db_{\Cc}(k_X,\C^\times q) $$
  such that $\SS(\Phi^{\Cc}_{\me{K}}(\me{F}))=\chi(\SS(\me{F}))$ is in generic 
  position at $q$. Then $\Phi_{\me{K}}^{\Cc}(\me{F})$ is functorially isomorphic 
  in $\Db_{\Cc}(k_X,\C^\times q)$ to an object of $\Db_{\Cc}(\,*\,)_{\pi(q)}$. 
  Hence many problems in $\Db_{\Cc}(k_X,\C^\times p)$ can be reduced to the study of 
  germs of complex constructible sheaves in generic position.
\end{remark}

\section{Microlocal perverse sheaves in $\Db_{\Rc}(k_X,S)$}

\subsection{Andronikof's prestack of microlocally perverse sheaves}

In this section we will first recall Andronikof's definition of the prestack of 
microlocal perverse sheaves. In \cite{An2} he suggests a definiton of a 
microlocal perverse sheaf on an arbitrary subset $S$ of $T^*X$ which is based 
on the microlocal characterisation of perverse sheaves (see Proposition 
\ref{pervsheaf} below). However he is not very precise concerning 
$\C^{\times}$-conicity. In particular the proof of his main tool (Proposition 3.2) 
is incomplete under the given assumptions.\\ 
Hence we will recall a slightly more precise version of his prestack and we 
will restrict ourselves from the beginning to $\C^{\times}$-conic sets. Moreover 
we will add complete proofs to the statements of \cite{An2} which hold in this 
case. Also while Andronikof restricts his studies to points and $\C^\times$-orbits
we will work from the beginning with the entire prestack.\\
In \cite{KS2} we find the following microlocal characterization of
perverse sheaves.
\begin{prop}\label{pervsheaf}
  Let $\me{F}\in \on{D}^{\on{b}}_{\C\text{-}c}(k_X)$. Then we have equivalence
  between
\begin{itemize}
   \item[(P1)] $\me{F}$ is a perverse sheaf.
   \item[(P2)] For every non-singular point $p$ of $\on{SS}(\me{F})$
     such that the projection $\pi:\,\on{SS}(\me{F})\ra X$ has constant
     rank on a neighborhood of $p$, there exist a submanifold
     $Y\subset X$ and an object $M\in\on{Mod}(k)$ such that 
     $\me{F}\simeq M_Y[\dim Y]$ in $\Db(k_X,p)$.
   \item[(P3)] The assertion (P2) is true for some point $p$ of any irreducible
     component of $\SS(\me{F})$.
\end{itemize}
\end{prop}
This naturally leads to the following
\begin{defi}
  Let $S\subset \dot T^*X$ be a $\C^{\times}$-conic subset. 
  \begin{itemize}
      \item[(i)] A sheaf $\me{F}\in\Db_{\Rc}(k_X)$ is called microlocally 
          perverse on $S$, if it is microlocally $\C$-constructible on $S$ 
          and there exists an open neighborhood $U$ of $S$ such that for every 
          non-singular point $p$ of $\on{SS}(\me{F})\cap U$ such that the 
          projection $\pi:\,\on{SS}(\me{F})\ra X$ has constant rank on a 
          neighborhood of $p$, there exists a submanifold $Y\subset X$ and an object
          $M\in\on{Mod}(k)$ such that $\me{F}\simeq M_Y[\dim Y]$ in $\Db(k_X,p)$.
      \item[(ii)] We will denote by $\Db_{\on{perv}}(k_X,S)$ the full subcategory
          of $\Db_{\Cc}(k_X,S)$ whose objects are perverse on $S$ (i.e. which 
          may be represented by a microlocally perverse sheaf on $S$).
      \item[(iii)] In particular if  $\Omega\subset \dot T^*X$ is an open 
          $\C^{\times}$-conic subset, we get the category 
          $\on{D}^b_{\on{perv}}(k_X,\Omega)$. Clearly this 
          defines a prestack of categories on $P^*X$, denoted by 
          $\Db_{\on{perv}}(k_X,\,*\,)$. This is Andronikof's prestack of microlocal 
          perverse sheaves. 
      \item[(iv)] Let $\Lambda\subset T^*X$ be an $\R^+$-conic Lagrangian variety 
          that is $\C^\times$-conic in a neighborhood of $\C^\times p$. Then we set
          $$ \Db_{\on{perv},\Lambda}(k_X,\C^\times p)=
             \Db_{\Cc,\Lambda}(k_X,\C^\times p)\cap
             \Db_{\on{perv}}(k_X,\C^\times p). $$ 
\end{itemize}
\end{defi}
\begin{prop}
  Let $p\in\dot T^*X$. The stalk of Andronikof's prestack at $\gamma(p)$ 
  is precisely the category $\on{D}^b_{\on{perv}}(k_X,\C^{\times} p)$.
\end{prop}
\begin{proof}
Consider the functor
$$  \twodilim{\C^{\times} p\subset U\subset \dot T^*X \above 0pt
   \text{$U$ $\C^\times$-conic}}{\on{D}^b_{\on{perv}}(k_X,U)} 
      \lra \on{D}^b_{\on{perv}}(k_X,\C^{\times} p). $$
It is essentially surjective by definition of 
$\on{D}^b_{\on{perv}}(k_X,\C^{\times} p)$.\\
The proof that it is fully faithful is analogous to Proposition 
\ref{germformula} and Corollary \ref{cgermformula}.
\end{proof}
\begin{prop}
  The duality functor 
  $$ \on{D}:\ \Db_{\Rc}(k_X) \lra \Db_{\Rc}(k_X) \quad ;\quad
      \me{F}\mapsto \Dr\me{H}om(\me{F},\omega_X) $$
  induces contravariant equivalences of prestacks
  $$ \on{D}:\ \Db_{\Cc}(k_X,\,*\,) \lra \Db_{\Cc}(k_X,\,*\,), $$  
  $$ \on{D}:\ \Db_{\on{perv}}(k_X,\,*\,) \lra \Db_{\on{perv}}(k_X,\,*\,). $$  
\end{prop}
\begin{proof}
The first functor is well-defined because $\SS(\on{D}\me{F})=\SS(\me{F})$.\\
Let $p\in\SS(\me{F})$ such that $\me{F}\simeq M_Y[\dim X]$ in $\Db(k_X,p)$. Then
$\on{D}\me{F}\simeq \on{D}M_Y[\dim Y]$ in $\Db(k_X,p)$ since $D$ is an anti-equivalence
in $\Db(k_X,p)$. But $\on{D}M_Y[\dim Y]\simeq M_Y[\dim_Y]$. Hence the second functor 
is well-defined.\\
The two functors are equivalences because $\on{D}^2=\on{Id}$ in $\Db_{\Rc}(k_X)$.
\end{proof}
Andronikof's prestack is invariant by quantized contact transformations:
\begin{prop}\label{perversekerneltheorem}
   Let $X,Y$ be two complex manifolds, $\Omega_X\subset \dot T^*X$, $\Omega_Y\subset 
   \dot T^*Y$ open $\C^\times$-conic subsets and
   $$ \chi:\ \Omega_X\overset{\sim}{\lra} \Omega_Y $$ 
   a contact transformation. Set
$$ \Lambda=\Big\{((y;\eta),(x;\xi))\in\Omega_Y\times\Omega_X^a\ |\ (y,\eta)=
   \chi(x,-\xi)\Big\}. $$
   Then $\Lambda$ is $\C^\times$-conic. Let $p_X\in\Omega_X$ and $p_Y=\chi(p_X)$.\\
   There exist open neighborhoods $X'$ of $\pi(p_X)$, $Y'$ of $\pi(p_Y)$,
   $\C^\times$-conic open neighborhoods $\Omega_X'$ of $\C^\times p_X$, $\Omega_Y'$ of 
   $\C^\times p_Y$ with $\Omega_X'\subset T^*X'\cap \Omega_X$, 
   $\Omega_Y'\subset T^*Y'\cap \Omega_Y$ and a kernel
   $\me{K}\in\Db_{\Rc}(k_{Y'\times X'})$ satisfying 
   \begin{itemize}
     \item[(1)] $\chi$ induces a homogenous contact transformation 
         $\Omega_X'\overset{\sim}{\ra} \Omega_Y'$,
     \item[(2)] $$\Big((\Omega_Y''\times T^*X')\cup (T^*Y'\times {\Omega_X''}^a)\Big)
       \cap\SS(\me{K})\subset \Lambda \cap (\Omega_Y''\times {\Omega_X''}^a)$$ 
       for every $\C^\times$-conic open subset $\Omega_{X}''\subset\Omega_{X}'$ and 
       $\Omega_{Y}''=\chi(\Omega_{X}'')$,
     \item[(3)] the functor $\Phi_{\me{K}}^{\Cc}$ induces an equivalence of prestacks
       $$\Phi_{\me{K}}^{\on{perv}}:\, 
              \chi_*\Db_{\on{perv}}(k_{X'},\,*\,)|_{\Omega_{X}'} \lra 
             \Db_{\on{perv}}(k_{Y'},\,*\,)|_{\Omega_{Y}'}.$$
     \item[(4)] $\SS(\Phi_{\me{K}}^{\on{perv}}(\me{F}))\cap \Omega_{Y}''=
                   \chi(\SS(\me{F})\cap \Omega_{X}'')$
     \item[(5)] $\chi_*\mu hom(\me{F},\me{G})|_{\Omega_X'}\simeq 
        \mu hom(\Phi^{\on{perv}}_{\me{K}}(\me{F}),
        \Phi_{\me{K}}^{\on{perv}}(\me{G}))|_{\Omega_Y'}$.
   \end{itemize}
\end{prop}
\begin{proof}
It is enough to show the proposition in the case in which $\Lambda$ is the conormal
bundle to a smooth hypersurface $S$ and we can choose $\me{K}\simeq k_S$. Then
the fact that $\Phi_{\me{K}}^{\Cc}$ preserves microlocal perverse sheaves follows
from Proposition 7.4.6 of \cite{KS3}.
\end{proof}

\subsection{The abelian category $\Db_{\on{perv}}(k_X,\C^{\times} p)$} 

In general, one does not know much about the category $\Db_{\on{perv}}(k_X,S)$ even if
$S$ is $\C^\times$-conic. In \cite{An2}, it is announced that if
$p\in\dot T^*X$, then $\Db_{\on{perv}}(X,p)$ and 
$\Db_{\on{perv}}(X,\C^\times p)$ are abelian. While we do not know if this is true
for $\Db_{\on{perv}}(X,p)$, we will give a proof here for the category 
$\Db_{\on{perv}}(X,\C^\times p)$. The main tools for the proof have been developped
in Section 6.

We will fix a point $p\in\dot T^*X$ and an $\R^+$-conic Lagrangian subvariety 
$\Lambda$ which is $\C^\times$-conic in a neighborhood of $\C^{\times}p$. 
\begin{prop}
  Suppose that $\Lambda$ is in generic position at $p$. Then the functor $\Phi$ of 
  Theorem \ref{dafunctorphi} induces a commutative diagram of fully faithful functors
  $$ \xymatrix{
     {\Db_{\Cc,\Lambda}(k_X,\C^{\times} p)} \ar[r]^{\Phi} & {\Db_{\Cc}(k_X)_{\pi(p)}}\\ 
     {\Db_{\on{perv},\Lambda}(k_X,\C^{\times} p)} \ar[r]\ar[u] & 
      {\Perv(X)_{\pi(p)}.} \ar[u] } $$
\end{prop}
\begin{proof}
It is sufficent to prove that if $\me{F}\in\Db_{\on{perv}}(k_X,\C^\times p)$ then
$\Phi(\me{F})$ is perverse in a neighborhood of $\pi(p)$. Since we have
$$ \SS(\Phi(\me{F}))\cap \dot\pi^{-1}(V)=\SS(\me{F})\cap (V\times \gamma)$$
for some small neighborhood $V\times\gamma$ of $\C^\times p$, we get the result
by the characterization of perverse sheaves (cf. Proposition \ref{pervsheaf}).
\end{proof}
Also note that $\Db_{\on{perv},\Lambda}(k_X,\C^{\times} p)$ is additive. Hence 
$\Db_{\on{perv},\Lambda}(k_X,\C^{\times} p)$ is a full additive subcategory of 
$\Perv(X)_{\pi(p)}$. In order to prove that it is actually a full abelian 
subcategory, we will need two lemmas:
\begin{lemma}\label{ssofkernels}
  Let $\varphi:\,\me{F}\ra\me{G}$ be a morphism of perverse sheaves such that
  $\SS(\me{F})\cup \SS(\me{G})\subset \Lambda$ in a neighborhood of $p$. Then
  $\SS(\ker(\varphi))\cup\SS(\coker(\varphi))\subset\Lambda$ in a neighborhood at 
  $p$. 
\end{lemma}
\begin{proof}
Recall (\cite{KS3}, Exercise X.6) that the micro-support of a $\C$-constructible 
sheaf $\me{F}\in\Db_{\Cc}(k_X)$ can be calculated as 
$$ \SS(\me{F})=\bigcup_{i\in \Z}{}^{\on{p}}\!\on{H}^i(\me{F}) $$
where ${}^{\on{p}}\!\on{H}^i(\me{F})$ denotes the $i^{th}$ perverse cohomology 
sheaf of $\me{F}$.\\
Let $\varphi:\,\me{F}\ra\me{G}$ be a morphism of perverse sheaves. We 
embed it into a distinguished triangle
$$ \xymatrix{
     {\me{F}} \ar[r] & {\me{G}} \ar[r] & {\me{H}} \ar[r]^+ & } $$
Then we consider the canonical distinguished triangle
$$  \xymatrix{ 
  {{}^{\on{p}}\tau^{\leqs -1}(\me{H})} \ar[r] & {\me{H}} \ar[r] &
    {{}^{\on{p}}\tau^{\geqs 0}(\me{H})} \ar[r]^+ & } $$
where ${}^{\on{p}}\tau^{\leqs -1},{}^{\on{p}}\tau^{\geqs 0}$ denote the perverse 
truncation functors. Then ${}^{\on{p}}\tau^{\leqs -1}(\me{H})[-1]$ is the kernel 
and ${}^{\on{p}}\tau^{\geqs 0}(\me{H})$ is the cokernel of $\varphi$. Therefore
$$ \SS(\ker\varphi)\cup \SS(\coker\varphi)=\SS(\me{H})\subset
    \SS(\me{F})\cup \SS(\me{G}). $$
\end{proof}
\begin{lemma}\label{stability}
  Let $\Lambda$ be in generic position at $p$.\\
  Consider a morphism $\me{F}\ra\me{G}$ of $\Db_{\on{perv},\Lambda}(k_X,\C^\times p)$.
  Let $\Phi(\me{F})\ra\Phi(\me{G})$ be the corresponing morphism in 
  $\Perv_{\pi(p)}$ and $\me{L}\ra \Phi(\me{F})$ a kernel (resp.
  $\Phi(\me{G})\ra\me{L}'$ a cokernel) in $\Perv_{\pi(p)}$. Then
  $\Phi(\me{L})\ra\me{L}$ (resp. $\Phi(\me{L}')\ra\me{L}'$) is an isomorphism in 
  $\Perv_{\pi(p)}$.
\end{lemma}
\begin{proof}
The proof follows the same idea as the proof of Lemma \ref{compatibilitycondition}.\\
Fix $U,\gamma$ such that $\Phi(\me{F})\simeq \Phi_{U,\gamma}(\me{F})$,
$\Phi(\me{G})\simeq \Phi_{U,\gamma}(\me{G})$ and $\Phi(\me{L})\simeq 
\Phi_{U,\gamma}(\me{L})$. Embed $\Phi_{U,\gamma}(\me{L})\ra\me{L}$ in a distinguished 
triangle
$$ \xymatrix{
    {\Phi_{U,\gamma}(\me{L})} \ar[r] & {\me{L}} \ar[r] & {\me{H}} \ar[r]^+ & } $$
Then there is an open neighborhood $V\subset U$ of $\pi(p)$ such that
$$ \SS(\Phi_{U,\gamma}(\me{L}))\cap\dot\pi^{-1}(V)=
   \SS(\me{L})\cap V\times\gamma. $$
Moreover, by Lemma \ref{ssofkernels}, we have
$$ \SS(\me{L})\cap \dot \pi^{-1}(V)\subset
   \SS(\Phi_{U,\gamma}(\me{F}))\cup \SS(\Phi_{U,\gamma}(\me{G})) \cap
   \dot\pi^{-1}(V)\subset \Lambda \cap V\times\gamma. $$
Hence if $V$ is sufficently small we get
$$ \SS(\me{L})\cap \dot \pi^{-1}(V)\subset \Lambda \cap V\times\gamma. $$
Since $\Phi_{U,\gamma}(\me{L})\ra\me{L}$ is an isomorphism on $U\times\gamma$ we get
that $\SS(\me{H})\cap V\times\gamma=\varnothing$ and therefore 
$\SS(\me{H})\subset T^*_XX$ on $V$. Then 
$$ \SS(\me{H})\cap V \subset (\SS(\me{L})\cup \SS(\Phi_{U,\gamma}(\me{L})))
  \cap V \subset \pi(\Lambda)\cap V.$$
Since $\Lambda$ is in generic position, locally at $\pi(p)$ the set $\pi(\Lambda)$ is
a closed hypersurface. Therefore the involutivity of the micro-support implies that
$\me{H}\simeq 0$ in a neighborhood of $\pi(p)$.\\
The proof for the cokernel is similar.
\end{proof}
\begin{prop}\label{thecategoryisabelian}
  The additive category $\Db_{\on{perv},\Lambda}(k_X,\C^{\times} p)$ 
  is equivalent to an abelian subcategory of $\Perv_{\pi(p)}$.
\end{prop}
\begin{proof}
By Lemma \ref{stability} the full additive subcategory $\Db_{\on{perv},\Lambda}
(k_X,\C^{\times} p)$ is stable by kernels and cokernels. Since it is a full 
subcategory, it is abelian.
\end{proof}
\begin{lemma}
  Let $\Lambda$ be a $\C^{\times}$-conic Lagrangian variety (we do not ask
  $\Lambda$ to be in generic position at $p$). Then $\Db_{\on{perv},\Lambda}(k_X,
  \C^{\times} p)$ is abelian.\\
  Moreover if $\Lambda'$ is another $\C^\times$-conic Lagrangian variety with
  $\Lambda\subset\Lambda'$, then the natural functor
  $$ \Db_{\on{perv},\Lambda}(k_X,\C^{\times} p) \lra 
     \Db_{\on{perv},\Lambda'}(k_X,\C^{\times} p) $$
  is exact. 
\end{lemma}
\begin{proof}
Consider $\Lambda\subset\Lambda'$. Let $\chi$ be a canonical transformation
such that $\chi(\Lambda')$ is in generic position at $q=\chi(p)$. Then 
$\chi(\Lambda)$ is also at generic position at $q$ and we get a diagram
$$ \xymatrix{
   {\Db_{\on{perv},\Lambda}(k_X,\C^{\times} p)} \ar[r]^(.45)\sim  \ar[d]_{\alpha}  &
   {\Db_{\on{perv},\chi(\Lambda)}(k_X,\C^{\times} q)}\ar[d]^{\beta} \ar[r] &
   {\Perv(X)_{\pi(q)}} \ar[d]^=  \\ 
   {\Db_{\on{perv},\Lambda'}(k_X,\C^{\times} p)} \ar[r]^(.45)\sim  &  
   {\Db_{\on{perv},\chi(\Lambda')}(k_X,\C^{\times} q)} \ar[r] &
     {\Perv(X)_{\pi(q)}}.  } $$
This diagram is commutative up to isomorphism. The horizontal functors are exact and
fully faithful. By Proposition \ref{thecategoryisabelian} the categories 
$\Db_{\on{perv},\chi(\Lambda)}(k_X,\C^{\times} q)$ and
$\Db_{\on{perv},\chi(\Lambda')}(k_X,\C^{\times} q)$ are abelian subcategories of
$\Perv(X)_{\pi(q)}$. Hence $\beta$ (and therefore $\alpha$) is exact.
\end{proof}
\begin{prop}\label{abelianstalks}
  The category $\Db_{\on{perv}}(k_X,\C^{\times} p)$ is abelian. Moreover,  for every 
  germ of a $\C^\times$-conic Lagrangian variety $\Lambda\subset T^*X$ defined in 
  a neighborhood of $\C^{\times} p$, the inclusion functor
  $$ \Db_{\on{perv},\Lambda}(k_X,\C^{\times} p) \lra
     \Db_{\on{perv}}(k_X,\C^{\times} p) $$
  is exact.
\end{prop}
\begin{proof}
We have
$$ \twodilim{\Lambda\supset\C^{\times} p}{\Db_{\on{perv},\Lambda}(k_X,\C^{\times} p)}
   \overset{\sim}{\lra}\Db_{\on{perv}}(k_X,\C^{\times} p). $$
Filtered 2-colimits of abelian categories with exact restriction functors are abelian.
\end{proof}
Now let us prove the following ``lifting property'' for kernels and cokernels
of microlocal perverse sheaves:
\begin{prop}\label{liftingkernels}
  Let $\varphi:\,\me{F}\ra\me{G}$ be a morphism in 
  $\Db_{\on{perv}}(k_X,\C^{\times} p)$. Then 
  there exists a neighborhood $V$ of $\C^{\times} p$, objects 
  $\me{K},\me{K}'\in\Db_{\on{perv}}(k_X,V)$ and morphisms $\me{K}\ra\me{F}$, 
  $\me{G}\ra\me{K}'$ in $\Db_{\on{perv}}(k_X,V)$ such that these morphisms 
  induce kernel and cokernel of $\varphi$ in 
  $\Db_{\on{perv}}(k_X,\C^{\times} q)$ for all $q\in V$.
\end{prop}
\begin{proof}
Choose a Lagrangian variety $\Lambda$ such that $\SS(\me{F})\cup 
\SS(\me{G})\subset \Lambda$ in a neighborhood of $\C^\times p$.\\
By the generic position theorem, we can find a contact transformation defined in a 
$\C^\times$-conic open neighborhood of $\C^\times p$
$$ \chi:\ (T^*X,\C^\times p) \lra (T^*X,\C^\times p') $$
such that $\chi(\Lambda)$ is in generic position at $p'$. Then $\chi(\Lambda)$
is isomorphic to the conormal bundle to a closed hypersurface in a neighborhood of
$\C^\times p'$. In particular, there exists a $\C^\times$-conic open neighborhood $V'$
of $\chi(\Lambda)$ such that $\chi(\Lambda)$ is in generic position at any point of
$\chi(\Lambda)\cap V'$.\\
By Proposition \ref{perversekerneltheorem}, we can assume that 
$\Lambda$ is in generic position at any point in a $\C^\times$-conic neighborhood 
$\Omega$ of $\C^\times p$.\\
Consider the functor
$$ \Phi:\ \Db_{\on{perv},\Lambda}(k_X,\C^\times p) \lra \Perv_{\pi(p)}.$$
Choose $(U,\gamma)$ such that the morphism $\Phi(\me{F})\ra\Phi(\me{G})$ is defined 
as $\Phi_{U,\gamma}(\me{F})\ra \Phi_{U,\gamma}(\me{G})$ in $\Perv(V)$ for some small
neighborhood $V$ of $\pi(p)$. Let $\me{L}\ra\Phi_{U,\gamma}(\me{F})$ (resp.
$\Phi_{U,\gamma}(\me{G})\ra\me{L}')$ be a kernel (resp. a cokernel) in $\Perv(V)$.
Choose a $\C^\times$-conic open neighborhood of
$\C^\times p$ such that $\Omega'\subset\Omega\cap V\times\gamma$. We will show that
$\me{L}\ra\Phi_{U,\gamma}(\me{F})$ (resp. $\Phi_{U,\gamma}(\me{G})\ra\me{L}')$ is a 
kernel (resp. a cokernel) in $\Db_{\on{perv},\Lambda}(k_X,\C^\times q)$ for all $q$
in $\Omega'$.\\
We have distinguished triangles in $\Db(k_X)$
$$ \xymatrix{
   {\Phi_{U,\gamma}(\me{F})} \ar[r] & {\Phi_{U,\gamma}(\me{G})} \ar[r] & {\me{H}} 
 \ar[r]^+ &  }$$ 
and 
$$ \xymatrix{
     {\me{L}[1]} \ar[r] & {\me{H}} \ar[r] & {\me{L}'} \ar[r]^+ & }$$
where $\me{L}$ is the kernel (resp. $\me{L}'$ the cokernel) of  
$\Phi_{U,\gamma}(\me{F})\ra \Phi_{U,\gamma}(\me{G})$ in $\Perv(V)$.
Now let $q$ be a point of $\Lambda\cap \Omega'$. Since $\Lambda$ is in generic 
position at $q$ we get the functor
$$ \Phi^q:\ \Db_{\on{perv},\Lambda}(k_X,\C^\times q) \lra
     \Perv_{\pi(q)}. $$
Applying $\Phi^q$ to the two triangles, we get disinguished triangles
$$ \xymatrix{
   {\Phi^q\Phi_{U,\gamma}(\me{F})} \ar[r] & {\Phi^q\Phi_{U,\gamma}(\me{G})} \ar[r] & 
   {\Phi^q(\me{H})} \ar[r]^+ & }$$
and  
$$ \xymatrix{
    { \Phi^q\me{L}[1]} \ar[r] & {\Phi^q\me{H}} \ar[r] & {\Phi^q\me{L}'} \ar[r]^+ & }$$
Since $\Phi^q(\me{L})$ and $\Phi^q(\me{L}')$ are perverse in a neighborhood of
$\pi(q)$ we get (by construction of kernels (resp. cokernels) in $\Perv_{\pi(q)}$)
that $\Phi^q(\me{L}) \ra  \Phi^q\Phi_{U,\gamma}(\me{F})$ (resp. 
$\Phi^q\Phi_{U,\gamma}(\me{G})\ra\Phi^q(\me{L}')$) is a kernel (resp. cokernel) in 
$\Db_{\on{perv},\Lambda}(k_X,\C^\times q)$. 
Finally since $\Phi^q$ is a cut-off functor we have $\Phi^q(\me{L})\simeq \me{L}$ in 
$\Db_{\on{perv},\Lambda}(k_X,\C^\times q)$.
\end{proof} 
\begin{thm}\label{dafirsttheorem}
  The stack associated to Andronikof's prestack of microlocal perverse sheaves is 
  abelian.
\end{thm}
\begin{proof}
We have shown that the stalks of this additive prestack are 
abelian categories (Proposition \ref{abelianstalks}). Further we have shown that
kernels and cokernels in the stalks may be lifted to small open neighborhoods 
(Proposition \ref{liftingkernels}). Therefore the conditions of Proposition 
\ref{abstack} are satisfied and the stack of microlocal perverse sheaves is abelian.
\end{proof}

\section{Microlocal perverse sheaves on $P^*X$}

\subsection{The stack of microlocal perverse sheaves}

We are now ready to give a first definition of the stack of microlocal perverse 
sheaves.
\begin{defi}
  The stack of microlocal perverse sheaves on $P^*X$ is the stack associated to
  Andronikof's prestack.
\end{defi}
By Theorem \ref{dafirsttheorem} we know that the stack of microlocal perverse
sheaves is abelian. Furthermore since the underlying prestack is invariant by quantized
contact transformations (by Proposition \ref{perversekerneltheorem}) we easily get
that the stack of microlocal perverse sheaves is invariant by quantized contact
transformations.\\
We will now give an explicit description of microlocal perverse 
sheaves in terms of ind-sheaves. For this purpose we will construct a 
subprestack 
$$\mu Perv\subset \gamma_*(\Db(\I(k_*))|_{\dot T^*X}).$$  
Here, $\gamma_*(\Db(\I(k_*))|_{\dot T^*X})$ is the prestack of bounded derived 
categories of ind-sheaves on $\C^\times$-conic open subsets of 
$\dot T^*X$. Then we will show that $\mu Perv$ is actually a stack and construct
a morphism $\mu:\, \Db_{\on{perv}}(k_X,*)\ra \mu Perv$ and prove that it 
induces equivalences in the stalks. Hence $\mu Perv$ can be identified to
the stack associated to $\Db_{\on{perv}}(k_X,*)$. In particular it is an abelian stack.
Then we will be able to define the microlocal Riemann-Hilbert functors using the
theory of ind-sheaves.\\
First note that the functor $\mu$ induces a functor of prestacks
\begin{equation} \label{functormuperv}
 \mu:\ \Db_{\on{perv}}(k_X,\,*\,) \lra \gamma_*\Db(\I(k_{\,*\,})) 
\end{equation}
where we consider $\Db(\I(k_{\,*\,}))$ as a prestack on $\dot T^*X$. 
\begin{defi}
  Let $\Omega\subset P^*X$ be an open subset.
  \begin{itemize}
     \item[(1)] An object $\me{F}\in\Db(\I(k_{\gamma^{-1}\Omega}))$ is 
        microlocally perverse (on $\Omega$) if it is locally in the image of the 
        functor \eqref{functormuperv}, i.e. if for all $p\in\gamma^{-1}(\Omega)$ 
        there exists a $\C^\times$-conic neighborhood $V\supset \C^\times p$ 
        and an object $\me{G}\in \on{D}^{\on{b}}_{perv}(k_X,V)$ such that
        $\mu\me{G}|_{V}\simeq \me{F}|_{V}$.
     \item[(2)] We denote by $\mu Perv(\Omega)$ the full subcategory of
        $\Db(\I(k_{\gamma^{-1}\Omega}))$ whose objects are microlocally 
        perverse.
   \end{itemize}
\end{defi}
\begin{remark}
  \em The functor \eqref{functormuperv} induces a functor in the stalks
  $$ \mu:\ \Db_{\on{perv}}(k_X,\C^\times p) \lra 
  \gamma_*\Db(\I(k_{\,*\,}))_{\gamma(p)}$$
  and the definition of a microlocal perverse sheaf is clearly equivalent to
  \begin{itemize}
     \item[(1')] An object $\me{F}\in\Db(\I(k_{\gamma^{-1}\Omega}))$ is 
        microlocally perverse (on $\Omega$) if for all $p\in\gamma^{-1}(\Omega)$ 
        there exists an object $\me{G}\in \on{D}^{\on{b}}_{perv}(k_X,\C^\times p)$ such 
        that $\mu\me{G}\simeq \me{F}$ in $\gamma_*\Db(\I(k_{\,*\,}))_{\gamma(p)}$.
  \end{itemize}
  Here, one shall keep in mind that the natural functor
$$ \gamma_*\Db(\I(k_*))_{\gamma(p)} \lra
     \Db(\I(k_{\C^\times p})) $$
  is not fully faithful. Therefore germs of microlocal perverse sheaves 
  should not be interpreted as complexes of ind-sheaves on $\C^\times p$.
  \end{remark}
  \begin{remark}
  \em Note that it is possible (but less obvious) to define the functor 
  $$\mu:\,\Db_{\on{perv}}(k_X,\C^\times p) \lra \gamma_*\Db(\I(k_{\,*\,}))_{\gamma(p)}$$
  directly without using the prestack $\Db_{\on{perv}}(k_X,\,*\,)$. Hence, we see by 
  definition (1') that in order to define microlocal perverse sheaves, it is not 
  necessary to construct the categories $\Db_{\on{perv}}(k_X,S)$ for any 
  $\C^{\times}$-conic subset $S\subset T^*X$ but only the categories 
  $\Db_{\on{perv}}(k_X,\C^\times p)$ for any $p\in \dot T^*X$. 
  Recall that the categories $\Db_{\Rc}(k_X,\C^\times p)$ can be identified with 
  the full subcategory of $\Db(k_X,\C^\times p)$ whose objects may be represented by 
  $\R$-constructible sheaves. Hence the construction of microlocal perverse sheaves is 
  independent of the conceptual question of the ``correct'' definition of the 
  microlocalization of $\R$-constructible sheaves (cf. Section 4).\\
  However, the prestack $\Db_{\on{perv}}(k_X,\,*\,)$ can sometimes be useful to define 
  functors on microlocal perverse sheaves and therefore we decided to include it in 
  our presentation. \em
\end{remark}
Clearly $\mu Perv(\Omega)$ is an additive subcategory of 
$\Db(\I(k_{\gamma^{-1}(\Omega)}))$, and the correspondance 
$$  P^*X\supset \Omega\mapsto \mu Perv(\Omega) $$ 
defines an additive prestack on $P^*X$. Hence $\mu Perv$ is a full additive 
subprestack of $\gamma_*\Db(\I(k_*))$, and we can see directly from the construction 
that it is defined by a local property. Moreover, by definition, for any open subset 
$\Omega\subset P^*X$ the functor $\mu$ induces a natural functor 
$$ \Db_{\on{perv}}(k_X,\gamma^{-1}(\Omega)) \lra
     \mu Perv(\Omega). $$
These functors define a functor of prestacks
\begin{equation}\label{mupervequivalence}
  \mu:\ \Db_{\on{perv}}(k_X,\,*\,) \lra 
      \mu Perv.
\end{equation}
We want to show that the functor \eqref{mupervequivalence} induces equivalences of 
categories in the stalks. Let $p\in\dot T^*X$. The definition of $\mu Perv$ 
immediately implies that 
$$ \mu:\ \Db_{\on{perv}}(k_X,\C^\times p) \lra \mu Perv_{\gamma(p)}$$ 
is essentially surjective. In order to prove that it is an equivalence of categories, 
we first have to calculate the morphisms in $\mu Perv_{\gamma(p)}$.
\begin{prop}\label{corollary91}\emph{ }
   \begin{itemize}
      \item[(i)]  Let $\me{F},\me{G}\in \mu Perv(\Omega)$ and 
       $p\in\gamma^{-1}(\Omega)$. Let $\tilde{\me{F}}$ and $\tilde{\me{G}}$ be
       two objects of $\on{D}^{\on{b}}_{perv}(k_X,\C^{\times} p)$ such that 
       $\mu\tilde{\me{F}}\simeq \me{F}$ and $\mu\tilde{\me{G}}\simeq \me{G}$ in
      $\mu Perv_{\gamma(p)}$. Then we have
       $$ \on{Hom}_{_{\mu Perv_{\gamma(p)}}}(\me{F},\me{G})
           \simeq \on{H}^0\!(\C^{\times} p,\mu hom
            \big(\tilde{\me{F}},\tilde{\me{G}})\big). $$
     \item[(ii)] Let $\me{F},\me{G}\in\Db_{\on{perv}}(k_X,\C^{\times} p)$. Then
 $$ \on{Hom}_{_{\mu Perv_{_{\gamma(p)}}}}(\mu\me{F},\mu\me{G})
     \simeq \on{Hom}_{_{\Db_{\on{perv}}(k_X,\C^{\times} p)}}(\me{F},\me{G}). $$
     \item[(iii)] The functor $\mu$ induces a canonical equivalence of categories
  $$  \Db_{\on{perv}}(k_X,\C^{\times} p) \overset{\sim}{\lra} 
      \mu Perv_{\gamma(p)}. $$
  \end{itemize}
\end{prop}
\begin{proof}
Note that the functor of (iii) is obviously essentially suejective. Hence (iii)
follows from (ii). Moreover (ii) follows from (i) and Proposition \ref{cutoffprop}.\\ 
Let us prove (i). We have
\begin{align*}
   \on{Hom}_{_{\mu Perv_{\gamma(p)}}}(\me{F},\me{G})
  &\simeq \inlim{\gamma(p)\in V\subset P^*X}
     {\on{Hom}_{\Db(\I(k_{\gamma^{-1}(V)}))}
     (\me{F}|_{_{\gamma^{-1}(V)}},\me{G}|_{_{\gamma^{-1}(V)}})}\\
  &\simeq \inlim{\gamma(p)\in V\subset P^*X}
     {\on{Hom}_{\Db(\I(k_{\gamma^{-1}(V)}))}
     (\mu\tilde{\me{F}}|_{_{\gamma^{-1}(V)}},
     \mu\tilde{\me{G}}|_{_{\gamma^{-1}(V)}})} \\
  & \simeq \inlim{\gamma(p)\in V\subset P^*X}
   {\on{H}^0\!\big(V,\on{RHom}(\mu\tilde{\me{F}},\mu\tilde{\me{G}})\big)}\\
  &\simeq  \inlim{\gamma(p)\in V\subset P^*X}{\on{H}^0\!\Gamma(V,
        \mu hom(\tilde{\me{F}},\tilde{\me{G}})})\\
  &\simeq \on{H}^0\big(\C^{\times} p,\mu hom(\tilde{\me{F}},
       \tilde{\me{G}})\big).
\end{align*}
\end{proof}
Proposition \ref{corollary91} implies that if $\mu Perv$ is a stack then it is 
equivalent to the stack associated to $\Db_{\on{perv}}(k_X,\,*\,)$. In order to prove 
that $\mu Perv$ is a stack we will apply the general results of Section 2.\\
Let us recall the following well-known proposition with proof.
\begin{prop}
  Let $\me{F},\me{G}$ be two perverse sheaves on $X$. Then 
  $\mu hom(\me{F},\me{G})[d_X]$ is a perverse sheaf. In particular it is
  concentrated in positive degrees.
\end{prop}
\begin{proof}
According to \cite{KS3}, Corollary 10.3.20, $\mu hom(\me{F},\me{G})[d_X]$ is a 
perverse sheaf on $T^*X$. Hence for any complex analytic subset
$S$ of $T^*X$ we have 
$$ \on{H}^{j+d_X}_S(\mu hom(\me{F},\me{G}))|_{S}\simeq 0 $$
if $j<-\dim S$. Now recall that
$$ \supp(\mu hom(\me{F},\me{G})) \subset \SS(\me{F})\cup \SS(\me{G}). $$
Since $\me{F},\me{G}$ are perverse sheaves their micro-supports are Lagrangian
subsets of $T^*X$, hence are of dimension $d_X$. Therefore
$$ \on{H}^{j}(\mu hom(\me{F},\me{G}))\simeq 0 $$
for $j<0$.  
\end{proof} 
\begin{prop}
  The prestack $\mu Perv$ of microlocal perverse sheaves is separated.
\end{prop}
\begin{proof}
Let $\me{F},\me{G}$ be two microlocal perverse sheaves of
$\mu Perv(\Omega)$. Recall that $\me{F},\me{G}\in\Db(\I(k_{\gamma^{-1}(\Omega)}))$
and that we have by definition
$$ \me{H}om_{\mu Perv}(\me{F},\me{G}) \simeq
   \gamma_*\me{H}om_{\Db(\I(k_{\,*\,}))}(\me{F},\me{G}). $$
Hence it is sufficent to prove that $\me{H}om_{\Db(\I(k_{\,*\,}))}(\me{F},\me{G})$ is 
a sheaf.\\
We will first show that the complex $\Dr\me{H}om(\me{F},\me{G})$ is 
concentrated in positive degrees. This is a local question, hence we may 
assume that $\me{F}\simeq\mu\me{\tilde{F}}$, $\me{G}\simeq\mu\me{\tilde{G}}$ for
two objects $\me{\tilde F},\me{\tilde{G}}$ of 
$\on{D}^{\on{b}}_{perv}(k_X,\gamma^{-1}(\Omega))$. Since 
$$ \Dr\me{H}om(\mu(\me{\tilde F}),\mu(\me{\tilde G}))\simeq
   \mu hom(\me{\tilde F},\me{\tilde G})$$
we are reduced to study $\mu hom(\me{\tilde F},\me{\tilde G})_p$ for
any $p\in\gamma^{-1}(\Omega)$. By invariance of quantized contact 
transformation (cf. Theorem \ref{perversekerneltheorem}) we may assume that 
$\me{\tilde F}$ and $\me{\tilde G}$ are perverse sheaves on a neighborhood of 
$\pi(p)$. Hence $\mu hom(\me{\tilde F},\me{\tilde G})_p$ is concentrated in
positive degrees by the last proposition.\\
Therefore $\Dr\me{H}om(\me{F},\me{G})$ is concentrated in positive degrees.
Then $\mc{H}om_{\Db(\I(k_{\,*\,}))}(\me{F},\me{G})$ is a sheaf since
$$ \me{H}om_{\Db(\I(k_{\,*\,}))}(\me{F},\me{G})(V)\simeq
   \on{H}^0\!(V,\Dr \me{H}om(\me{F},\me{G}))\simeq
   \Gamma(V,\on{H}^0\Dr\me{H}om(\me{F},\me{G})). $$
\end{proof}
\begin{thm}
   The prestack $\mu Perv$ on $P^*X$ is an abelian stack. Moreover the functor
   $\mu$ induces an equivalence of abelian stacks
   $$ \Db_{\on{perv}}(k_X,\,*\,)^\ddagger \overset{\sim}{\lra}
      \mu Perv, $$
   where $\Db_{\on{perv}}(k_X,\,*\,)^\ddagger$ denotes the stack associated to
   $\Db_{\on{perv}}(k_X,\,*\,)$.
\end{thm}
\begin{proof}
Since microlocal perverse sheaves form a separated subprestack of
the prestack of ind-sheaves and are obviously defined by a local property, 
they form a stack.\\
Moreover, Proposition \eqref{corollary91} states that the functor of prestacks
$$ \mu:\  \Db_{\on{perv}}(k_X,\,*\,) \lra \mu Perv $$
induces equivalences of categories in the stalks. Hence $\mu$ identifies
$\mu Perv$ with the stack associated to $\Db_{\on{perv}}(k_X,\,*\,)$.
\end{proof}

\subsection{Autoduality}

\begin{prop}
  Let $\me{F},\me{G}\in \mu Perv(\Omega)$. Then
  $$ \Dr \me{H}om(\me{F},\me{G})[d_X] $$
  is a perverse sheaf on $\gamma^{-1}(\Omega)$.
\end{prop}
\begin{proof}
Locally we can find $\widetilde{\me{F}}$,$\widetilde{\me{G}}\in \Db_{\Rc}(k_X)$ 
such that
$$ \Dr \me{H}om(\me{F},\me{G})\simeq
   \mu hom(\widetilde{\me{F}},\widetilde{\me{G}}) $$
By invariance of quantized contact transformations we may assume that
$\widetilde{\me{F}},\widetilde{\me{G}}$ are perverse sheaves. Then the result
follows fromthe fact that $\mu hom(\widetilde{\me{F}},\widetilde{\me{G}})[d_X]$ 
is a perverse sheaf.
\end{proof}
\begin{prop}
  There is a contravariant functor 
  $$ D:\ \mu Perv \lra \mu Perv $$
  such that 
  \begin{itemize}
    \item[(i)] $D\circ D\simeq \on{Id}$.
    \item[(ii)] $\Dr\me{H}om(\me{F},\me{G})\simeq \Dr \me{H}om(D\me{G},D\me{F})$.
    \item[(iii)] If $\me{F}\in\mu Perv(\Omega)$ is isomorphic to $\mu\widetilde
       {\me{F}}$ and $\widetilde{\me{F}}\in\Db_{\on{perv}}(k_X,\gamma^{-1}(\Omega))$ 
       then we have a natural isomorphism $D\me{F}\simeq 
       \mu\Dr\me{H}om(\me{F},\omega_X)$.
    \item[(iv)] The stack $\mu Perv$ is autodual, i.e. it is equivalent to its 
      opposite stack.
\end{itemize}
\end{prop}
\begin{proof}
Recall that the functor $\on{D}=\Dr\me{H}om(\,.\,,\omega_X)$ induces a contravariant
equivalence of prestacks
$$ \on{D}:\ \Db_{\on{perv}}(k_X,\,*\,) \lra \Db_{\on{perv}}(k_X,\,*\,).$$
Hence we get a contravariant equivalence $\on{D}$ on the stack associated to 
$\Db_{\on{perv}}(k_X,\,*\,)$ which satisfies by definition (i) and (iii).\\
Let $\me{F},\me{G}\in\mu Perv(\Omega)$. It is enough to prove (ii) locally. Hence we 
may assume that there are objects $\widetilde{\me{F}},\widetilde{\me{G}}\in
\Db_{\on{perv}}(k_X,\gamma^{-1}\Omega)$ such that $\mu\widetilde{\me{F}}\simeq\me{F}$ 
and $\mu\widetilde{\me{G}}\simeq\me{G}$. Recall (cf. Exercise IV.4 of \cite{KS3}) that
on $\gamma^{-1}\Omega$ we have
$$  \mu hom (\widetilde{\me{F}},\widetilde{\me{G}})\simeq
   \mu hom (\on{D}\widetilde{\me{G}},\on{D}\widetilde{\me{F}}). $$
Then on $\gamma^{-1}\Omega$ we get
$$ \Dr \me{H}om(\mu \widetilde{\me{F}},\mu\widetilde{\me{G}})\simeq
   \mu hom (\widetilde{\me{F}},\widetilde{\me{G}})\simeq
   \mu hom (\on{D}\widetilde{\me{G}},\on{D}\widetilde{\me{F}}) \simeq
   \Dr \me{H}om(\mu \on{D}\widetilde{\me{G}},\mu\on{D}\widetilde{\me{F}}) \simeq
   \Dr \me{H}om(\on{D}\me{G},\on{D}\me{F}). $$
\end{proof}

\subsection{Integral transforms for microlocal perverse sheaves}

Let $p\in\Omega_X\subset T^*X$, $q\in\Omega_Y\subset T^*Y$ where $\Omega_X,\Omega_Y$
are $\C^\times$-conic open subsets. Suppose that we are given a contact transformation
$$ \chi:\ \Omega_X\overset{\sim}{\lra} \Omega_Y $$
such that $\chi(p)=q$. We have seen in Section 6 that after shrinking $\Omega_X$ 
and $\Omega_Y$, we can establish an equivalence of prestacks
$$ \Phi_{\me{K}}^{\on{perv}}:\ \chi_*\Db_{\on{perv}}(k_X,\,*\,)|_{\Omega_X} 
   \overset{\sim}{\lra} 
   \Db_{\on{perv}}(Y,\,*\,)|_{\Omega_Y}. $$
Hence we get an equivalence of the associated stacks of microlocal perverse sheaves:
\begin{equation}\label{equationalpha} 
  \chi_*\mu Perv^X|_{\gamma(\Omega_X)} \overset{\sim}{\lra} 
   \mu Perv^Y|_{\gamma(\Omega_Y)} 
\end{equation}
together with a commutative diagram
$$ \xymatrix{
    {\chi_*\Db_{\on{perv}}(k_X,\,*\,)|_{\Omega_X}} 
    \ar[r]^{\Phi_{\me{K}}^{\on{perv}}}_{\sim} \ar[d]  &
     {\Db_{\on{perv}}(Y,\,*\,)|_{\Omega_Y}} \ar[d] \\
    {\chi_*\mu Perv^X|_{\gamma(\Omega_X)}} \ar[r] &
   {\mu Perv^Y|_{\gamma(\Omega_Y)}} } $$
The existence of this diagram is sufficent for many applications. However, in Section 9,
we need to know that the equivalence \eqref{equationalpha} is given by integral 
transform with the kernel $\mu\me{K}$. \\ 
Let us fix some notations first. We will first consider real manifolds $X,Y,Z$ 
and the diagram 
$$ \xymatrix{ & X\times Y\times Z \ar[dl]_{q_{12}}
   \ar[d]^{q_{13}}\ar[dr]^{q_{23}} & \\
      X\times Y & X\times Z & Y\times Z. } $$
Composition (or integral transforms) will be considered in its ind-version.
For $\me{F}\in\Db(\I(k_{X\times Y}))$ and $\me{G}\in\Db(\I(k_{Y\times Z}))$
we set
$$ \me{F}\circ \me{G}=\Dr q_{13!!}(q_{12}^{-1}\me{F}
  \otimes q_{23}^{-1}\me{G}),$$
even if $\me{F},\me{G}$ are actually complexes of classical sheaves. In this case
$\alpha(\me{F}\circ\me{G})$ gives the classical composition. Note that there is 
always a natural morphism
$$ \me{F}\circ \me{G} \lra \alpha(\me{F}\circ\me{G}). $$
Furthermore, we will consider a variant. Consider
$$ q_{12}^a:\,T^*X\times T^*Y\times T^*Z\ra T^*X\times T^*Y \quad ; \quad
   \big((x;\xi),(y;\eta),(z;\zeta)\big)\mapsto \big((x;\xi),(y,\eta)\big). $$
Then set for  $\me{F}\in\Db(\I(k_{T^*X\times T^*Y}))$ and $\me{G}\in\Db(\I(k_{T^*Y
\times T^*Z}))$
$$ \me{F}\overset{a}{\circ}\me{G}=
   \Dr q_{13!!}(q_{12}^{a-1}\me{F}\otimes q_{23}^{-1}\me{G}).$$
Note that the operation $\overset{a}{\circ}$ is associative up to 
natural isomorphism.\\
Recall that the stack $\mu Perv$ is embedded into the prestack 
$\gamma_*\Db(\I(k_{\,*\,}))$.\\
Hence, for any kernel $\me{K}\in\Db(k_{Y\times X})$ the functor
$$ \Phi_{\mu\me{K}}^a=\mu\me{K}\overset{a}{\circ}:
    \ \gamma_*\Db(\I(k_{\,*\,}))|_{\Omega_X} \lra
      \gamma_*\Db(\I(k_{\,*\,}))|_{\Omega_Y} $$
is well defined. Our interest in this section is to find kernels $\me{K}$ such that 
this functor preserves the subprestack $\mu Perv$. For this purpose, we may use the 
fact that microlocal perverse sheaves are defined by a local property. Hence, if 
$\me{F}\in\mu Perv^X$ then $\Phi_{\mu\me{K}}^a(\me{F})\in \mu Perv^Y$ if and only it 
is microlocally perverse on $\C^\times p$ for every $p\in \dot T^*Y$.\\ 
Take an object $\me{F}\in\mu Perv^X$. Locally $\me{F}\simeq \mu\widetilde{\me{F}}$
for some object $\widetilde{\me{F}}\in\Db_{\on{perv}}(k_X,\C^\times p)$. Hence locally
$$ \mu\me{K}\overset{a}{\circ}\me{F} \simeq \mu \me{K}\overset{a}{\circ} 
   \mu\widetilde{\me{F}}. $$
Since $\me{K}\circ\widetilde{\me{F}}\in\Db_{\on{perv}}(Y,\C^\times q)$, if we know that
there is an isomorphism
\begin{equation}\label{miccompform}
 \mu\me{K}\overset{a}{\circ}\mu\widetilde{\me{F}} \simeq
   \mu(\me{K}\circ\widetilde{\me{F}}) 
\end{equation}
in a neighborhood of $\C^\times q$ then we have shown that the object 
$\mu\me{K}\overset{a}{\circ}\me{F}$ lies in $\mu Perv^Y$.\\  
Kashiwara's composition formula states that the morphism \eqref{miccompform} exists
and gives a criterion for it to be an isomorphism. Hence it allows us to a certain 
extent to translate results from the microlocal study of sheaves from \cite{KS3} 
to the study of ind-sheaves and in particular to microlocal perverse sheaves.
Let us recall Kashiwara's theorem.
\begin{thm}[Microlocal composition of kernels]\label{mainthm}\emph{ }\\
  Let $\me{K}_1\in\Db(\I(k_{X\times Y}))$ and $\me{K}_2\in\Db(\I(k_{Y\times Z}))$.
  \begin{itemize}
    \item[(1)] There is a natural morphism
       \begin{equation}\label{compmorph}
          \mu_{X\times Y}\me{K}_1\overset{a}{\circ}
          \mu_{Y\times Z}\me{K}_2 \lra
          \mu_{X\times Z}(\me{K}_1\circ \me{K}_2). 
       \end{equation}
    \item[(2)] Assume the non-characteristic condition\begin{footnote}
        {For two sets $S_1\subset T^*X\times T^*Y$ and $S_2\subset T^*Y\times T^*Z$
        we denote by $S_1\overset{a}{\ctimes{T^*Y}}S_2$ the cartesian product of 
        $q_{2}^a|_{_{S_1}}:\,S_1\ra T^*Y$ and $q_{1}|_{_{S_2}}:\,S_2\ra T^*Y$, hence
       $$S_1\overset{a}{\ctimes{T^*Y}}S_2=\Big\{((x;\xi_x),(y;\xi_Y),(z;\xi_Z)\in
          T^*X\times T^*Y\times T^*Z\ |\ ((x;\xi_x),(y;\xi_Y))\in S_1
          \ \ ((y;-\xi_Y),(z;\xi_Z))\in S_2\Big\}. $$ }\end{footnote}
        $$ \SS_0(\me{K}_1)\overset{a}{\ctimes{T^*Y}}\SS_0(\me{K}_2) \cap
           (T^*_XX\times T^*Y\times T^*_ZZ) \subset
              T^*_XX\times T^*_YY\times T^*_ZZ, $$
        then the morphism 
       $$ \on{K}_{X\times Z}\circ \big(\mu_{X\times Y}\me{K}_1\overset{a}{\circ}
           \mu_{Y\times Z}\me{K}_2\big) \lra
          \mu_{X\times Z}(\me{K}_1\circ \me{K}_2) $$
        is an isomorphism and 
        $$ \mu_{X\times Y}\me{K}_1\overset{a}{\circ}
           \mu_{Y\times Z}\me{K}_2 \lra
          \mu_{X\times Z}(\me{K}_1\circ \me{K}_2) $$
        is an isomorphism outside $\ol{p_{13}\big(\SS_0(\me{K}_1)
         \overset{a}{\ctimes{T^*Y}}\SS_0(\me{K}_2)\cap T^*X\times T^*_YY
         \times T^*Z\big)}$.
  \end{itemize}
\end{thm}
\begin{remark}
  \em The ``ideal situation'' to apply Kashiwara's theorem is given by condition (1)
  or (2) below.
  \begin{itemize}
     \item[(1)] $(T^*_XX\times T^*Y\cup T^*X\times T^*_YY)\cap \SS(\me{K}_1)\subset
                  T^*_XX\times T^*_YY,$
     \item[(2)] $(T^*_YY\times T^*Z\cup T^*Y\times T^*_ZZ)\cap \SS(\me{K}_2)\subset
                  T^*_YY\times T^*_ZZ.$
  \end{itemize}
  Then the natural morphism (\ref{compmorph}) is an isomorphism outside the zero 
  section.\\
  In particular, let $\me{F}\in\Db(\I(k_X))$ and consider $\me{K}$ such that
  $$ \SS(\me{K})\cap\Big(T^*_YY\times T^*X\cup T^*Y\cup T^*_XX\Big)\subset
        T^*_YY\times T^*_XX $$
  Then the morphism (\ref{compmorph})
   $$ \mu(\me{K})\overset{a}{\circ}\mu\me{F}\lra
      \mu(\me{K}\circ \me{F}) $$
  is an isomophism outside the zero section. For instance we get the following 
  lemma. \em
\end{remark}
\begin{lemma}
   Let $\me{F}\in\Db(\I(k_X))$. The morphism
   $$ \mu(k_{\Delta})\overset{a}{\circ}\mu\me{F}\lra
      \mu(k_{\Delta}\circ \me{F})\simeq \mu\me{F} $$
   is an isomophism outside the zero section.
\end{lemma}
We will mostly be interested in the following more general situation:
\begin{defi}
  Let us denote by $\widetilde{\on{N}}(Y,X,\Omega_Y,\Omega_X)$ the full subcategory of 
  $\Db(k_{Y\times X}, \Omega_Y\times T^*X)$ such that 
  \begin{itemize}
    \item[(i)] $\SS(\me{K})\cap (\Omega_Y\times T^*X\cup T^*Y\times \Omega_X^a)\subset
          \Omega_Y\times \Omega_X^a$
    \item[(ii)] $p_1:\,\SS(\me{K})\cap \Omega_Y\times T^*X\ra \Omega_Y$ is proper.
  \end{itemize}
\end{defi}
Note that $\widetilde{\on{N}}(Y,X,\Omega_Y,\Omega_X)\subset 
\on{N}(Y,X,\Omega_Y,\Omega_X)$\begin{footnote}{The category
$\on{N}(Y,X,\Omega_Y,\Omega_X)$ has been recalled at the beginning of Section 5.1.}
\end{footnote} and  that the kernels produced in Theorem \ref{kerneltheorem} are objects 
of $\widetilde{\on{N}}(Y,X,\Omega_Y,\Omega_X)$. We get
\begin{prop}\label{isocondition}
  Suppose that $\Omega_X\cap T^*_XX=\varnothing$ and $\Omega_Y\cap T^*_YY=
  \varnothing$.\\
  Let $\me{K}\in\widetilde{\on{N}}(Y,X,\Omega_Y,\Omega_X)$ and $\me{F}\in\Db(k_X)$ 
  such that $\SS(\me{F})\cap \dot T^*X\subset \Omega_X$.
  Then the natural morphism 
  $$ \mu\me{K}\overset{a}{\circ}\mu\me{F} \lra \mu(\me{K}\circ\me{F}) $$
  is an isomorphism on $\Omega_Y$.
\end{prop}
\begin{proof}
We have by hypothesis
\begin{align*}
   \SS(\me{K})\cap (T^*Y\times \SS(\me{F})^a)&\subset 
   (\SS(\me{K}) \cap T^*Y\times \Omega_X^a) \cup (\SS(\me{K})\cap T^*Y\times T^*_XX)\\
   &\subset  \Omega_Y\times \Omega_X^a  \cup T^*Y\times T^*_XX. 
\end{align*}
Intersecting both sides with $T^*_YY\times T^*X$ we get
$$ \SS(\me{K})\cap (T^*Y\times \SS(\me{F})^a)\cap (T^*_YY\times T^*X) \subset
    T^*_XX\times T^*_YY. $$
Hence the non-characterstic condition is satisfied.\\
Moreover by assumption $\SS(\me{K}) \cap (\Omega_Y\times T^*X)\subset \Omega_Y\times 
\Omega_X^a$, hence
$$ \SS(\me{K})\cap (T^*Y\times \SS(\me{F})^a)\cap \Omega_Y\times T^*_XX=\varnothing $$
and 
$$ \SS(\me{K})\cap (T^*Y\times \SS(\me{F})^a)\cap T^*Y\times T^*_XX \subset
    \complement \Omega_Y\times T^*_XX $$
Therefore the morphism is an isomorphism outside $\complement \Omega_Y$, hence on 
$\Omega_Y$.
\end{proof}
\begin{lemma}
  Consider $\me{K}\in\on{N}(Y,X,\Omega_Y,\Omega_X)$ and $\me{F}\in\Db(k_X)$. Then we 
  have a chain of natural isomorphisms: 
  $$(\mu{\me{K}})|_{\Omega_Y\times \Omega_X^a}\overset{a}{\circ}(\mu\me{F})|_{\Omega_X}
     \simeq
     (\mu{\me{K}})^{(1,a)}|_{\Omega_Y\times \Omega_X}\circ(\mu\me{F})|_{\Omega_X}\simeq
       \big(\mu\me{K}\overset{a}{\circ}\mu\me{F}\big)|_{\Omega_Y}.$$
\end{lemma}
\begin{proof}
We will show that the three terms are isomorphic to
$$ (\mu\me{K})|_{\Omega_Y\times T^*X}\overset{a}{\circ}\me{F}. $$
For the right hand side this is trivial.\\
For the two terms on the left hand side denote by 
$j:\,\Omega_Y\times \Omega_X^a\hookrightarrow \Omega_Y\times X$ the inclusion map. Then we have  
$$ (\mu\me{K})|_{\Omega_Y\times T^*X}\simeq
  \Dr j_{!!}j^{-1}\mu\me{K}|_{\Omega_Y\times T^*X} $$
by the hypothesis on $\SS(\me{K})=\supp(\mu\me{K})$ and the result follows easily.
\end{proof}
\begin{remark}
  \em Let $\me{K}\in\widetilde{\on{N}}(Y,X,\Omega_Y,\Omega_X)$. Consider the morphism
  \begin{equation} \label{moralpha}
     (\mu\me{K}\overset{a}{\circ}\mu\me{F})|_{\Omega_Y}\lra \mu(\me{K}\circ
      \me{F})|_{\Omega_Y}\lra \mu(\alpha(\me{K}\circ\me{F}))|_{\Omega_Y}  
  \end{equation}
  The term on the right only depends on the image of $\me{F}$ in $\Db(k_X,\Omega_X)$. By
  the last lemma, the same result holds for the term on the left side\begin{footnote}{
  The assumptions on $\me{K}$ might imply that the second morphism of the 
  composition is an isomorphism, but I cannot prove it.}\end{footnote}. In the 
  following we will be interested in situations when the composition 
  \eqref{moralpha} is an isomorphism.\\ 
  Hence, in the situation of Proposition \ref{isocondition}, we may replace $\me{F}$ 
  by any object isomorphic to $\me{F}$ in $\Db(k_X,\Omega_X)$. However, in general, it 
  is not possible to find an object $\me{F}'$ isomorphic to $\me{F}$ in 
  $\Db(k_X,\Omega_X)$ such that $\SS(\me{F}')\subset \Omega_X$. Therefore we will have 
  to assume a stability condition on 
  $\me{K}$ and the existence of a suitable cut-off functor.\em
\end{remark}
\begin{prop}\label{ingosprop}
 Let $X,Y$ be affine and $\me{K}\in\widetilde{\on{N}}(Y,X,\Omega_Y,\Omega_X)$.\\
 Let $\delta\subset X^*$ be an open cone and $U\subset X$ be a relatively 
 compact open set such that
  $$ (\ol{U}\times \ol{\delta})\cap \dot T^*X \subset \Omega_X.$$
 Set $\Omega_X'=U\times\delta$ and suppose that there exists an open subset 
 $\Omega_Y'\subset \Omega_Y$ such that 
 $\me{K}\in\widetilde{\on{N}}(Y,X,\Omega_Y',\Omega_X')$.\\
 Then the natural morphism (\ref{moralpha})
 $$ \mu\me{K}\overset{a}{\circ}\mu\me{F} \lra \mu(\alpha(\me{K}\circ\me{F})) $$
 is an isomorphism on $\Omega_Y'$. 
\end{prop}
\begin{proof}
There is a cut-off functor $\Phi_{U,\delta}$ such that $\SS(\Phi_{U,\delta}(\me{F}))
\subset \ol{U}\times \ol{\delta}$ and $\Phi_{U,\delta}(\me{F})\ra \me{F}$ is an 
ismorphism in $\Db(k_X,U\times\delta)$.\\
Then the morphism
$$ \mu\me{K}\overset{a}{\circ}\mu(\Phi_{U,\delta}(\me{F})) \lra
   \mu (\me{K}\circ\Phi_{U,\delta}(\me{F})) $$
is an isomorphism on $\Omega_Y'$. Moreover $\supp(\Phi_{U,\gamma}(\me{F}))
\subset \ol{U}$ which is compact. Hence $\me{K}\circ\Phi_{U,\gamma}(\me{F})\ra
\alpha(\me{K}\circ\Phi_{U,\gamma}(\me{F}))$ is an isomorphism. 
Therefore 
$$ \mu\me{K}\overset{a}{\circ}\mu(\Phi_{U,\delta}(\me{F})) \lra
   \mu (\alpha(\me{K}\circ\Phi_{U,\delta}(\me{F}))) $$
is an isomorphism on $\Omega_Y'$ and is isomorphic to  
 $$ \mu\me{K}\overset{a}{\circ}\mu\me{F} \lra \mu(\alpha(\me{K}\circ\me{F})). $$
\end{proof}
Summarizing the situation, we have the following diagram:
$$ \xymatrix{
       {\Db(k_X,\Omega_X)} \ar[r]^{\mu} \ar[d]_{\Phi_{\me{K}}} & 
       {\Db(\I(k_{\Omega_X}))} \ar[d]^{\mu(\me{K})|_{\Omega_Y\times \Omega_X^a}
        \overset{a}{\circ}} \df[dl] \\
       {\Db(k_Y,\Omega_Y)} \ar[r]_{\mu} & {\Db(\I(k_{\Omega_Y}))}     } $$
where the diagonal arrow indicates that the diagram is not commutative, but 
there is a natural transformation which for any $\me{F}\in\Db(k_X,\Omega_X)$ is 
induced by
$$ \mu(\me{K})\overset{a}{\circ}\mu(\me{F}) \lra \mu(\me{K}\circ \me{F})
   \lra \mu(\Phi_{\me{K}}(\me{F})). $$
It is an isomorphism under suitable conditions on $\me{K}$. In particular, consider
the situation of Theorem \ref{kerneltheorem}, we get a diagram of functors of prestacks
$$ \xymatrix{ 
       \chi_*\Db(k_X,\,*)|_{\Omega_X} \ar[r]^\mu \ar[d]_{\Phi_{\me{K}}} &
       \Db(\I(k_{\,*\,}))|_{\Omega_X} \ar[d]^{\mu\me{K}\overset{a}{\circ}}  \df[dl] \\
       \Db(k_Y,\,*)|_{\Omega_Y} \ar[r]^\mu  &
       \Db(\I(k_{\,*\,}))|_{\Omega_Y} } $$ 
that commutes on small open sets of $\chi(U\times\gamma)$\begin{footnote}{To apply 
Proposition \ref{ingosprop}, we actually have to consider an $\R^+$-homogenous 
variant of Theorem \ref{kerneltheorem} where we replace the open sets by $\R^+$-conic 
open sets. The theorem still holds since there is actually no difference between 
$\Db(k_X,S)$ and $\Db(k_X,\R^+S)$ by the $\R^+$-conicity of the micro-support. Hence
the diagram should be considered on the sphere bundle rather then on $\dot T^*X$.}
\end{footnote}. In particular the diagram commutes (up to natural isomorphism)
in the stalks.\\
Now let us consider the complex contact transformations and place ourselves
in the $\C^\times$-conic situation of Theorem \ref{homogenouskerneltheorem}. 
Hence all manifolds are complex, the open sets $\C^\times$-conic, and the kernel
$\me{K}$ is has $\R$-constructible cohomology sheaves. We get a diagram of functors 
of prestacks
$$ \xymatrix{ 
       \gamma_*\Db(k_X,\,*)|_{\Omega_X} \ar[r]^\mu \ar[d]_{\Phi_{\me{K}}} &
       \gamma_*\Db(\I(k_{\,*\,}))|_{\Omega_X} \ar[d]^{\mu\me{K}\overset{a}{\circ}}  \\
       \gamma_*\Db(k_Y,\,*)|_{\Omega_Y} \ar[r]^\mu  &
       \gamma_*\Db(\I(k_{\,*\,}))|_{\Omega_Y} } $$ 
that commutes up to natural isomorphism in the stalks. Hence let us look at this 
diagram in the stalks, restrict to microlocal perverse sheaves and immediately get 
that for $p\in\Omega_X$ and $q=\chi(p)\in\Omega_Y$ the diagram 
\begin{equation}\label{littlediagram} \xymatrix{
        \Db_{\on{perv}}(k_X,\C^\times p)  \ar[r]^\mu    
        \ar[d]_{\Phi_{\me{K}}} & \mu Perv_{\gamma(p)} 
         \ar[d]^{\mu\me{K}\overset{a}{\circ}} \\
        \Db(k_Y,\C^\times q) \ar[r]^\mu & \Db(\I(k_{*}))_{\gamma(q)}  }
\end{equation}
is well defined and commutes up to natural isomorphism.
\begin{thm}\label{thm123}
  The functor of prestacks
  $$ \mu\me{K}\overset{a}{\circ}:\ \gamma_*\Db(\I(k_{\,*\,}))|_{\gamma(\Omega_X)} \lra
       \gamma_*\Db(\I(k_{\,*\,}))|_{\gamma(\Omega_Y)} $$
  induces an equivalence of stacks
 $$ \mu Perv^X|_{\gamma(\Omega_X)} \overset{\sim}{\lra} 
    \mu Perv^Y|_{\gamma(\Omega_Y)} $$
  such that the diagram
  $$ \xymatrix{ 
       \Db_{\on{perv}}(k_X,\,*)|_{\Omega_X} \ar[r]^\mu \ar[d]_{\Phi_{\me{K}}} &
       \mu Perv^X|_{\Omega_X} \ar[d]^{\mu\me{K}\overset{a}{\circ}}  \\
       \Db(k_Y,\,*)|_{\Omega_Y} \ar[r]^\mu  &
       \mu Perv^Y|_{\Omega_Y} } $$ 
  commutes up to natural isomorphism and induces in the stalks the following
  diagram:
   $$ \xymatrix{ 
       \Db_{\on{perv}}(k_X,\C^\times p) \ar[r]^\mu \ar[d]_{\Phi_{\me{K}}^{\on{perv}}} &
       \mu Perv^X_{\gamma(p)} \ar[d]^{\mu\me{K}\overset{a}{\circ}}  \\
       \Db_{\on{perv}}(k_Y,\C^\times q) \ar[r]^\mu  &
       \mu Perv^Y_{\gamma(q)}. } $$ 
\end{thm}
\begin{proof}
Consider the diagram \eqref{littlediagram}. Let 
$\me{F}\in\Db_{\on{perv}}(k_X,\C^\times p)$. Recall that for a sufficently small
relatively compact open neighborhood $V$ of $\pi(p)$ the object 
$\me{K}\circ\me{F}_{\ol{V}}$ is a well-defined object of 
$\Db_{\on{perv}}(k_Y,\C^\times q)$ and isomorphic to $\me{K}\circ\me{F}$ in 
$\Db(k_Y,\C^\times q)$. Since $\Db_{\on{perv}}(k_Y,\C^\times q)$ is a full subcategory
of $\Db(k_Y,\C^\times q)$, we get that $\mu\me{K}\overset{a}{\circ}\mu\me{F}$
is an object of $\mu Perv^Y_{\gamma(q)}$. Hence we get the induced functor of stacks
 $$ \mu Perv^X|_{\gamma(\Omega_X)} \lra \mu Perv^Y|_{\gamma(\Omega_Y)}. $$
In the same way we verify that $\mu\me{K}^*\overset{a}{\circ}$ is a well defined
functor of stacks in the opposite direction and it is easily verified that
$\mu\me{K}^*\overset{a}{\circ}$ is a quasi-inverse\begin{footnote}{
Note that the difficulties encountered in Theorem \ref{kerneltheoremrcons} to
prove the equivalence which were related to the fact that $\Db_{Rc}(k_X,*)$ is not a 
full subprestack of $\Db(k_X,*)$ do not appear here since all the morphisms are 
well-defined in $\Db(\I(k_*))$ and $\mu Perv$ is a full substack of 
$\gamma_*\Db(\I(k_*))$.}\end{footnote}.
\end{proof} 
Now suppose that we have contact transformations
$$ \chi_1:\ \Omega_X\overset{\sim}{\ra} \Omega_Y \qquad \qquad
   \chi_2:\ \Omega_Y\overset{\sim}{\ra} \Omega_Z $$
and kernels $\me{K}_1\in\Db_{\Rc}(\I(k_{Y\times X}))$, 
$\me{K}_2\in\Db_{\Rc}(\I(k_{Z\times Y})$ such that Theorem \ref{kerneltheorem} is valid. 
Then we know that $\Phi^a_{\mu\me{K}_1}=\mu\me{K}_1\overset{a}{\circ}$ and 
$\Phi^a_{\mu\me{K}_2}=\mu\me{K}_2\overset{a}{\circ}$ are well-defined on microlocally
perverse sheaves. However note that $\me{K}_2\circ\me{K}_1$ is not 
$\R$-constructible and we do not know if
$$ \mu\me{K}_2\overset{a}{\circ}\mu\me{K}_1 \lra \mu(\me{K}_2\circ\me{K}_1) $$
is ani somorphism in $\Db(\I(k_{T^*Z\times T^*X}))$.
Nevertheless we get
\begin{prop}\label{prop826}
  The functor
 $$ \mu(\me{K}_2\circ\me{K}_1)\overset{a}{\circ}:\ \mu Perv^X \lra \mu Perv^Z $$
  is well-defined and naturally isomorphic to $\Phi^a_{\mu\me{K}_1}\circ
  \Phi^a_{\mu\me{K}_1}$.
\end{prop}
\begin{proof}
The strategy is similar to the proof of Theorem \ref{thm123}. First we can show that
$\mu(\me{K}_2\circ\me{K}_1)\overset{a}{\circ}:\ \mu Perv^X \lra \mu Perv^Z$ is 
well-defined by looking in the stalks and using the fact that 
$\Db_{\on{perv}}(Y,\C^\times q)$ is a full subcategory of $\Db(k_Y,\C^\times q)$. The
same argument shows the isomorphism.
\end{proof}
\begin{remark}
  \emph{Proposition \ref{prop826} is an important step towards the definition of 
  the stack of 
  microlocal perverse sheaves on a complex contact manifold. Such a manifold is 
  locally isomorphic to an open subset $\Omega_X\subset P^*X$ and the transition maps 
  are contact transformations. Using the proposition we can locally associate an 
  equivalence of stacks. But the choice of the kernel $\me{K}$ is neither unique nor 
  canonical.} 
\end{remark} 

\section{Microlocal Riemann-Hilbert correspondence}

The classical Riemann-Hilbert correspondance states that on a complex variety 
$X$ the solution functor $\Dr\me{H}om_{\me{D}_X}(\,\cdot\,,\me{O}_X)$ defines 
an equivalence between the stack of regular holonomic $\me{D}_X$-modules and the stack 
of perverse sheaves. A quasi-inverse of this functor has been constructed 
explicitly by Kashiwara (cf. \cite{K3}).\\ 
We denote by $\mc{H}ol\mc{R}eg(\me{D}_X)$ the abelian category of regular 
holonomic $\me{D}_X$-modules.\\
Let $\me{M}$ be a holonomic $\me{D}$-module. We set
$$ \on{Sol}(\me{M})=\Dr \me{H}om_{\me{D}_X}(\me{M},\me{O}_X). $$
Let $\me{F}$ be an object of $\Db_{\Rc}(k_X)$ and set
$$ \on{RH}(\me{F})=\thom(\me{F},\me{O}_X). $$
Then the Riemann-Hilbert corresponance can be stated as follows.
\begin{thm}
  The functors $\on{RH}$ and $\on{Sol}$ define quasi-inverse equivalences of
  abelian stacks
   $$ \xymatrix@C=3cm{ {\on{Perv}(k_X)} \ar@<+3pt>[r]^{\on{RH}} & 
       {\mc{H}ol\mc{R}eg(\me{D}_X).} 
        \ar@<+3pt>[l]^{\on{Sol}} } $$
\end{thm} 

The microlocal Riemann-Hilbert correspondance should therefore establish
an equivalence between the stack of microlocal perverse sheaves and
the stack of regular holonomic $\me{E}_X$-modules.

\subsection{The ind-objects $\mu\me{O}_X$ and $\mu\me{O}^t_X$}

The microlocalization of the ring of holomorphic functions defines an object
of $\I(\C_X)$. More precisely, we have
\begin{prop}
  The ind-sheaf
  $$ \mu\me{O}_X|_{\dot T^*X}$$
  is concentrated in degree $d_X$. 
\end{prop}
Let $\me{F}\in\Db_{\Rc}(\C_X)$. We will study the microlocal solution complex
$$ \mu hom(\me{F},\me{O}_X)\simeq \Dr \me{H}om(\pi^{-1}\me{F},\mu\me{O}_X) $$
The stalks of this complex have been studied in \cite{KS3} and we will show that
some results can be extended to open neighborhoods.\\
The ``ring'' $\me{O}^t_X\in \Db(\I(\C_X))$ of temperate holomorphic functions has 
been defined in \cite{KS2}. It is defined from the ring 
$\me{D}b^t_X\in \Db(\I(\C_X))$ as
$$ \me{O}^t_X=\Dr\me{H}om_{\beta(\me{D}_{\ol{X}})}(\beta(\me{O}_{\ol{X}}),
   \me{D}b^t_X). $$
We will not recall the construction of $\me{D}b^t_X$ here. Recall that $\me{O}^t_X$
is only defined in the derived category. It is not concentrated in a single degree.\\
The link with Kashiwara's functor $T\me{H}om$ is given by the formula
\begin{equation}\label{temperedhom} 
  \Dr\me{H}om(\me{F},\me{O}^t_X)=T\me{H}om(\me{F},\me{O}_X) 
\end{equation}
where $\me{F}\in\Db_{\Rc}(\C_X)$.\\
By microlization we get an object $\mu\me{O}^t_X\in\Db(\I(\C_{T^*X}))$. It is not known
if this object is concentrated in a single degree or not.\\
In \cite{KS2}, the full functoriality of $\me{O}^t$ is established. We will 
only need the following result:
\begin{prop}
  Let $f:\,X\ra Y$ be a smooth map between complex varieties. Then there is a
  canonical isomorphism in $\Db(\I(\beta(\me{D}_X)))$:   
  $$  \Dr\me{IH}om_{\beta\me{D}_X}(\beta\me{D}_{X\ra Y},\me{O}^t_X)
      \simeq f^{-1}\me{O}^t_Y. $$
\end{prop}
First we will prove that the microlocalization of the formula \eqref{temperedhom} 
holds, i.e.
\begin{prop}\label{temperedsolutions}
  Let $\me{F}\in\Db_{\Rc}(\C_X)$. Then
  $$ \Dr\me{H}om(\mu\me{F},\mu\me{O}^t)=t\mu hom(\me{F},\me{O}_X). $$
\end{prop}
\begin{proof}
The proof is similar to the proof of the formula
$$ \Dr \me{H}om(\mu\me{F},\mu\me{G})\simeq \mu hom(\me{F},\me{G}) $$
from \cite{K4} where $\me{F},\me{G}\in\Db(k_X)$.\\
The normal deformation of the diagonal in $X\times X$ can be visualized by the 
following diagram
$$\xymatrix{{TX}\ar[r]^(.3){\sim} & {T_{\Delta_X}(X\times X)}
   \ar[d]^{\tau_X} \ar@<-1pt>@{ (->}[r]^(.6){s} & 
    {\widetilde{X\times X}}\ar[d]_p\ar@<-2pt>@/_3.5ex/[dd]_(.7){p_1}    
   \ar@<+2pt>@/^3.5ex/[dd]^(.7){p_2} & {\Omega} 
   \ar@<1pt>@{ )->}[l]_(.3){j}\ar[ld]^(.4){\widetilde p} \\   
   & {X}\ar@<-1pt>@{ (->}[r]_{\Delta_X}\ar[dr]_{id_{_X}} & 
   {X\times X} \ar@<+2pt>[d]^(.4){q_2}\ar@<-2pt>[d]_(.4){q_1} & {} \\
    & {} & {X} & {.} }$$
Note that $\tilde p,p_1$ and $p_2$ are smooth but $p$ is not. 
Also, the square is not cartesian. We use the same notations as in
Section 5 when we deform the diagonal in $T^*X\times T^*X$. We will need the 
following lemma
\begin{lemma}
  Let $\me{F}\in\Db_{\Rc}(\C_{T^*X})$ and $\me{G}\in\Db(I(\C_{T^*X}))$. Then
  $$ \Dr \me{H}om(\me{F},\mu \me{G})\simeq
  \omega_{X}^{-1}\Big(s^{-1}\Dr \me{H}om((p_2^{-1}\me{F})_{\Omega},
      p^{-1}q_1^!\me{G})\Big)^{\wedge}.$$
\end{lemma}
\begin{proof}
Set
$$ P=\big\{(x;v)\in TT^*X\ |\ \langle v,\omega(x) \rangle \leqs 0\big\}, $$
$$ P'=\big\{(x;v)\in TT^*X\ |\ \langle v,\omega(x) \rangle \geqs 0\big\}. $$
We have
\begin{align*}
  \Dr\me{H}om(\me{F},\mu\me{G})&\simeq
     \Dr\me{H}om(\me{F},\Dr q_{1!!}(q_2^{-1}\me{G}\otimes 
      \Dr p_{!!}(\C_{\ol{\Omega}})
     \otimes \beta(\C_P\otimes \Dr s_*\omega_{TT^*X|T^*X}) \\
   &\simeq \Dr \me{H}om\me{F},\Dr p_{1!!}(p_2^{-1}\me{G}\otimes 
       \C_{\ol{\Omega}}
     \otimes \beta(\C_P\otimes \Dr s_*\omega_{TT^*X|T^*X}))\\
   &\simeq \Dr \me{H}om\me{F},\Dr p_{2!!}(p_1^{-1}\me{G}\otimes 
       \C_{\ol{\Omega}}
     \otimes \beta(\C_{P'}\otimes \Dr s_*\omega_{TT^*X|T^*X}))\\
  &\simeq \Dr p_{2!}\Dr \me{H}om(p_2^{-1}\me{F},p_1^{-1}\me{G}\otimes 
       \C_{\ol{\Omega}}
     \otimes \beta(\C_{P'}\otimes \Dr s_*\omega_{TT^*X|T^*X}))\\
  &\simeq \Dr p_{2!}\Big(\Dr \me{H}om(p_2^{-1}\me{F},p_1^{-1}\me{G}\otimes 
     \C_{\ol{\Omega}})
     \otimes \C_{P'}\otimes \Dr s_*\omega_{TT^*X|T^*X}) \Big)
\end{align*}
Now note that locally on $\widetilde{T^*X\times T^*X}$ the set $\Omega$ is 
convex and 
$$ \SS_0(p_1^{-1}\me{G})
\cap T^*_{\Omega}(\widetilde{T^*X\times T^*X})
  \subset T^*_{\widetilde{T^*X\times T^*X}}(\widetilde{T^*X\times T^*X}).$$
Hence
$$ \C_{\ol{\Omega}}\otimes p_1^{-1}\me{G}\simeq
   \Dr \me{H}om(\me{C}_{\Omega},\C_X)\otimes p_1^{-1}\me{G}\simeq
   \Dr \me{IH}om(\C_{\Omega},p_1^{-1}\me{G}). $$
Moreover note that 
$$ \omega_{TT^*X|T^*X}\simeq \tau_{T^*X}^{-1}
   \omega_{\Delta|T^*X\times T^*X}^{\otimes -1}
   \simeq s^{-1}\tilde p^{-1} q_1^!\C_{T^*X}. $$
Hence
\begin{align*}
  \Dr\me{H}om(\me{F},\mu\me{G})&\simeq
  \Dr p_{2!}\Big(\Dr \me{H}om(p_2^{-1}\me{F}\otimes \C_{\Omega},p_1^{-1}\me{G})
     \otimes \C_{P'}\otimes \Dr s_*\omega_{TT^*X|T^*X}) \Big)\\
  &\simeq \Dr p_{2!}\Dr s_!\Big( s^{-1} \Dr\me{H}om(p_2^{-1}\me{F}
    \otimes \C_{\Omega},
     p^{-1}q_1^!\me{G})\otimes \C_{P'}\Big)\\
  &\simeq \Dr \tau_!\Big( s^{-1} \Dr\me{H}om(p_2^{-1}\me{F}\otimes \C_{\Omega},
     p^{-1}q_1^!\me{G})\otimes \C_{P'}\Big)\\
  &\simeq \omega^{-1}_X \Big(s^{-1}\Dr\me{H}om((p_2^{-1}\me{F})_{\Omega},
     p^{-1}q_1^!\me{G})\Big)^{\wedge}. 
\end{align*}
This shows the lemma.
\end{proof}
Now consider the situation of the lemma but take $\me{F}\in\Db(\C_X)$. Then a standard 
calculation shows that 
$$ \Dr\me{H}om(\pi^{-1}\me{F},\mu\me{G})=
      \Big(s^{-1}\Dr\me{H}om\big((p_2^{-1}\me{F})_{\Omega},
      p^{-1}q_1^!\me{G}\big)\Big)^{\wedge}. $$
When we apply this result to $\me{O}^t_X$ and $\me{F}\in\Db_{\Rc}(k_X)$ we get
\begin{align*}
   \Dr\me{H}om&(\pi^{-1}\me{F},\mu\me{O}^t)\simeq
      \Big(s^{-1}\Dr\me{H}om\big((p_2^{-1}\me{F})_{\Omega_X},
       p^{-1}q_1^!\me{O}^t\big)\Big)^{\wedge} \\
     &\simeq \Big(s^{-1}\Dr\me{H}om_{\me{D}_{\widetilde{X\times X}}}
        \big(\me{D}_{\widetilde{X\times X}\overset{p_1}{\ra} X},
        \thom((p_2^{-1}\me{F})_{\Omega_X}, 
           \me{O}_{\widetilde{X\times X}}\big)\otimes s^{-1}p^{-1}q_1^!\C_X
          \Big)^{\wedge} \\
     &\simeq \Big(s^{-1}\me{D}_{X\overset{p_1}{\leftarrow}
     \widetilde{X\times X}}
         \underset{\me{D}_{\widetilde{X\times X}}}{\otimes} 
         \thom((p_2^{-1}\me{F})_{\Omega}, 
        \me{O}_{\widetilde{X\times X}}\big)[1]\Big)^{\wedge} \\
    &\simeq t\mu hom(\me{F},\me{O}_X).
\end{align*}
\end{proof}
Following \cite{KS1}, we set for a locally free $\me{O}_X$-module $\me{L}$:
$$ \me{L}^t=\me{O}^t\overset{L}{\underset{\beta\me{O}_X}{\otimes}}\beta\me{L}. $$
\begin{lemma}\label{stupidlemma}
   Let $\me{L}$ a locally free $\me{O}_X$-module of finite rank. Then there are natural 
   isomorphisms
  $$ \mu\me{L} \lra \mu\me{O}_X\overset{\on{L}}{\underset{\pi^{-1}\beta
     \me{O}_X}{\otimes}}\pi^{-1}\beta\me{L}, $$
  $$ \mu\me{L}^t \lra \mu\me{O}^t_X\overset{\on{L}}{\underset{\pi^{-1}\beta
      \me{O}_X}{\otimes}}\pi^{-1}\beta\me{L}. $$
\end{lemma}
\begin{proof}
Let us show the second isomorphism, the proof of the first being similar
since $\me{L}\simeq \beta\me{L}\underset{\beta\me{O}_X}{\otimes}\me{O}_X$.
By definition
$$ \mu\big(\me{O}^t_X\underset{\beta\me{O}_X}{\otimes}\beta\me{L}\big)
    \simeq \Dr q_{1!!}\Big(\on{K}_X\otimes q_2^{-1}\pi^{-1}
   \big(\me{O}^t_X\underset{\beta\me{O}_X}{\otimes}\beta\me{L}\big)\Big). $$
But since $\on{K}_X\simeq \on{K}_X\otimes \beta\C_{\Delta}$ and
$$ \beta\C_{\Delta}\otimes \big(q_2^{-1}\pi^{-1}\me{O}^t_X \underset{\beta q_2^{-1}
   \pi^{-1}\me{O}_X}{\overset{\on{L}}{\otimes}}\beta q_2^{-1}\pi^{-1}\me{L}\big)\simeq
   \big((\beta\C_{\Delta}\otimes q_2^{-1}\pi^{-1}\me{O}^t_X) \underset{\beta q_1^{-1}
   \pi^{-1}\me{O}_X}{\overset{\on{L}}{\otimes}}\beta q_1^{-1}\pi^{-1}\me{L}\big) $$
we get the second isomorphism.
\end{proof}

\subsection{$\me{E}_X$-modules}

The ring $\me{E}_X$ of microdifferential operators on $T^*X$ has been defined
in \cite{SKK}. For a short introduction to the theory of $\me{E}_X$-modules 
we refer to \cite{K5}, for a more detailed study see \cite{S}. The ring 
$\me{E}_X$ has many ``good'' properties, for instance it is coherent and Noetherian.
For our purpose it is convenient to consider microdifferential 
operators outside the zero-section, hence when we write $\me{E}_X$, we 
consider $\me{E}_X|_{_{\dot T^*X}}$.

We will consider two variants of this ring. In the sequel we will identify $T^*X$ 
with $T^*_{\Delta_X}(X\times X)$  by the map
$$ \delta^a:\ T^*X \overset{\sim}{\lra} T^*_{\Delta_X}(X\times X) \hookrightarrow 
   T^*X\times T^*X \quad ;\quad (x;\xi) \mapsto ((x;\xi),(x;-\xi)).$$
The functor $\delta^{a-1}$ is often omitted for complexes with support on $T^*X$.\\
First one defines the ring $\me{E}_X^{\R}$ on $T^*X$ as
$$ \me{E}_X^{\R}\simeq \on{H}^{d_X}\mu hom(\C_{\Delta_X},
   \me{O}_{X\times X}^{(d_X,0)})\simeq \on{H}^{d_X}\Dr\me{H}om(\pi^{-1}
   \C_{\Delta_X},\mu\me{O}_X^{(d_X,0)}). $$
In \cite{An2} Andronikof introduced the ring $\me{E}_{X}^{\R,f}$ on $T^*X$ of tempered
microdifferential operators as
$$ \me{E}_X^{\R,f}=\on{H}^{d_X}\big(t\mu hom(\C_{\Delta},
   \me{O}_{X\times X})\underset{\me{O}_{X\times X}}{\overset{\on{L}}{\otimes}}
   \me{O}^{(d_X,0)}_{X\times X}\big)$$ 
where we omit the functor $(\pi\times\pi)^{-1}$ in the notation.
\begin{prop}
We have 
$$ \me{E}_X^{\R,f}\simeq \on{H}^{d_X}\Dr\me{H}om(\pi^{-1}\C_{\Delta_X},
     \mu\me{O}^{t,(d_X,0)}_{X\times X}).  $$
\end{prop}
\begin{proof}
Follows from Proposition \ref{temperedsolutions} and Lemma \ref{stupidlemma}.
\end{proof}   
Recall some basic properties:
\begin{itemize}
  \item[(i)]  We have $\me{E}_X\simeq \gamma^{-1}\Dr \gamma_*\me{E}_X^{\R}\simeq
    \gamma^{-1}\Dr \gamma_*\me{E}_{X}^{\R,f}$. In particular 
    $\gamma^{-1}\Dr \gamma_*\me{E}_X\simeq \me{E}_X$. Moreover
    $\Dr^i \gamma_*\me{E}_X\simeq 0$ for $i\neq 0$.
  \item[(ii)] The rings $\me{E}^{\R}_X$ and $\me{E}_X^{\R,f}$ are faithfully flat
    over $\me{E}_X$.
  \item[(iii)] The ring $\me{E}_X$ (and therefore $\me{E}_X^{\R},\me{E}_X^{\R,f}$)
    is a $\pi^{-1}\me{D}_X$-modules.
\end{itemize}

A priori, $\me{E}_X$-modules are defined on the cotangent space $T^*X$. 
But by (i) coherent $\me{E}_X$-modules (hence in particular regular holonomic 
$\me{E}_X$-modules) are conic objects, hence it is often convenient to work 
on the projective bundle $P^*X$ or on $\C^\times$-conic sets.\\
Let $\me{M}$ be an $\me{E}_X$-module. Then its support $\supp(\me{M})$ is called 
its characteristic variety. If $\me{M}$ is a $\me{D}_X$-module, then the characteristic 
variety of $\me{E}_X\otimes_{\pi^{-1}\me{D}_X}\pi^{-1}\me{M}$ coincides with the 
characteristic variety of $\me{M}$ as a $\me{D}_X$-module. The main result about the 
characteristic variety of $\me{E}_X$-modules is:
\begin{prop}
  Let $\me{M}$ be a coherent $\me{E}_X$-module. Then its charcteristic variety is 
  a closed analytic, involutive, $\C^\times$-conic subset of $T^*X$.
\end{prop}
\begin{defi}
  Let $\me{M}$ be a coherent $\me{E}_X$-module.\\
  One says that $\me{M}$ is holonomic if its charcteristic variety is Lagrangian.
\end{defi}
\begin{remark}
  \em Hence a holonomic $\me{D}_X$-module defines a holonomic $\me{E}_X$-module.
\end{remark}
Regular holonomic $\me{E}_X$-modules (or ``holonomic systems with regular 
singularities'') have been studied in \cite{KK}. A $\me{D}_X$-module is regular 
holonomic if and only if its associated $\me{E}_X$-module is regular holonomic. 
We do not recall the definition here since regular holonomic systems can be 
characterized by the Riemann-Hilbert correspondence and refer to loc. cit.
\begin{thm}
  Let $\me{M}$ be a regular holonomic $\me{E}_X$-module such that its characteristic 
  variety is in generic position at a point $p\in\dot T^*X$. Then there exists a
  regular holonomic $\me{D}_X$-modules $\widetilde{\me{M}}$ such that
  $$ \me{M} \simeq \me{E}_X\underset{\pi^{-1}\me{D}_X}{\otimes}
       \widetilde{\me{M}}. $$
\end{thm}
Regular holonomic $\me{E}_X$-modules form a stack of abelian categories. This 
stack is invariant by quantized contact transformations. Hence modulo a contact
transformation a regular holonomic $\me{E}_X$-module is locally isomorphic to
$$ \me{E}_X \underset{\pi^{-1}\me{D}_X}{\otimes} \on{RH}(\me{F}) $$
where $\me{F}$ is a perverse sheaf. For our purpose this can be taken as a definition
of a regular holonomic system. \\
Finally recall the following theorem.
\begin{prop}(\cite{An2}, Th\'eor\`eme 4.2.6 and Proposition 5.6.1)
  Let $\me{F}\in\Db_{\Cc}(k_X)$. Then
  $$ t\mu hom(\me{F},\me{O}_X)\simeq \me{E}_X^{\R,f}
     \otimes_{\pi^{-1}\me{D}_X}^L \pi^{-1}\thom(\me{F},\me{O}_X). $$  
  If $\me{F}$ is perverse then $t\mu hom(\me{F},\me{O}_X)$ is concentrated in
  degree 0 and 
   $$ \on{H}^0 t\mu hom(\me{F},\me{O}_X)\simeq \me{E}_X^{\R,f}
     \otimes_{\pi^{-1}\me{D}_X}\pi^{-1}\thom(\me{F},\me{O}_X). $$ 
\end{prop}

\subsection{The microdifferential structure of $\mu\me{O}_X$ and $\mu\me{O}_X^t$}

Recall that an object $\me{A}\in \I(k_X)$ is called a $k_X$-algebra if there exist
morphisms 
$$ k_X\ra \me{A} \qquad\qquad \me{A}\otimes \me{A} \lra \me{A} $$
that satisfy the usual conditions of unit and associativity (for more details, see 
for instance \cite{KS2}, Section 5.4). For example, If $\me{A}$ is a classical 
$k_X$-algebra in $\mc{M}od(k_X)$ then $\beta\me{A}$ is a $k_X$-algebra in $\I(k_X)$.\\
Let $\me{A}\in \I(k_X)$ be a $k_X$-algebra. A left $\me{A}$-module in $\I(k_X)$ is 
given by an object $\me{M}\in\I(k_X)$ and a structure morphism
$$   \me{A}\otimes \me{M} \lra \me{M} $$
that satisfies the usual compatibiulity conditions with the structure morphisms of 
$\me{A}$ (for more details, see loc. cit.).\\
Similarly, one defines the notion of a $k_X$-algebra $\me{A}$ in $\Db(\I(k_X))$ and 
the notion of a left  $\me{A}$-module in $\Db(\I(k_X))$. Note that if a $k_X$-algebra 
$\me{A}$ in $\Db(\I(k_X))$ is concentrated in a single degree, then it defines a 
$k_X$-algebra in $\I(k_X)$. However, even if $\me{A}$ is concentrated in a single 
degree, an $\me{A}$-module in $\Db(\I(k_X))$ is in general not well-defined in the 
derived category of $\me{A}$-modules in $\I(k_X)$. In order to avoid confuision, one 
often calls an $\me{A}$-module in $\Db(\I(k_X))$ a formal $\me{A}$-module.\\ 
In this section we will show that $\mu\me{O}_X$ is a $\beta\me{E}_X^{\R}$-module in
$\I(k_X)$. The same strategy will show that $\mu\me{O}_X^t$ is a formal 
$\beta\me{E}_X^{\R,f}$-module in the derived category $\Db(\I(k_X))$. In particular
its cohomology ind-sheaves are $\beta\me{E}_X^{\R,f}$-modules.\\
Let $f:\,X\ra Y$ be a morphism of complex manifolds.
\begin{lemma}
  There is a natural morphism in $\Db(\I(\C_X))$
  $$ \varphi: \Dr f_{!!}\Omega_X[d_X]\lra \Omega_Y[d_Y] $$
  such that $\alpha(\varphi)$ is the classical integration morphism.
\end{lemma}
\begin{proof}
Consider the chain of natural isomorphisms:
\begin{align*}
  \Hom{\Db(\I(\C_Y))}(&\Dr f_{!!}\Omega_X[d_X],\Omega_Y[d_Y])
   \simeq \Hom{\Db(\I(\C_X))}(\Omega_X[d_X],f^!\Omega_Y[d_Y])\\
   & \simeq \Hom{\Db(\C_X)}(\Omega_X[d_X],f^!\Omega_Y[d_Y])
   \simeq \Hom{\Db(\C_Y)}(\Dr f_!\Omega_X[d_X],\Omega_Y[d_Y]) 
\end{align*}
The inverse image of the classical integration morphism gives $\varphi$.
\end{proof}
There is a natural tempered version of the integration morphism:
\begin{prop}\label{tempint}
  There is a natural morphism in $\Db(\I(\C_X))$
  $$ \Dr f_{!!} \Omega_X^t[d_X] \lra \Omega_{Y}^t[d_Y]. $$
\end{prop}
\begin{proof}
This is a simple version (not respecting $\me{D}$-module structures) of 
the morphism established in \cite{KS2}: 
$$ \Dr f_{!!}(\Omega_X^t\overset{L}{\underset{\beta\me{D}_X}{\otimes}}
    \beta\me{D}_{X\ra Y})[d_X] \lra \Omega_Y^t[d_Y]. $$
\end{proof}
Now consider complex manifolds $X,Y,Z$ of complex dimensions $d_X,d_Y,d_Z$ and the 
diagram
$$ \xymatrix{ & X\times Y\times Z \ar[dl]_{q_{12}}
   \ar[d]^{q_{13}}\ar[dr]^{q_{23}} & \\
      X\times Y & X\times Z & Y\times Z. } $$
\begin{lemma}\label{combintmor}
  The integration morphisms induce natural morphisms
 \begin{align*} 
  \mu&\me{O}^{(d_X,0)}_{X\times Y}\overset{a}{\circ} 
  \mu\me{O}^{(d_Y,0)}_{Y\times Z} \lra
  \mu\me{O}^{(d_X,0)}_{X\times Z}[-d_Y]\\
  \mu&\me{O}^{t,(d_X,0)}_{X\times Y}\overset{a}{\circ} 
  \mu\me{O}^{t,(d_Y,0)}_{Y\times Z} \lra
  \mu\me{O}^{t,(d_X,0)}_{X\times Z}[-d_Y]
 \end{align*}
\end{lemma}
\begin{proof}
The two constructions being similar (and just an ind-variant of the 
construction used in Lemma 11.4.3. of \cite{KS3}) we will only show how to 
define the second morphism.\\
First let us consturct the natural morphism
$$ \me{O}^{t,(d_X,0)}_{X\times Y}\circ \me{O}^{t,(d_Y,0)}_{Y\times Z} \lra
   \me{O}^{t,(d_X,0)}_{X\times Z}[-\dim_{\C}Y]. $$
It can be obtained as follows
\begin{align*} 
  \me{O}^{t,(d_X,0)}_{X\times Y}\circ \me{O}^{t,(d_Y,0)}_{Y\times Z} 
       & \simeq \Dr p_{13!!}\big(p_{12}^{-1}\me{O}^{t,(d_X,0)}_{X\times Y}
         \otimes p_{23}^{-1}\me{O}^{t,(d_Y,0)}_{Y\times Z}\big) \\
   & \ra\Dr p_{13!!} \me{O}^{t,(d_X,d_Y,0)}_{X\times Y\times Z} \\
   & \ra \me{O}^{t,(d_X,0)}_{X\times Z}[-d_Y]
\end{align*}
where the last morphism is the integration morphism of Proposition 
\ref{tempint}. Using the microlocal composition formula we get a morphism
$$ \mu\me{O}^{t,(d_X,0)}_{X\times Y}
    \overset{a}{\circ} \mu\me{O}^{t,(d_Y,0)}_{Y\times Z} \lra
   \mu\big( \me{O}^{t,(d_X,0)}_{X\times Y}\circ \me{O}^{t,(d_Y,0)}_{Y\times Z} 
   \big) \lra \mu \me{O}^{t,(d_X,0)}_{X\times Z}[-d_Y] .$$
Therefore we get the morphism of the proposition.
\end{proof}
\begin{prop}\label{compositionisgreat}
  Let $\me{K}_1\in\Db(\I(\C_{T^*X\times T^*Y}))$ and $\me{K}_2\in
  \Db(\I(\C_{T^*Y\times T^*Z}))$. There are natural morphisms
  \begin{align*}
    \Dr\me{IH}om\big(\me{K}_1,&\mu\me{O}_{X\times Y}^{(d_X,0)}\big) 
    \overset{a}{\circ}
    \Dr\me{IH}om\big(\me{K}_2,\mu\me{O}_{Y\times Z}^{(d_Y,0)}\big) \lra
    \Dr\me{IH}om\big(\me{K}_1\overset{a}{\circ}\me{K}_2,
     \mu\me{O}_{X\times Z}^{(d_X,0)}\big)[-d_Y]\\ 
    \Dr\me{IH}om\big(\me{K}_1,&\mu\me{O}_{X\times Y}^{t,(d_X,0)}\big)
    \overset{a}{\circ}
    \Dr\me{IH}om\big(\me{K}_2,\mu\me{O}_{Y\times Z}^{t,(d_Y,0)}\big) \lra
    \Dr\me{IH}om\big(\me{K}_1\overset{a}{\circ}\me{K}_2,
    \mu\me{O}_{X\times Z}^{t,(d_X,0)}\big)[-d_Y]
  \end{align*}
  These morphisms satisfy the obvious associativity condition (analogous to
  Lemma 11.4.3 in \cite{KS3}).
\end{prop}
\begin{proof}
Note that there is a natural morphism
$$ \Dr\me{IH}om(\me{F}_1,\me{G}_1) \overset{a}{\circ}
   \Dr\me{IH}om(\me{F}_2,\me{G}_2) \lra
   \Dr\me{IH}om(\me{F}_1\overset{a}{\circ} \me{F}_2,
                \me{G}_1\overset{a}{\circ} \me{G}_2). $$
Combining this morphism with the morphisms of Lemma \ref{combintmor}, we get
the desired arrows.\\
The associativity condition is tedious to write down. It is a straightforward
consequence of the corresponding associativity conditions of the morphisms
involved in the construction.  
\end{proof}
\begin{remark}
  \em Consider $\me{K}_1\in\Db(\C_{X\times Y})$ and $\me{K}_2\in\Db(\C_{Y\times Z})$.
  Applying the functor $\alpha$ to the first morphism of Proposition 
  \ref{compositionisgreat}, we get a morphism
  $$ \Dr p_{13!}\Big(p_{12}^{a-1}\mu hom(\me{K}_1,\me{O}_{X\times Y}^{(d_X,0)}) 
     \otimes p_{23}^{-1}\mu hom(\me{K}_2,\me{O}_{Y\times Z}^{(d_Y,0)})\Big) \lra
     \mu hom(\me{K}_1\circ\me{K}_2,\me{O}_{X\times Z}^{(d_X,0)})[-d_Y]. $$
  This is precisely the morphism of \cite{KS3}, Lemma 11.4.3.
\end{remark}
In the situation of Proposition \ref{compositionisgreat} consider $Z=\{pt\}$
and $\me{K}_2=\mu\me{O}_Y$ (resp. $\me{K}_2=\mu\me{O}^t_Y$). 
Applying $\alpha$ we get natural morphisms
$$ \Dr p_{1!}\Big(\Dr\me{H}om\big(\me{K},\mu\me{O}_{X\times Y}^{(d_X,0)}
   \big)\otimes p_2^{a-1}\Dr\me{H}om\big(\mu\me{O}_Y,
   \mu\Omega_{Y}\big)\Big)
   \lra\Dr\me{H}om\big(\me{K}\overset{a}{\circ}\mu\me{O}_Y,
     \mu\Omega_X\big)[-d_Y], $$
$$ \Dr p_{1!}\Big(\Dr\me{H}om\big(\me{K},\mu\me{O}_{X\times Y}^{t,(d_X,0)}
   \big)\otimes p_2^{a-1}\Dr\me{H}om\big(\mu\me{O}_Y^t,
   \mu\Omega_{Y}^t\big)\Big)
   \lra\Dr\me{H}om\big(\me{K}\overset{a}{\circ}\mu\me{O}_Y^t,
     \mu\Omega_X^t\big)[-d_Y]. $$
Hence we get
\begin{equation}\label{morph0815} 
   \Dr p_{1!}\Big(\Dr\me{H}om\big(\me{K},\mu\me{O}_{X\times Y}^{(0,d_Y)}
   \big)\otimes p_2^{a-1}\Dr\me{H}om\big(\mu\me{O}_Y,
   \mu\me{O}_{Y}\big)\Big)
   \lra\Dr\me{H}om\big(\me{K}[d_Y]\overset{a}{\circ}\mu\me{O}_Y,
     \mu\me{O}_X\big), 
\end{equation}
\begin{equation}\label{morph0816} 
   \Dr p_{1!}\Big(\Dr\me{H}om\big(\me{K},\mu\me{O}_{X\times Y}^{t,(0,d_Y)}
   \big)\otimes p_2^{a-1}\Dr\me{H}om\big(\mu\me{O}_Y^t,
   \mu\me{O}_{Y}^t\big)\Big)
   \lra\Dr\me{H}om\big(\me{K}[d_Y]\overset{a}{\circ}\mu\me{O}_Y^t,
     \mu\me{O}_X^t\big). 
\end{equation}
Taking $X=Y$ and $\me{K}=\pi^{-1}\C_{\Delta_X}$ we get the morphisms
\begin{equation*}
   \me{E}^{\R}_X \otimes \Dr\me{H}om(\mu\me{O}_X,\mu\me{O}_X) \lra
      \Dr\me{H}om(\mu\me{O}_X,\mu\me{O}_X),
\end{equation*}
\begin{equation*} 
    \me{E}^{\R,f}_X \otimes \Dr\me{H}om(\mu\me{O}_X^t,\mu\me{O}_X^t) \lra
      \Dr\me{H}om(\mu\me{O}_X^t,\mu\me{O}_X^t).
\end{equation*}
Note that for any $\me{F}\in\Db(\I(\C_X))$ the identity of $\me{F}$ defines a natural
morphism $\C_{\Delta_X}\ra\Dr\me{H}om(\me{F},\me{F})$. Hence we get 
the structure morphisms
\begin{equation}\label{morph0819} 
   \beta\me{E}^{\R}_X \otimes \mu\me{O}_X \lra
      \mu\me{O}_X,
\end{equation}
\begin{equation}\label{morph0820} 
    \beta\me{E}^{\R,f}_X \otimes \mu\me{O}_X^t \lra
       \mu\me{O}_X^t.
\end{equation} 
In order to prove that the two morphisms \eqref{morph0819} and 
\eqref{morph0820} define structures of formal modules in $\Db(\I(\C_X))$ one uses
the associativity of the construction in Proposition \eqref{compositionisgreat}.
Therefore we get: 
\begin{prop}
  The object $\mu\me{O}_X$ (resp. $\mu\Omega_X$) is a left (resp. right) 
  $\beta\me{E}_X^\R$-module in $\Db(\I(\C_X))$.\\
  The object $\mu\me{O}_X^t$ (resp. $\mu\Omega_X^t$) is a left (resp. right) 
  $\beta\me{E}_X^{\R,f}$-module in $\Db(\I(\C_X))$.
\end{prop}
\begin{cor}\label{modulestructures}
  Let $\me{F}\in\Db(\I(\C_{T^*X}))$. Then
  \begin{itemize}
    \item[(1)] the complex $\Dr\me{IH}om(\me{F},\mu\me{O}_X)$ (resp.
       $\Dr\me{IH}om(\me{F},\mu\me{O}_X^t)$) is a left  
       $\beta\me{E}_X^\R$-module (resp. $\beta\me{E}_X^{\R,f}$-module) 
       in $\Db(\I(\C_X))$, 
    \item[(2)] the complex $\Dr\me{IH}om(\me{F},\mu\Omega_X)$  (resp.
       $\Dr\me{IH}om(\me{F},\mu\Omega_X^t)$) is a right 
       $\beta\me{E}_X^\R$-module (resp. $\beta\me{E}_X^{\R,f}$-module) 
       in $\Db(\I(\C_X))$.
  \end{itemize} 
\end{cor}
\begin{cor}\label{mo}
  \begin{itemize}
  \item[(i)]
  The object $\mu\me{O}_X$ (resp. $\mu\Omega_X$) is a left (resp. right)
  $\beta\me{E}_X^{\R}$-module in $\I(\C_{T^*X})$.
  \item[(ii)]
  For any $\me{F}\in\Db(\I(\C_{T^*X}))$ the complex $\Dr \me{IH}om(\me{F},\mu\me{O}_X)$
  (resp. $\Dr\me{IH}om(\me{F},\mu\Omega_X)$) is well defined in the bounded derived 
  category of left (resp. right) $\beta\me{E}_X^{\R}$-modules.
  In particular if $\me{F}\in\Db(k_X)$ then $\mu hom(\me{F},\me{O}_X)$ (resp.
  $\mu hom(\me{F},\Omega_X)$) is well defined in the bounded derived category of
  left (resp. right) $\me{E}_X^{\R}$-modules.
  \item[(iii)] Let $\me{F}\in\Db(\C_X)$. Then $\mu hom(\me{F},\me{O}_X)$ is well 
  defined in the derived category of $\me{E}_X^\R$-modules.
  \item[(iv)] Let $\me{F}\in\Db(\I(\C_{T^*X}))$ such that 
  $\Dr\me{H}om(\me{F},\mu\me{O}^t_X)$ is 
  concentrated in a single degree (for instance if $\me{F}$ is a microlocal perverse
  sheaf (see Lemma \ref{pervmodule} below)). Then the natural morphism
  $$ \Dr\me{H}om(\me{F},\mu\me{O}_X^t) \lra 
     \Dr\me{H}om(\me{F},\mu\me{O}_X) $$
  is well defined in the derived category of $\me{E}_X$-modules.
  \end{itemize}
\end{cor}
\begin{prop}
  The natural morphism
  \begin{align*} 
     \Dr p_{13!!}\Big(p_{12}^{a-1}\Dr\me{IH}om(\me{K}_1,
        \mu\me{O}_{X\times Y}^{(d_X,0)}) 
     \otimes & p_{23}^{-1}\Dr\me{IH}om(\me{K}_2,
       \mu\me{O}_{Y\times Z}^{(d_Y,0)})\Big) \\
     & \lra
     \Dr\me{IH}om(\me{K}_1\overset{a}{\circ} 
     \me{K}_2,\mu\me{O}_{X\times Z}^{(d_X,0)})[-d_Y] 
  \end{align*}
  factors through
   $$ \Dr p_{13!!}\Big(p_{12}^{a-1}\Dr\me{IH}om(\me{K}_1,
       \mu\me{O}_{X\times Y}^{(d_X,0)}) 
      \underset{p_2^{-1}\beta\me{E}_X^{\R}}{\overset{\on{L}}{\otimes}} 
      p_{23}^{-1}\Dr\me{IH}om(\me{K}_2,\mu\me{O}_{Y\times Z}^{(d_Y,0)})\Big). $$
\end{prop}
\begin{proof}
This is a consequence of the associativity condition of Proposition 
\ref{compositionisgreat}.
\end{proof}

\subsection{Quantized contact transformations for $\mu\me{O}_X$ and $\mu\me{O}_X^t$}

In this section we will adapt Section 11.4 of \cite{KS3} to study the 
behaviour of $\mu\me{O}_X$ and $\mu\me{O}^t_X$ under complex contact 
transformations. We will restrict ourselves now to $\mu\me{O}_X$, the study of 
$\mu\me{O}^t_X$ being similar. We also did not include all proofs which are rather 
similar to Section 11.4 of loc. cit.\\
Our main interest is to prove that the complex
$$ \Dr\me{IH}om_{\beta\me{E}_X}(\beta\me{M},\mu\me{O}_X)$$  
is invariant under quantized contact transformations where $\me{M}$ is a coherent
$\me{E}_X$-module.\\
More precisely consider a contact transformation 
$$ \chi:\ \Omega_X \overset{\sim}{\lra} \Omega_Y $$
Then locally we may find a kernel $\me{K}$ such that $\Phi_{\mu\me{K}}^a$ induces
an equivalence of microlocal perverse sheaves and 
Since $\me{M}$ is coherent, the natural morphism
$$ \mu\me{K}\overset{a}{\circ}\Dr\me{IH}om_{\beta\me{E}_X}(\beta\me{M},\mu\me{O}_X) 
    \lra \Dr\me{IH}om_{\beta(\chi_*\me{E}_X)}(\beta(\chi_*\me{M}),\mu\me{K}
   \overset{a}{\circ}\mu\me{O}_X) $$
is an isomorphism. Hence we are reduced to construct an isomorphism
$$ \mu\me{K}\overset{a}{\circ}\mu\me{O}_Y \lra \mu\me{O}_X $$
that is compatible with the action of $\beta\me{E}_X$ (resp. $\beta\me{E}_Y$) on 
$\mu\me{O}_X$ (resp. $\mu\me{O}_Y$). Such a morphism has been constructed in 
$\Db(k_X,p)$ in Section 11.4 of \cite{KS3}.\\
 Recall the morphism \eqref{morph0815} 
\begin{equation*}
   \Dr p_{1!}\Big(\Dr\me{H}om\big(\me{K},\mu\me{O}_{X\times Y}^{(0,d_Y)}
   \big)\otimes p_2^{a-1}\Dr\me{H}om\big(\mu\me{O}_Y,
   \mu\me{O}_{Y}\big)\Big)
   \lra\Dr\me{H}om\big(\me{K}[d_Y]\overset{a}{\circ}\mu\me{O}_Y,
     \mu\me{O}_X\big). 
\end{equation*}
Note that if $\me{K}\simeq \mu \on{K}$ and $\supp \on{K}\ra X$ is proper, then we get 
the commutative diagram
$$ \xymatrix{
    \Dr p_{1!}\Big(\mu hom\big(\on{K},\me{O}_{X\times Y}^{(0,d_Y)}
   \big)\otimes p_2^{a-1}\mu hom\big(\me{O}_Y,
   \me{O}_{Y}\big)\Big) \ar[r] \ar[dr] 
    & \Dr\me{H}om\big(\mu\!\on{K}[d_Y]\overset{a}{\circ}\mu\me{O}_Y,
     \mu\me{O}_X\big) \\
   & \mu hom\big(\on{K}[d_Y]\circ\me{O}_Y,
     \me{O}_X\big) \ar[u] } $$
where the diagonal morphism is the classical morphism from \cite{KS3}. 

Now suppose that $\Dr\me{H}om\big(\me{K},\mu\me{O}_{X\times Y}^{(0,d_Y)}\big)$ is 
concentrated in positive degrees and that 
$$\supp \Dr\me{H}om\big(\me{K},
 \mu\me{O}_{X\times Y}^{(0,d_Y)}\big)\ra X $$ 
is proper. Then by taking the $0$-cohomology we get a morphsm
\begin{equation}\label{morph0821} 
   p_{1*}\Big(\on{H^0}\!\Dr\me{H}om\big(\me{K},\mu\me{O}_{X\times Y}^{(0,d_Y)}
   \big)\otimes p_2^{a-1}\on{H^0}\!\Dr\me{H}om\big(\mu\me{O}_Y,
   \mu\me{O}_{Y}\big)\Big)
   \lra \on{H}^0\!\Dr\me{H}om\big(\me{K}[d_Y]\overset{a}{\circ}\mu\me{O}_Y,
     \mu\me{O}_X\big) 
\end{equation}
Hence the identity of $\mu\me{O}_Y$ and any section 
$$ s\in\on{H}^0\!\Dr\me{H}om(\me{K},\mu\me{O}_{X\times Y}^{(0,d_Y)}). $$
defines a morphism
$$ \varphi(s):\,\me{K}[d_Y]\overset{a}{\circ}\mu\me{O}_Y \lra
   \mu\me{O}_X $$
In the sequel we will only consider kernels $\me{K}\in\Db(\C_{T^*X\times T^*Y})$
satisfying
\begin{itemize}
  \item[(i)] $\me{K}$ is $\R$-constructible,
  \item[(ii)] $(\Omega_Y\times T^*X\cup T^*Y\times \Omega_X^a)\cap \SS(\me{K})\subset
               \Lambda$,
  \item[(iii)] $\me{K}$ is simple\begin{footnote}{For the definition of simple sheaves
    see \cite{KS3}, Section 7.5.}\end{footnote} with shift $0$ along $\Lambda$.
\end{itemize}
Note that given a $\C^\times$-conic Lagrangian subvariety $\Lambda$ that is 
associated to a contact transformation, there always 
locally (on $P^*(X\times Y)$) exists $\me{K}$, i.e. for each 
$p\in\Omega_X$ there exists a $\C^\times$-conic open neighborhood $\Omega_X'$ of 
$\C^\times p$ such that (i), (ii) and (iii) are satisfied.  
Recall the morphism
\begin{align*} 
  p_{13*}^a\Big(p_{12}^{a-1}\on{H}^0\Dr\me{H}om(\me{K}_1,\mu
   \me{O}_{X\times Y}^{(0,d_Y)})&\underset{\beta\me{E}_Y^{\R}}{\otimes}
   p_{23}^{a-1}\Dr\me{H}om(\me{K}_2,\mu\me{O}_{Y\times Z}^{(0,d_Z)})\Big) \\
  & \lra \Dr\me{IH}om(\me{K}_1\overset{a}{\circ}\me{K}_2[d_Y],\mu
   \me{O}_{X\times Z}^{(0,d_Z)}. 
\end{align*}
Denote by $s\circ s'$ the image of $s\otimes s'$ by this morphism.\\
Then we get 
\begin{prop}
  $$ \varphi(s\circ s')=\varphi(s)\circ(\me{K}[n]\overset{a}{\circ}\alpha(s')) $$
\end{prop}
\begin{proof}
This follows from the fact that the morphism of Proposition \ref{compositionisgreat} 
satisfies to the obvious associativity condition.
\end{proof}
\begin{thm}\label{quantizationthm}
  For every $p\in \Omega_X$ there exists a $\C^\times$-conic open neighborhood 
  $\Omega_X'\subset\Omega_X$ of $\C^\times p$ auch that if we set $\Omega_Y'=
  \chi(\Omega_Y)$ we can find a section 
 $$ s\in \on{H}^0\Dr\me{H}om(\mu\me{K},\mu\me{O}_{X\times Y}^{(d_X,0)})|_{\Omega_{X}'
   \times \Omega_{Y}'} $$
  such that the morphisms
  \begin{eqnarray*}   
\me{E}_{X}^\R|_{\Omega_X} & \lra & 
   p_{1*}\on{H}^0\Dr\me{H}om(\mu\me{K},\mu\me{O}_{X\times Y}^{(d_X,0)})|_{\Omega_{X}'
   \times \Omega_{Y}'} \qquad ;\quad P\mapsto Ps \\
    \me{E}_Y^\R|_{\Omega_Y} & \lra &
   p_{2*}\on{H}^0\Dr\me{H}om(\mu\me{K},\mu\me{O}_{X\times Y}^{(d_X,0)})|_{\Omega_{X}'
   \times\Omega_{Y}'} \qquad ;\quad Q\mapsto sQ 
  \end{eqnarray*}
  are isomorphisms and we get an antiisomorphism
  $$ \chi_*\me{E}_X^\R|_{\Omega_X} \overset{\sim}{\lra}
       \me{E}_Y^\R|_{\Omega_Y} \quad ;\quad P\mapsto Q\ such\ that\ Ps=sQ $$ 
  For such a section $s$ the morphism
  $$  \varphi(s):\ \me{K}[d_Y]\overset{a}{\circ}\mu\me{O}_Y \lra
   \mu\me{O}_X $$
  is an isomorphism of $\beta\me{E}_X^{\R}$-modules.
\end{thm}
\begin{cor}
  Let $s$ be a section as in Theorem \ref{quantizationthm} and $\me{F}$ be a
  perverse sheaf. Then $\varphi(s)$ defines an isomorphism of $\me{E}_X^{\R}$-modules
 $$ \chi_*\Dr\me{H}om(\pi^{-1}\mu\me{F},\mu\me{O}_Y) \simeq
    \Dr\me{H}om(\mu\me{K}\overset{a}{\circ}\mu\me{F},
      \mu\me{O}_X). $$
\end{cor}
\begin{cor}
  Let $s$ be a section as in Theorem \ref{quantizationthm} and $\me{F}$ be a 
  coherent $\me{E}_Y$-module. Then $\varphi(s)$ defines an isomorphism
 $$ \Dr\me{IH}om_{\beta\me{E}_Y}(\beta\me{F},\mu\me{O}_Y) \simeq
    \Dr\me{IH}om_{\beta\me{E}_X}(\chi_*\beta\me{F},
      \mu\me{O}_X) $$
\end{cor}
The tempered version of Theorem \ref{quantizationthm} is similar with the only 
difference that the microdifferential structures are formal structures in the 
derived categories.
\begin{thm}\label{temperedquantizationthm}
  For every $p\in \Omega_X$ there exists a $\C^\times$-conic open neighborhood 
  $\Omega_X'\subset\Omega_X$ of $\C^\times p$ auch that if we set $\Omega_Y'=
  \chi(\Omega_Y)$ we can find a section 
 $$ s\in \on{H}^0\Dr\me{H}om(\mu\me{K},\mu\me{O}_{X\times Y}^{t,(d_X,0)})|_{\Omega_{X}'
   \times \Omega_{Y}'} $$
  such that the morphisms
  \begin{eqnarray*}   
\me{E}^{\R,f}_{X}|_{\Omega_X} & \lra & 
   p_{1*}\on{H}^0\Dr\me{H}om(\mu\me{K},\mu\me{O}_{X\times Y}^{t,(d_X,0)})|_{\Omega_{X}'
   \times \Omega_{Y}'} \qquad ;\quad P\mapsto Ps \\
    \me{E}_Y^{\R,f}|_{\Omega_Y} & \lra &
   p_{2*}\on{H}^0\Dr\me{H}om(\mu\me{K},\mu\me{O}_{X\times Y}^{t,(d_X,0)})|_{\Omega_{X}'
   \times\Omega_{Y}'} \qquad ;\quad Q\mapsto sQ 
  \end{eqnarray*}
  are isomorphisms and we get an antiisomorphism
  $$ \chi_*\me{E}_X^{\R,f}|_{\Omega_X} \overset{\sim}{\lra}
       \me{E}_Y^{\R,f}|_{\Omega_Y} \quad ;\quad P\mapsto Q\ such\ that\ Ps=sQ $$ 
  For such a section $s$ the morphism
  $$  \varphi^t(s):\ \me{K}[d_Y]\overset{a}{\circ}\mu\me{O}_Y^t \lra
   \mu\me{O}_X^t $$
  is an isomorphism of formal $\beta\me{E}_X^{\R,f}$-modules in $\Db(\I(k_X))$.
\end{thm}
\begin{cor}
  Let $s$ be a section as in Theorem \ref{temperedquantizationthm} and $\me{F}$ be a
  perverse sheaf. Then $\varphi(s)$ defines an isomorphism of $\me{E}_X^{\R,f}$-modules
 $$ \chi_*\Dr\me{H}om(\pi^{-1}\mu\me{F},\mu\me{O}_Y^t) \simeq
    \Dr\me{H}om(\mu\me{K}\overset{a}{\circ}\mu\me{F},
      \mu\me{O}_X^t). $$
\end{cor}

\subsection{Classical Riemann-Hilbert Theorem and Ind-sheaves}

Using the functor $\mu$ we can reformulate the isomorphisms of the classical
Riemann-Hilbert Theorem from the microlocal point of view. 
\begin{lemma}
\begin{itemize}
   \item[(1)] Let $\me{F}\in\Db_{\Cc}(k_X)$. Then
      \begin{eqnarray*}  
        \Dr\me{H}om(\mu \me{F},\mu \me{O}^t_X) & \simeq & \me{E}_{X}^{\R,f}
          \overset{\on{L}}{\underset{\pi^{-1}\me{D}_X}{\otimes}} 
          \pi^{-1}\on{RH}(\me{F}), \\
        \gamma^{-1}\Dr \gamma_*\Dr\me{H}om(\mu \me{F},\mu \me{O}^t) & \simeq & 
          \me{E}_{X}\overset{\on{L}}{\underset{\pi^{-1}\me{D}_X}{\otimes}} 
          \pi^{-1}\on{RH}(\me{F}).
      \end{eqnarray*}
   \item[(2)] Let $\me{M}$ be a coherent $\me{D}_X$-module. Then
       \begin{eqnarray*}
          \Dr \me{IH}om_{\beta\me{E}_X^{\R,f}}(\beta\me{E}_X^{\R,f}
           \overset{\on{L}}{\underset{\pi^{-1}\beta\me{D}_X}{\otimes}} 
         \beta(\pi^{-1}\me{M}),\mu \me{O}_X) & \simeq & \mu(\on{Sol}(\me{M})), \\
            \Dr \me{IH}om_{\beta\me{E}_X}(\beta\me{E}_X
           \overset{\on{L}}{\underset{\pi^{-1}\beta\me{D}_X}{\otimes}} 
         \beta(\pi^{-1}\me{M}),\mu \me{O}_X) & \simeq & \mu(\on{Sol}(\me{M})).
        \end{eqnarray*}
\end{itemize}
\end{lemma}
\begin{proof}
(1) is a reformulation of Theorem 4.2.6 of \cite{An2} in terms of ind-sheaves using 
the fact that $\gamma^{-1}\Dr\gamma_*\me{E}^{\R,f}_X\simeq\me{E}_X$ and (2) follows 
directly from the fact that $\me{M}$ is a coherent $\me{D}_X$-module.
\end{proof}
Therefore we can now formulate Riemann-Hilbert Theorem in terms of ind-sheaves:
\begin{prop}\label{RieHil}
\begin{itemize}
   \item[(1)] Let $\me{F}$ be a perverse sheaf on $X$. Then
       \begin{eqnarray*}
     \Dr\me{IH}om_{\beta\me{E}_X^{\R,f}}(\beta\Dr\me{H}om(\mu \me{F},
        \mu \me{O}^t_X),\mu\me{O}_X) & \simeq &  \mu\me{F} \\
     \Dr \me{IH}om_{\beta\me{E}_X}(\beta\gamma^{-1}\Dr\gamma_*\Dr\me{H}om(\mu\me{F},
        \mu \me{O}^t_X),\mu\me{O}_X) & \simeq & \mu\me{F} 
        \end{eqnarray*}
   \item[(2)] Let $\me{M}$ be a regular holonomic $\me{D}_X$-module. Then
       \begin{eqnarray*} 
         \Dr\me{H}om(\Dr \me{IH}om_{\beta\me{E}_X^{\R,f}}(
         \beta(\me{E}_X^{\R,f}\overset{\on{L}}{\underset{\pi^{-1}\me{D}_X}{\otimes}} 
         \pi^{-1}\me{M}),\mu \me{O}_X),\mu\me{O}^t_X) 
         & \simeq  &\me{E}_X^{\R,f}\overset{\on{L}}{\underset{\pi^{-1}\me{D}_X}
          {\otimes}}\pi^{-1}\me{M} \\
          \gamma^{-1}\Dr\gamma_*\Dr\me{H}om(\Dr \me{IH}om_{\beta\me{E}_X}
         (\beta(\me{E}_X\overset{\on{L}}{\underset{\pi^{-1}\me{D}_X}{\otimes}} 
          \pi^{-1}\me{M},\mu \me{O}_X),\mu\me{O}^t_X) 
         & \simeq & \me{E}_X\overset{\on{L}}{\underset{\pi^{-1}\me{D}_X}
          {\otimes}}\pi^{-1}\me{M}
       \end{eqnarray*}
\end{itemize}
\end{prop}
Now let us formulate the comparison theorem for regular holonomic 
$\me{E}_X$-modules in terms of ind-sheaves. The classical version (\cite{An1}, 
Proposition 5.6.3) states
\begin{prop}
  Let $\me{M}$ be a regular holonomic $\me{E}_X$-module and $\me{G}\in\Db_{\Rc}(k_X)$
  such that $t\mu hom(\me{G},\me{O}_X)$ is concentrated in a single degree. 
  Then the natural morphism
  $$ \Dr \me{H}om_{\me{E}_X}(\me{M},t\mu hom(\me{G},\me{O}_X)) \lra
     \Dr \me{H}om_{\me{E}_X}(\me{M},\mu hom (\me{G},\me{O}_X)) $$
  is an isomorphism.
\end{prop}
\begin{remark}
  \em Note that we do not have to assume that $\mu hom (\me{G},\me{O}_X)$ is 
  concentrated in degree zero because we now know that $\mu hom(\me{G},\me{O}_X)$
  is well defined in the derived category of microdifferential modules.\\
  It is slightly more complicated to deal with $\mu\me{O}_X^t$ since we do not know
  if $\mu\me{O}_X^t$ is well defined in the derived category of microdifferential
  modules. The analog of the last proposition should be the formula
   $$ \Dr\me{IH}om_{\beta\me{E}_X}(\beta\me{M},\mu \me{O}_X)
     \overset{\sim}{\lra} 
     \Dr\me{IH}om_{\beta\me{E}_X}(\beta\me{M},\mu\me{O}^t_X)  $$
  for any regular holonomic $\me{E}_X$-module $\me{M}$. However the second term of this
  isomorphism is unfortunately not (yet) well-defined. Therefore we only get the 
  following weaker statement:
\end{remark}
\begin{prop}
  Let $\me{M}$ be a regular holonomic $\me{E}_X$-module. Then there is a natural 
  morphism
  $$ \Dr\me{IH}om_{\beta\me{E}_X}(\beta\me{M},\mu \me{O}_X) 
      \lra \Dr\me{IH}om(\beta\me{M},\mu\me{O}^t_X). $$
\end{prop}
\begin{proof}
By the comparison theorem if $\me{G}\in\Db_{\on{perv}}(\C_X,U)$ we have a 
natural morphism
\begin{equation*}
     \Dr\me{H}om_{\me{E}_X}(\me{M},\mu hom(\me{G},\me{O}_X))    \simeq 
    \Dr\me{H}om_{\me{E}_X}(\me{M},t\mu hom(\me{G},\me{O}_X)) 
    \lra\Dr\me{H}om(\me{M},t\mu hom(\me{G},\me{O}_X)).
\end{equation*}
Since we have
\begin{align*}
  \Dr\me{H}om_{\me{E}_X}(\me{M},\mu hom(\me{G},\me{O}_X))& \simeq
   \Dr\me{H}om_{\me{E}_X}(\me{M},\Dr \me{H}om(\mu\me{G},\mu \me{O}_X)) \\
 &\simeq \Dr\me{H}om_{\beta\me{E}_X}(\beta \me{M},\Dr 
       \me{IH}om(\mu\me{G},\mu \me{O}_X))\\
 &\simeq \Dr\me{H}om(\mu\me{G},\Dr 
       \me{IH}om_{\beta\me{E}_X}(\beta \me{M},\mu \me{O}_X))
\end{align*}
and
\begin{align*}
  \Dr\me{H}om(\me{M},t\mu hom(\me{G},\me{O}_X))& \simeq
   \Dr\me{H}om(\me{M},\Dr \me{H}om(\mu\me{G},\mu \me{O}_X^t) \\
 &\simeq \Dr\me{H}om(\beta \me{M},\Dr 
       \me{IH}om(\mu\me{G},\mu \me{O}_X^t))\\
 &\simeq \Dr\me{H}om(\mu\me{G},\Dr 
       \me{IH}om(\beta \me{M},\mu \me{O}_X^t))
\end{align*}
we get a morphism
$$ \Dr\me{H}om(\mu\me{G},\Dr 
       \me{IH}om_{\beta\me{E}_X}(\beta \me{M},\mu \me{O}_X)) \lra  
    \Dr\me{H}om(\mu\me{G},\Dr 
       \me{IH}om(\beta \me{M},\mu \me{O}_X^t)). $$
Now recall that $\Dr\me{IH}om_{\beta\me{E}_X}(\me{M},\mu\me{O}_X)$ is a microlocal 
perverse sheaf. Hence locally it is of the form $\mu\me{G}$ for some object
$\me{G}\in\Db_{\on{perv}}(\C_X,U)$. Thus locally the identity morphism of 
$\Dr\me{IH}om_{\beta\me{E}_X}(\me{M},\mu\me{O}_X)$ defines the desired morphism and we 
may patch it because microlocal perverse sheaves form a stack.
\end{proof}
\begin{remark}
 \em If $\mu\me{O}^t_X$ was well-defined in the derived category of 
 $\beta\me{E}_X$-modules then the proof of the last Proposition would establish 
 the isomorphism
  $$ \Dr\me{IH}om_{\beta\me{E}_X}(\beta\me{M},\mu \me{O}_X) \overset{\sim}{\lra}
      \Dr\me{IH}om_{\beta\me{E}_X}(\beta\me{M},\mu\me{O}^t_X).  $$
\end{remark}

\subsection{Microlocal Riemann-Hilbert morphism}

\begin{lemma}\label{pervmodule}
  Let $\me{F}\in\mu Perv(\Omega)$ and set $U=\gamma^{-1}(\Omega)$. Then
  $$\Dr \me{H}om(\me{F},\mu\me{O}_X^t|_{U})$$ 
  is an $\me{E}_X^{\R,f}|_{U}$-module. Moreover
  $$ \gamma^{-1}\Dr\gamma_*\Dr \me{H}om(\me{F},\mu\me{O}_X^t|_{U}) $$
  is a regular holonomic $\me{E}_X|_{U}$-module.
\end{lemma}
\begin{proof}
First we will show that $\Dr \me{H}om(\me{F},\mu\me{O}_X^t|_{U})$ is a well-defined 
$\me{E}_X^{\R,f}|_{U}$-module.\\
By Proposition \ref{mo} it is enough to prove that the complex
 $$\Dr \me{H}om(\me{F},\mu\me{O}_X^t|_{U})$$ 
is concentrated in a single degree. This is a local problem. Thus we may assume 
that $\me{F}\simeq \mu\tilde{\me{F}}$ where $\widetilde{\me{F}}$ is an object of 
$\Db_{\on{perv}}(k_X,U)$. Therefore
$$ \Dr \me{H}om(\me{F},\mu\me{O}^t_{X}|_{U}) \simeq
    \Dr \me{H}om(\mu\widetilde{\me{F}}|_U,\mu\me{O}^t_{X}|_{U}) \simeq 
   t\mu hom(\widetilde{\me{F}},\me{O}_X)|_U $$
Since $t\mu hom(\tilde{\me{F}},\me{O}_X)|_U$ is invariant under quantized
contact transformations, we can suppose that $\widetilde{\me{F}}$ is a perverse sheaf. 
Then the complex $t\mu hom(\tilde{\me{F}},\me{O}_X)|_U$ is concentrated in 
degree $0$. Hence $\Dr \me{H}om(\me{F},\mu\me{O}_X^t|_{U})$ is a well-defined 
$\me{E}_X^{\R,f}|_{U}$-module.\\
Hence $\gamma^{-1}\Dr\gamma_*\Dr \me{H}om(\me{F},\mu\me{O}_X^t|_{U})$ is an 
$\me{E}_X$-module. Let us show that it is regular holonomic. This is again a local
question, invariant by quantized contact transformations. Therefore we may
assume that $\me{F}\simeq\mu\widetilde{\me{F}}$ for a perverse sheaf 
$\widetilde{\me{F}}$. Recall that
  $$ t\mu hom(\tilde{\me{F}},\me{O}_{X})|_U \simeq
      \big(\me{E}^{\R,f}_X \underset{\pi^{-1} \me{D}_X}{\otimes}
      \pi^{-1}\thom(\tilde{\me{F}},\me{O}_X)\big)|_U.$$
By the Riemann-Hilbert theorem $\thom(\tilde{\me{F}},\me{O}_X)$ is a regular holonomic 
$\me{D}_X$-module. Hence $\gamma^{-1}\Dr\gamma_*\Dr \me{H}om(\me{F},
\mu\me{O}_X^t|_{U})$ is regular holonomic.
\end{proof}
\begin{lemma}
  Let $\me{M}\in\me{H}ol\me{R}eg(U)$. Then
  $$ \Dr \me{IH}om_{\beta\me{E}_X|_U}(\beta\me{M},\mu\me{O}_X|_{U}) $$
  is an object of $\mu Perv(\Omega)$.
\end{lemma}
\begin{proof}
This is a local problem, invariant by quantized contact transformations.
Therefore we may assume that $\me{M}$ is isomorphic to 
$\me{E}_X\underset{\pi^{-1}\me{D}_X}\otimes\pi^{-1}\tilde{\me{M}}$ on $U$ where 
$\tilde{\me{M}}$ is a regular holonomic $\me{D}_X$-module. Then
$$ \Dr \me{IH}om_{\beta\me{E}_X|_U}(\beta\me{M},\mu\me{O}_X|_U) \simeq
    \Dr \me{IH}om_{\beta\me{E}_X}(\me{E}_X\underset{\pi^{-1}\me{D}_X}\otimes  
       \pi^{-1}\tilde{\me{M}},\mu\me{O}_X)|_U \simeq
    \mu \Dr\me{H}om_{\me{D}_X}(\tilde{\me{M}},\me{O}_X)|_U $$
By the Riemann-Hilbert Theorem $\Dr\me{H}om_{\me{D}_X}(\tilde{\me{M}},\me{O}_X)$ 
is a perverse sheaf.  
\end{proof}
For any open subset $\Omega\subset P^*X$ let us define the microlocal 
Riemann-Hilbert correspondance:
$$ \xymatrix@C=3cm{ {\mu Perv(\Omega)} \ar@<+3pt>[r]^{\murh} & 
       {\mc{H}ol\mc{R}eg(\me{E}_X|_{\gamma^{-1}\Omega}).} 
        \ar@<+3pt>[l]^{\musol} } $$
by the formulas
$$ \musol(\me{M})= 
  \Dr\me{IH}om_{\beta(\me{E}_X|_{\gamma^{-1}\Omega})}(\beta(\me{M}),
    \mu \me{O}_X|_{\gamma^{-1}\Omega})$$
$$ \murh(\me{F})=\gamma^{-1}_{\Omega}\Dr \gamma_{\Omega *}
  (\Dr\me{H}om(\me{F},\mu \me{O}^t|_{\gamma^{-1}(\Omega)})) $$
where $\gamma_{\Omega}$ is the restriction of $\gamma$ to $\gamma^{-1}\Omega$.\\
The functors $\musol$ and $\murh$ are obviously functors of stacks.

\begin{lemma}\label{morphism2}
   There is a natural morphism
  \begin{equation}\label{morphism2abc}
     \on{Id} \lra \musol\circ\murh. 
  \end{equation}
\end{lemma}
\begin{proof}
Let $\me{F}$ be a microlocal perverse sheaf. We will define the morphism 
\eqref{morphism2abc} by a natural element of
$$ \on{Hom}_{\Db(\I(k_X))}(\me{F},\Dr\me{IH}om_{\beta\me{E}_X}(\beta
   \gamma^{-1}\Dr\gamma_*\Dr\me{H}om(\me{F},\mu\me{O}^t_X),\mu\me{O}_X). $$
Note that 
\begin{align*}
   \on{Hom}_{\Db(\I(k_X))}(\me{F},\Dr\me{IH}om_{\beta\me{E}_X}&(\beta
   \gamma^{-1}\Dr\gamma_*\Dr\me{H}om(\me{F},\mu\me{O}^t_X),\mu\me{O}_X) \\
   & \simeq
   \on{Hom}_{\Db(\I(\beta\me{E}_X))}(\me{F}\otimes\beta\gamma^{-1}\Dr\gamma_*
   \Dr\me{H}om(\me{F},\mu\me{O}^t_X),\mu\me{O}_X). 
\end{align*}
Now recall that
$$ \beta\gamma^{-1}\Dr\gamma_*\Dr\me{H}om(\me{F},\mu\me{O}^t_X) \lra
   \beta\gamma^{-1}\Dr\gamma_*\Dr\me{H}om(\me{F},\mu\me{O}_X) $$
is $\beta\me{E}_X$-linear. Hence the natural morphism in $\Db(\I(\beta\me{E}_X))$
$$ \me{F}\otimes\beta\gamma^{-1}\Dr\gamma_*\Dr\me{H}om(\me{F},\mu\me{O}_X) 
   \lra \mu\me{O}_X $$
defines the morphism of the lemma.
\end{proof}
\begin{thm}
  The functors $\musol$ and  $\murh$ define quasi-inverse equivalences of 
  stacks
  $$ \xymatrix@C=3cm{ {\mu Perv} \ar@<+3pt>[r]^{\murh} & 
       {\gamma_*\mc{H}ol\mc{R}eg(\me{E}_X).} 
        \ar@<+3pt>[l]^{\musol} } $$
\end{thm}
\begin{proof}
First let us show that the morphism of Lemma \ref{morphism2} is an isomorphisms.\\
Let $\me{F}$ be a microlocal perverse sheaf defined in a neighborhood at $p$.
Then there exists $\widetilde{\me{F}}\in\Db_{\on{perv}}(k_X,\C^\times p)$ such that
$\me{F}\simeq\mu\widetilde{\me{F}}$. Let $\chi$ be a contact transformation such 
that $\chi(\SS(\widetilde{\me{F}}))$ is in generic position at $\chi(p)$. Then
$$
 \Phi_{\me{K}}\Dr\me{IH}om_{\beta\me{E}_X}(\beta\Dr\me{H}om(\mu\widetilde{\me{F}},
   \mu\me{O}^t_X,),\mu\me{O}_X))\simeq 
 \Dr\me{IH}om_{\beta\me{E}_X}(\beta\Dr\me{H}om(\Phi_{\me{K}}\widetilde{\me{F}},
   \mu\me{O}^t_X,),\mu\me{O}_X))
$$
Since $\Phi_{\me{K}}(\widetilde{\me{F}})$ is isomorphic to a perverse sheaf in a 
neighborhood of $\pi(p)$ the isomorphism follows from the second part
of Proposition \ref{RieHil}.\\
Now let us show that the functor $\mu\on{RH}$ is an equivalence of stacks. It is 
sufficent to prove this locally. In order to prove that $\mu\on{RH}$ is essentially 
surjective, we show that $\mu \on{RH}$ and $\mu\on{Sol}$ are inverse to each other on 
the level of objects. We have already seen that $\mu\on{Sol}\circ\mu\on{RH}(\me{F})
\simeq \me{F}$ for any microlocal perverse sheaf $\me{F}$.\\
Let $\me{M}$ be a regular holonomic $\me{E}_X$-module defined in a neighborhood of
$p\in\dot T^*X$. Let $\chi$ be a contact transformation such that $\chi_*\me{M}$
is in generic position at $\chi(p)$. Then
\begin{align*}
   \chi_*\Dr \gamma^{-1}\Dr\gamma_*&\Dr\me{H}om(\Dr\me{IH}om_{\beta\me{E}_X}( 
   \beta(\me{M}),\mu\me{O}_X),\mu\me{O}_X^t) \\
   & \simeq
   \gamma^{-1}\Dr\gamma_*\Dr\me{H}om(\Phi_{\me{K}}^\mu\Dr\me{IH}om_{\beta\me{E}_X}( 
   \beta(\me{M}),\mu\me{O}_X),\mu\me{O}_X^t) \\
   & \simeq  
   \gamma^{-1}\Dr\gamma_*\Dr\me{H}om(\Dr\me{IH}om_{\beta\me{E}_X}( 
   \beta\chi_*\me{M},\mu\me{O}_X),\mu\me{O}_X^t)
\end{align*}
But since $\chi_*\me{M}$ is in generic position there exists a regular holonomic
$\me{D}_X$-module $\widetilde{\me{M}}$ such that 
$$\me{M}\simeq \me{E}_X\overset{\on{L}}
{\underset{\pi^{-1}\me{D}_X}{\otimes}}\pi^{-1}\widetilde{\me{M}}. $$
Hence $\mu\on{RH}$ is essentially surjective.\\
Let us show that $\mu \on{RH}$ is fully faithful. Let $\me{F},\me{G}$ be microlocal
perverse sheaves. By invariance under quantized contact transformations we may assume
that there exists perverse sheaves $\widetilde{\me{F}},\widetilde{\me{G}}$ such that
$\me{F}\simeq \mu\widetilde{\me{F}}$ and $\me{G}\simeq\widetilde{\me{G}}$. Then the
fact that $\mu\on{RH}$ is fully faithful follows from the well-known formula
$$ \mu hom(\me{F},\me{G})\simeq \Dr\me{H}om(\on{RH}(\me{G}),\on{RH}(\me{F})). $$
\end{proof}

\appendix

\section{Stacks on Topological Spaces}

In this section we will make use of the language of 2-functors and 2-colimits. But since
we only work with the 2-category $\mc{CAT}$ of all (small) categories, we will not
recall the abstract theory and refer to \cite{St},\cite{McL}.

\subsection{Prestacks}

A prestack is a ``presheaf of categories up to equivalence''. More precisely,
let $X$ be a topological space and denote by $\mc{T}(X)$ the category of open 
sets of $X$. A prestack on $X$ is just a 2-functor $\mc{T}(X)^\circ\ra\mc{CAT}$. 
Therefore many authors call a $2$-functor with domain $\mc{I}$ a prestack on $\mc{I}$.\\
Since the category $\mc{T}(X)$ is particularly simple, we will recall here a detailed
description.
\begin{defi}
  A \textbf{prestack} $(\me{C},\rho,\Phi)$ on $X$
  consists of 
  \begin{itemize}
    \item[(1)] a small category $\me{C}(U)$ for every open set 
        $U\subset X$,
    \item[(2)] a functor $\rho_{_{VU}}:\,\me{C}(U)\rightarrow
        \me{C}(V)$ for any two open sets $V\subset U\subset X$, called
        restriction functor or just restriction, 
    \item[(3)] a natural equivalence $\Phi_{_{WVU}}:\,\rho_{_{WU}}
        \overset{\sim}{\ra}\rho_{_{WV}}\rho_{_{VU}}$ for every three
        open sets $W\subset V\subset U\subset X$.
  \end{itemize}
  This data should satisfy 
  \begin{itemize}
    \item[(P2F1)] $\rho_{_{UU}}=Id_{_{\me{C}(U)}}$ for every open set 
      $U\subset X$.
    \item[(P2F2)] $\Phi_{_{UUU}}=Id_{_{Id_{_{\me{C}(U)}}}}$ for every
      open set $U\subset X$.
    \item[(P2F3)]  For every four open sets $T\subset W\subset V\subset U
      \subset X$ the equation
      $$\Big(\rho_{_{TW}}\bullet\Phi_{_{WVU}}\Big)\circ\Phi_{_{TWU}}=
        \Big(\Phi_{_{TWV}}\bullet\rho_{_{VU}}\Big)\circ \Phi_{_{TVU}}, $$ 
      holds, i. e. the following diagram commutes:
      $$ \xymatrix@C=2cm@R=1.2cm{
            \rho_{_{TU}} \ar[r]^{\Phi_{_{TWU}}}_{\sim}
            \ar[d]_{\Phi_{_{TVU}}}^{\sim} & \rho_{_{TW}}\rho_{_{WU}} 
           \ar[d]^{\rho_{_{TW}}\bullet\Phi_{_{WVU}}}_{\sim} \\
            \rho_{_{TV}}\rho_{_{VU}} \ar[r]_{\Phi_{_{TWV}}\bullet
            \rho_{_{VU}}}^{\sim} &
            \rho_{_{TW}}\rho_{_{WV}}\rho_{_{VU}}  . }$$
  \end{itemize}
  We shall mostly denote a prestack $(\me{C},\rho,\Phi)$ by $\me{C}$
  for short. Moreover it is often convenient to denote the restriction
  functor $\rho_{VU}$ by $i_{VU}^{-1}$.
\end{defi}
Hence a prestack is just a contravariant 2-functor
$\me{C}:\,\mc{T}(X)\ra\mc{CAT}$ with strict identities. We therefore
immediately get the notion of a functor of prestacks (being a
2-natural transformation of the underlying 2-functors) and the notion
of a natural transformation of functors of prestacks (being a
modification of the underlying 2-natural transformations). In particular we 
get the concept of an equivalence of prestacks and we may 
define the (2-)category $\mc{PST}(X)$ of prestacks on $X$.
\begin{remark}
  \emph{Let $\me{C}$ be a prestack on $X$, $U\subset X$ an open subset and 
  $A,B\in\Ob{\me{C}(U)}$. For $V\subset U$, we set
 $$  \mc{H}\mathit{om}_{_{\me{C}|_{_U}}}(A,B)(V)=
     \Hom{\me{C}(V)}(A|_{_V},B|_{_V})$$
  If $W\subset V\subset U$ the restriction functor $\rho_{WV}$ and the natural 
  equivalence $\Phi_{WVU}$ define a restriction map and one easily verifies 
  that $\mc{H}\mathit{om}_{_{\me{C}|_{_U}}}(A,B)$ is a presheaf.}
\end{remark}
Let us add some notations
\begin{defi}
  \begin{itemize}
     \item[(1)] A prestack $\me{C}$ is called additive if for any $U\subset X$
        the category $\me{C}(U)$ is additive and the restriction functors
        are additive.
     \item[(2)] An additive prestack $\me{C}$ is called triangulated 
        if for any $U\subset X$ the category $\me{C}(U)$ is triangulated
        and the restriction functors are exact. 
     \item[(3)] An additive prestack $\me{C}$ is called abelian
        if for any $U\subset X$ the category $\me{C}(U)$ is abelian
        and the restriction functors are exact. 
  \end{itemize}
\end{defi}
We then get the obvious concept of an additive (resp. exact) functor between additive
(resp. triangulated or abelian) prestacks.

\subsection{Operations on prestacks}

Let $f:\,X\ra Y$ be a continous map.
\begin{prop}
  Let $\me{C}$ be a prestack on $X$. Then there is a natural prestack
  $f_*\me{C}$ on $Y$ such that for any open set $V\subset Y$ we have
  a canonical equivalence $f_*\me{C}(V)\simeq \me{C}(f^{-1}(V))$.
\end{prop}
\begin{prop}
  Let $\me{C}$ be a prestack on $Y$. Then there is a natural prestack
  $f^{-1}_{p}\me{C}$ on $X$ such that for any open set $U\subset X$ we have
  a canonical equivalence $f^{-1}_p\me{C}(U)\simeq \twodilim{f(U)\subset V}
  {\me{C}(V)}$.
\end{prop}
\begin{prop}
  The operations $f_*$ and $f^{-1}_p$ are (2-)adjoint to each other, i.e.
  there is a (2-)natural equivalence of categories
  $$ \ul{\mr{Hom}}_{_{\mc{PST}(X)}}
        (f^{-1}_p\me{C},\me{D}) \simeq \ul{\mr{Hom}}_{_{\mc{PST}(Y)}}
        (\me{C},f_*\me{D}). $$  
\end{prop}

\subsection{Stalks}

Since prestacks (and stacks) are often treated in the more general framework of
sites, we will describe here in detail the notion of a stalk of a prestack on a
topological space. Of course if $p\in X$ and $i:\,\{p\}
\hookrightarrow X$ is the inclusion, then the stalk $\me{C}_p$ of a prestack 
$\me{C}$ at $p$ is nothing but $i^{-1}_p\me{C}$. 

\begin{defi}
  Let $\me{C}$ be a prestack on $X$ and $p\in X$ a point.\\
  Consider the category $\mc{T}_p(X)$ of open sets that contain the
  point $p$. Note that since the set of open sets containig $p$ is
  stable by union and intersection the category $\mc{T}_p(X)$ is
  filtered and cofiltered.\\
  The prestack $\me{C}$ induces a 2-functor
  $\alpha_p:\,\mc{T}_p(X)^{\circ}\ra\mc{CAT}$.\\
  We set 
  $$ \me{C}_p=\twodilim{U\ni p}{\me{C}(U)}=
      \twodilim{U\in T_p(X)}{\alpha_p(U)} $$
  and call $\me{C}_p$ the \textbf{stalk} of $\me{C}$ at $p$ or the 
  \textbf{category of germs} of $\me{C}$ at $p$. 
\end{defi}
Hence the stalk of a prestack is defined up to canonical equivalence of
categories. It can easily be described using the explicit construction
of $2$-colimits. We get
$$ \Ob{\me{C}_p}=\Big\{(U,A)\ |\ \text{$p\in U\subset X$ open and 
    $A\in\Ob{\me{C}(U)}$}\Big\}=\bigsqcup_{p\in U\subset X}{\Ob{\me{C}(U)}}. $$
Let $(U,A)$,$(V,B)$ be two objects of $\me{C}_p$. Then
$$ \Hom{\me{C}_p}((U,A),(V,B))=
   \underset{p\in W\subset U\cap V}{\varinjlim}{\Hom{\me{C}(W)}(A|_{_W},
   B|_{_W})}. $$
Hence a morphism $f:\,(U,A)\ra (V,B)$ is defined on a small 
neighborhood $W\subset U\cap V$ of $p$. In particular we get
\begin{prop}
  Let $\me{C}$ be a prestack on $X$, $p\in X$ a point, $U\subset X$ an
  open set containing $p$ and $A,B\in\Ob{\me{C}(U)}$ two objects. Then
  we have a canonical isomorphism
  $$ \mc{H}\mathit{om}_{_{\me{C}|_{_U}}}(A,B)_p \overset{\sim}{\lra}
     \on{Hom}_{_{\me{C}_p}}(A,B). $$
  This isomorphism is compatible with the composition maps in $\me{C}$
  in the following way:\\
  Let $U\subset X$ be an open set containing $p$ and 
  $A,B,C\in\Ob{\me{C}(U)}$.\\
  Then the morphism of sheaves
  $$  \mc{H}\mathit{om}_{_{\me{C}|_{_U}}}(A,B)\times
      \mc{H}\mathit{om}_{_{\me{C}|_{_U}}}(B,C) \lra
      \mc{H}\mathit{om}_{_{\me{C}|_{_U}}}(A,C) $$
  induces in its stalks the composition in the stalk:
  $$  \on{Hom}_{_{\me{C}_p}}(A,B)\times
      \on{Hom}_{_{\me{C}_p}}(B,C) \lra 
      \on{Hom}_{_{\me{C}_p}}(A,C). $$
\end{prop}
Hence germs of morphisms may be seen as morphisms in the category of
germs. The compatibility with the composition map has an obvious corollary:
\begin{cor}\label{liftingisoingerms}
  Let $\me{C}$ be a prestack on $X$, $p\in X$ a point, $U,V\subset X$ two
  open sets containing $p$ and $A\in\Ob{\me{C}(U)}$,
  $B\in\Ob{\me{C}(V)}$ two objects. Then $A$ and $B$ are isomorphic in
  $\me{C}_p$ if and only if they are isomorphic on an open
  neighborhood of $p$.
\end{cor}
Let us observe that an object $A\in\Ob{\me{C}(U)}$ (with $p\in U$) is
isomorphic to all its restrictions to sets $V\in p$ but there is no 
equivalence relation imposed on the objects.\\
If $A\in\Ob{\me{C}(U)}$ we will still denote by $A$ its image in
$\mc{C}_p$.  If $f:\,A\ra B$ is a morphism in $\me{C}(U)$ then
we note $f_p:\,A\ra B$ its image in $\me{C}_p$. The reason why we do not write $A_p$
is given by the following remark.
\begin{remark}
  \emph{Consider a sheaf of rings $\me{A}$ and the stack $\mc{MOD}(\me{A})$.
  One shall beware that the natural functor
  $$ \mc{MOD}(\me{A})_p \lra \mc{MOD}(\me{A}_p) $$
  is not an equivalence of categories because the morphism
  $$ \mr{Hom}(\me{F},\me{G})_p \lra \mr{Hom}_{\me{A}_p}(\me{F}_p,\me{G}_p)$$
  is not an isomorphism in general.}
\end{remark}
\begin{prop}
  \begin{itemize}
     \item[(1)] If $\me{C}$ is additive then its stalks are additive 
       categories and the natural functors into the stalks are additive.
     \item[(2)] If $\me{C}$ is triangulated then its stalks are triangulated
       categories and the natural functors into the stalks are exact.
     \item[(3)] If $\me{C}$ is abelian then its stalks are abelian
       categories and the natural functors into the stalks are exact.
   \end{itemize}
\end{prop}

\subsection{Stacks}

\begin{defi}
  A prestack $\me{C}$ on $X$ is \textbf{separated} if for all open subsets
  $U\subset X$ and all objects $A,B\in\Ob{\me{C}(U)}$ the presheaf
  $$    \mc{H}\mathit{om}_{_{\me{C}|_{_U}}}(A,B) $$
  is a sheaf.
\end{defi}
\begin{defi}\label{patchcond}
  A prestack $\me{C}$ on $X$ is a \textbf{stack} if the following two
  conditions are satisfied
  \begin{itemize}
  \item[(i)] The prestack $\me{C}$ is separated
  \item[(ii)] Let $U=\bigcup_{i\in I}{U_i}$ be an open covering of an open
    subset $U\subset X$ and suppose that we are given the following data
    \begin{itemize}
    \item[(a)] for every $i\in I$ an object  $A_i\in\me{C}(U_i)$ 
    \item[(b)] for every $i,j\in I$ an isomorphism 
        $\sigma_{ij}:A_j|_{U_{ij}}\overset{\sim}{\ra}A_i|_{U_{ij}}$ 
        such that for any $i,j,k\in I$ the equation $\sigma_{jk}\circ
        \sigma_{ij}=\sigma_{ik}$ holds on $U_{ijk}$.
    \end{itemize}
    \end{itemize}
  Then there exists an object $A\in \me{C}(U)$ and isomorphisms
  $\rho_i:\,A|_{U_i}\ra A_i$ such that $\sigma_{ij}\circ \rho_j=\rho_i$.
\end{defi}
\begin{prop}
  Let $\me{C},\me{C}'$ be two stacks on $X$. Consider a functor
  $F:\,\me{C}\ra \me{C}'$. Then we have 
  \begin{itemize}
    \item[(1)] $F$ is faithful if and only if $F_p$ is faithful for 
        all $p\in X$.
    \item[(2)] $F$ is fully faithful if and only if $F_p$ is fully
      faithful for all $p\in X$. 
    \item[(3)] $F$ is an equivalence of stacks if and only if
      $F_p$ is an equivalence of categories for all $p\in X$.
  \end{itemize}
\end{prop}
\begin{proof}
We know that $F$ is faithful (resp. fully faithful) if and only if 
the morphisms of sheaves
$$ \mc{H}\mathit{om}_{_{\me{C}|_{_{U}}}}(A,B)\lra
   \mc{H}\mathit{om}_{_{\me{C}'|_{_{U}}}}(F(A),F(B)) $$
are monomorphisms (resp. isomorphisms). Since these two properties are
verified in the stalks, we immediately get (1) and (2).\\ 
Let's prove (3). Note that the condition is clearly necessary.\\
Now suppose that $F_p$ is an equivalence of categories for all $p\in X$.\\ 
By (2) we know that $F$ is fully faithful. Hence ist is
sufficent to show that $F$ is essentially surjective for any open
set $U\subset X$.\\
Let $A'\in\Ob{\me{C}'(U)}$. For any point $p\in U$ there is an object
$A_p\in \me{C}_p $ such that $F_p(A_p)\simeq A'_p$. Hence by corollary
\eqref{liftingisoingerms} there is an open neighborhood $V(p)$ of $p$, an 
object $A(p)\in\Ob{\big(\me{C}(V(p))\big)}$ and an isomorphism 
$\vartheta(p):\,F(A(p))\simeq A'|_{_{V(p)}}$.\\
These isomorphisms define a cocycle that patches together the objects
$F(A(p))$ to an object isomorphic to $A'$.
Since $F$ is fully faithful this cocycle can be lifted to a cocycle in
$\me{C}'$ where the $A(p)$ patch together to an object $A$ such that
$F(A)$ is isomorphic to $A'$. 
\end{proof}

\subsection{The stack associated to a prestack}

In this paragraph we will describe the stack associated to a prestack 
on a topological space. As is the case of the sheaf associated to a presheaf
this can be done explicitly and is less complicated than on an arbitrary 
site.

\begin{prop}
  Let $\me{C}$ be a prestack on $X$. Then there exists a separated
  prestack $\me{C}^\dagger$ on $X$ and a canonical functor
  $\eta_{\me{C}}^\dagger:\me{C}\ra\me{C}^\dagger$ that induces an 
  equivalence of categories on the stalks, such that any morphism 
  $\me{C}\ra\me{D}$ into a separated prestack $\me{D}$ factors uniquely 
  through $\me{C}^{\dagger}$ (up to unique equivalence).\\
  Moreover for any functor $F:\,\me{C}\ra\me{D}$ there exists a functor
  $F^\dagger:\,\me{C}^\dagger\ra\me{D}^\dagger$ such that the diagram
  commutes
  $$\xymatrix{
       {\me{C}} \ar[r]^F \ar[d] & {\me{D}} \ar[d] \\
       {\me{C}^\dagger} \ar[r]_{F^\dagger} & {\me{D}^\dagger}. }$$
\end{prop}
\em Proof. \em \\
Let $U\subset X$ be an open subset. Let us define a category
$\me{C}^\dagger(U)$. Set 
$$ \Ob{\me{C}^\dagger(U)}=\Ob{\me{C}(U)}$$
now let $A,B\in\Ob{\me{C}(U)}$ be two objects. Put
$$ \Hom{\me{C}^\dagger(U)}(A,B)=
  \Gamma(U,\mc{H}\mathit{om}_{_{\me{C}|_{_U}}}(A,B)^\dagger),$$
where $\mc{H}\mathit{om}_{_{\me{C}|_{_U}}}(A,B)^\dagger$ is the sheaf 
associated to $\mc{H}\mathit{om}_{_{\me{C}|_{_U}}}(A,B)$. Note that we have 
a canonical map $\Hom{\me{C}(U)}(A,B)\ra \Hom{\me{C}^\dagger(U)}(A,B)$
which is just the natural morphism from the presheaf $\Hom{\me{C}(U)}(A,B)$ 
into its associated sheaf.\\ 
The map 
$$ \mc{H}\mathit{om}_{_{\me{C}|_{_U}}}(A,B)\times
   \mc{H}\mathit{om}_{_{\me{C}|_{_U}}}(B,C) \lra
   \mc{H}\mathit{om}_{_{\me{C}|_{_U}}}(A,C)  $$
induces the composition in $\me{C}^\dagger$:
$$ \mc{H}\mathit{om}_{_{\me{C}^\dagger|_{_U}}}(A,B)\times
   \mc{H}\mathit{om}_{_{\me{C}^\dagger|_{_U}}}(B,C) \lra
   \mc{H}\mathit{om}_{_{\me{C}^\dagger|_{_U}}}(A,C).  $$
The restriction functors of $\me{C}^\dagger$ and the equivalences can easily 
be constructed by the universal property of the sheaf associated to a 
preshaef which also implies that all the axioms are verified.\\
The universal property of the separated prestack associated to a prestack also 
follows from the universal property of the sheaf associated to a presheaf.
$\hspace{0.6cm}_{\texttt{}\square}$
\begin{thm}
  Let $\me{C}$ be a prestack on $X$. Then there exists a stack 
  $\me{C}^\ddagger$ on $X$ together with a canonical functor of
  prestacks $\eta_{\me{C}}^\ddagger:\,\me{C}\ra\me{C}^\ddagger$ such that 
  any morphism $\me{C}\ra\me{D}$ into some stack $\me{D}$ factors uniquely
  through $\me{C}^{\ddagger}$.\\
  Moreover for any stack $\me{D}$ and any morphism of prestacks 
  $F:\,\me{C}\ra\me{D}$ there exists a canonical functor of stacks 
  $F^\ddagger$ such that the following diagram is commutative
  $$ \xymatrix{
        {\me{C}} \ar[r]^F \ar[d]_{\eta_{\me{C}}^\ddagger}  & 
        {\me{D}} \ar[d]^{\eta_{\me{C}}^\ddagger} \\
        {\me{C}^\ddagger} \ar[r]_{F^\ddagger} &  {\me{D}^\ddagger}. }$$
  Finally the functor $\eta_{\me{C}}^\ddagger:\,\me{C}\ra\me{C}^\ddagger$ 
  induces equivalences of categories on the stalks at 
  every point of $X$.
\end{thm}
\begin{proof}
By the proposition we may assume that $\me{C}$ is a separated prestack,
i.e. all associated presheaves are actually sheaves.\\ 
Let $U\subset X$ be an open subset. We have to define a category
$\mc{C}^\ddagger(U)$.\\
Consider families $A=\{(A_p,U_p^A)\}_{p\in U}$ where $U^A_p$ is an
open neighborhood of $p$ with $A_p\in \me{C}(U^A_p)$ and families of 
morphisms $\theta^A=\{\theta^A_{pq}\}_{p,q\in U}$ where 
$\theta^A_{pq}:\,A_q|_{_{U^A_{pq}}}\overset{\sim}{\ra}
A_p|_{_{U_{pq}^A}}$ is an isomorphism for all $p,q\in U$ (here 
$U^A_{pq}=U_p^A\cap U^A_q$) satisfying the cocycle condition. During the 
proof let us call a pair $(A,\theta^A)$ a cocycle on $U$.\\
We shall now define morphisms of cocycles. A morphism
$f:\,(A,\theta^A)\ra (B,\theta^B)$ consists of a family of germs
morphisms $f_p:\,(A_p,U^A_p)\ra (B_p,U^B_p)$ of $\on{Mor}{\me{C}_p}$ such
that for any point $p\in U$ there is an open set $U^f_p$ on which
$f_p$ is represented as a morphism $f_p:\,A_p|_{_{U^f_p}}\ra 
B_p|_{_{U^f_p}}$ (where $U^f_p\subset U_p^{AB}=U_p^A\cap U_p^B$ 
is an open neighborhood of $p$) satisfying the following 
compatibility condition: the diagram
 $$ \xymatrix@C=2cm@R=1.2cm{ 
   A_q|_{_{U^f_{pq}}} \ar[r]^{\theta^A_{pq}\|_{_{U^f_{pq}}}} 
  \ar[d]_{f_q\|_{_{U^f_{pq}}}} &   A_p|_{_{U^f_{pq}}}  
  \ar[d]^{f_p\|_{_{U^f_{pq}}}} \\
  B_q|_{_{U^f_{pq}}} \ar[r]_{\theta^B_{pq}\|_{_{ U^f_{pq}}}} & 
    B_p|_{_{U^f_{pq}}}  }$$
should be commutative for all $p,q\in U$.\\
Now define $\me{C}^{\ddagger}(U)$ to be the category of cocycles. The obvious
restriction maps define a prestack $\me{C}^\ddagger$ on $U$ (which is 
actually a presheaf). It is now tedious but straightforward, that 
$\me{C}^\ddagger$ is a stack that satisfies the universal property. Note 
that we need the assumption that $\me{C}$ is separated when proving the 
patching condition.
\end{proof}
\begin{cor}
  Let $\me{C}$ be a stack. Then there exists a stack $\me{C}'$, canonically 
  isomorphic to $\me{C}$, which is also a presheaf of categories. 
\end{cor}
\begin{cor}
  Let $\me{C}$ be a prestack and $F:\,\me{C}\ra\me{D}$ be a morphism into a
  stack $\me{D}$. Suppose that $F$ induces equivalences of categories in
  the stalks. Then $\me{D}$ is equivalent to the stack associated to
  $\me{C}$.
\end{cor}

\subsection{Patching of stacks}
\begin{thm}[Patching Theorem]\emph{ }\\
  Let $X=\bigcup_{i\in I}{U_i}$ be an open covering. Suppose that for
  any $i\in I$ we are given a stack $\me{C}_i$ and for any $i,j\in I$
  an equivalence of stacks
  $\theta_{ij}:\,\me{C}_j|_{_{U_{ij}}}\overset{\sim}{\ra}\me{C}_j|_{_{U_{ij}}}$
  satisfying the cocycle relation
  $ \theta_{ij}|_{_{U_{ijk}}} \circ \theta_{jk}|_{_{U_{ijk}}}=
   \theta_{ik}|_{_{U_{ijk}}}$.\\
  Then there exists a stack $\me{C}$ on $X$, unique up to equivalence
  of stacks and equivalences $\theta_i:\,\me{C}|_{_{U_i}}
  \overset{\sim}{\ra}\me{C}_i$ staisfying $\theta_{ij}\circ
  \theta_j|_{_{U_{ij}}}=\theta_i|_{_{U_{ij}}}$.
\end{thm}
    
\section{The functor of ind-micolocalization}

In this appendix, we first recall some definitions and statements of the theory
of analytic ind-sheaves from \cite{KS2}. Then we define Kashiwara's functor of 
ind-microlocalization and give (without proof) some basic properties that we used
in Sections 8 and 9. 

\subsection{Ind-sheaves}

If $\mc{C}$ is a category, one embeds $\mc{C}$ into the category of presheaves 
(of sets) on $\mc{C}$ by the fully faithful Yoneda-functor:
$$ \mc{C} \lra \widehat{\mc{C}}\qquad ;\qquad A\mapsto
   \on{Hom}_{\mc{C}}(\,\cdot\,,A) $$
where $\widehat{\mc{C}}$ is the category of contravariant functors 
$\mc{C}\ra\mc{S}et$. An object in the essential image of the Yoneda-functor is
called representable.\\ 
Note that $\widehat{\mc{C}}$ admits all small colimits since the category $\mc{S}et$
does but even if $\mc{C}$ admits colimits the Yoneda-functor does not commute with 
them.\\
One denotes by $\on{Ind}\mc{C}$ the full subcategory of $\widehat{\mc{C}}$ formed
by small filtered colimits of representable objects and calls it the category of 
ind-objects of $\mc{C}$. Then $\on{Ind}\mc{C}$ admits all small filtered colimits.\\
If $\mc{C}$ is abelian then $\on{Ind}\mc{C}$ is abelian and the Yoneda-functor
induces an exact fully faithful functor $\mc{C}\ra\on{Ind}\mc{C}$.\\
Now let $X$ be a locally compact topological space and fix a field $k$. One sets
$$ \I(k_X)=\on{Ind}\mc{M}od^c(k_X) $$
where $\mc{M}od^c(k_X)$ denotes the full subcategory of $\mc{M}od(k_X)$ formed by
sheaves with compact support. We call $\I(k_X)$ the category of ind-sheaves 
(of $k$-vector spaces). One can show that the prestack $X\supset U\mapsto \I(k_U)$
is a proper stack, in particular it is an abelian stack.\\
There are three important basic functors for ind-sheaves
\begin{eqnarray*}
  \iota&:&\ \mc{M}od(k_X) \lra \I(k_X) \quad; \quad \me{F}
        \mapsto \underset{U\subset\subset X}{\text{\textquotedblleft}\varinjlim
        \text{\textquotedblright}}{\ \me{F}_U}  \\
  \alpha&:&\ \I(k_X) \lra \mc{M}od(k_X) \quad ; \quad  
        \underset{i\in I}{\text{\textquotedblleft}\varinjlim\text{\textquotedblright}}
        {\ \me{F}_i} \mapsto \underset{i\in I}{\varinjlim}{\ \me{F}_i} \\
  \beta&:&\ \mc{M}od(k_X) \lra \I(k_X) \qquad \text{left adjoint to $\alpha$}  
\end{eqnarray*}
where we write $\text{\textquotedblleft}\varinjlim\text{\textquotedblright}$ for colimits
in the category $\I(k_X)$. All three functors induce functors of stacks.
\begin{prop}
  \begin{itemize}
     \item[(i)] The functor $\iota$ is fully faithful and exact.
     \item[(ii)] The functor $\alpha$ is exact.
     \item[(iii)] The functor $\beta$ is fully faithful and exact.
     \item[(iv)] The triple $(\beta,\alpha,\iota)$ is a triple of adjoint functors, 
        i.e. $\beta$ is left adjoint to $\alpha$ and $\alpha$ is left adjoint to
        $\iota$.
  \end{itemize}
\end{prop}
Note that since the functors $\iota,\alpha,\beta$ are exact they are well-defined in
the derived categories, guard the adjoint properties and $\alpha,\iota$ are still
fully faithful. An object $\me{F}\in\Db(k_X)$ is identified with $\iota\me{F}$ in
$\Db(\I(k_X))$.\\
There are internal operations on ind-sheaves
$$ (\,\cdot\,)\otimes(\,\cdot\,) \qquad \qquad \text{and} \qquad \qquad
   \me{IH}om(\,\cdot\,,\,\cdot\,) $$
and an external 
$$ \me{H}om(\,\cdot\,,\,\cdot\,):\ \I(k_X)\times \I(k_X) \lra \mc{M}od(k_X). $$ 
Moreover for any continous map $f:\,X\ra Y$ between locally compact spaces we get
the external operations
$$ f^{-1},f_*,f_{!!}, $$
where the notation $f_{!!}$ indicates that $\iota f_!\not\simeq f_{!!}\iota$.\\
While $\otimes$ and $f^{-1}$ are exact the other functors have a right derived functor
and pass to the derived category where we can define Poincar\'e-Verdier duality, i.e.
we have a right adjoint $f^!$ to $\Dr f_{!!}$ and we get the usual formalism of
Grothendiek's six operations. We will not recall here the various natural isomorphisms
relating these functors and refer to \cite{KS2} but let us summarize the commutation
properties with $\iota,\alpha,\beta$: 
\begin{prop}
  \begin{itemize}
     \item[(i)] The functor $\iota$ commutes to $\otimes, f^{-1},f^!,\Dr f_*$.
     \item[(ii)] The functor $\alpha$ commutes to $\otimes, f^{-1},\Dr f_*,\Dr f_{!!}$
        and $\alpha\Dr\me{IH}om(\,\cdot\,,\,\cdot\,)\simeq
        \Dr\me{H}om(\,\cdot\,,\,\cdot\,).$
     \item[(iii)] The functor $\beta$ commutes to $\otimes,f^{-1}$.
   \end{itemize}
\end{prop}
Finally let us state the following Proposition which has no counterpart in calssical
sheaf theory:
\begin{prop}
  Let $\me{F},\me{G}\in\Db(k_X)$ and $\me{M}\in\Db(\I(k_X))$. Then there is
  a natural isomorphism
  $$ \Dr\me{IH}om(\me{F},\me{M})\otimes\beta\me{G} \overset{\sim}{\lra}
     \Dr\me{IH}om(\me{F},\me{M}\otimes\beta\me{G}). $$ 
\end{prop}
 
\subsection{Microlocalization of ind-sheaves}
In \cite{K4}, Kashiwara establishes the following theorem
\begin{thm}
  There is a functor 
  $$ \mu:\ \Db(\I(k_X)) \lra \Db(\I(k_{T^*X})) $$
  such that for any $\me{F},\me{G}\in\Db(k_X)$ we have a natural isomorphism
  $$ \Dr\me{H}om(\mu\me{F},\mu\me{G})\simeq \Dr\me{H}om(\pi^{-1}\me{F},\mu\me{G}) 
     \simeq \mu hom(\me{F},\me{G}). $$
\end{thm}
\begin{remark}
  \em Note that if $\me{F}\in\Db(k_X)$, then
  $$ \supp(\mu\me{F})=\supp(\Dr\me{H}om(\mu\me{F},\mu\me{F}))=
     \supp (\mu hom(\me{F},\me{F}))=\SS(\me{F}). $$
\end{remark}
The construction of $\mu$ is rather straight-forward if we want to have the property
of the Theorem (cf. Proposition \ref{temperedsolutions} in Section 9.1). 
Let us recall the definition.\\
The normal deformation of the diagonal in $T^*X\times T^*X$ can be visualized 
by the following diagram
$$\xymatrix{{TT^*X}\ar[r]^(.3){\sim} & {T_{\Delta_{T^*X}}(T^*X\times T^*X)}
   \ar[d]^{\tau_{T^*X}} \ar@<-1pt>@{ (->}[r]^(.6){s} & 
    {\widetilde{T^*X\times T^*X}}\ar[d]_p   & {\Omega} 
   \ar@<1pt>@{ )->}[l]_(.3){j}\ar[ld]^(.4){\widetilde p} \\   
   & {T^*X}\ar@<-1pt>@{ (->}[r]_{\Delta_{T^*X}}  & 
   {T^*X\times T^*X} & {} }$$
Note that $\tilde p$ is smooth but $p$ is not. Also, the 
square is not cartesian. Set
$$
  \on{K}_X=\Dr p_{!!}\big(k_{\ol{\Omega}}\otimes\beta(k_P)\big)\otimes 
   \beta(\omega_{\Delta_{T^*X}|T^*X\times T^*X}^{\otimes -1})
$$
where the set $P\subset TT^*X$ is defined by
$$ P=\Big\{(x,\xi; v_x,v_\xi)\ |\ \langle v_x,\xi \rangle \geqs 0  \Big\}. $$
\begin{defi}
  Kashiwara's functor of microlocalisation is defined on $T^*X$ as
  $$  \mu:\ \on{D}^{\on{b}}(\I(k_{_{X}}))\lra 
           \on{D}^{\on{b}}(\I(k_{_{T^*X}}))\ ;\ \me{F}\mapsto
           \mu\me{F}=\on{K}_X\circ \pi^{-1}\me{F}. $$
\end{defi}
\begin{lemma}
  Let $S\subset T^*X$ be an arbitrary subset. Then $\mu$ defines functors
  $$ \mu:\ \Db(k_X,S) \lra  \Db(\I(k_S)).      $$
  If one considers these functors for open subsets $U\subset T^*X$, they 
  define functors of prestacks.
\end{lemma}
\begin{proof}
It is enough to show the existence of the first functor. If 
$\mu(\me{F})|_{_S}\simeq 0$, then $\supp(\mu(\me{F}))\cap S=\varnothing$, 
hence $\SS(\me{F})\cap S=\varnothing$ and $\mu(\,\cdot\,)|_{_S}$ factors 
through $\Db(k_X,S)$. 
\end{proof}
The following Proposition is used in Section 9.1.
\begin{prop}
  Let $\me{F}\in\Db_{\Rc}(k_X)$ and $\me{G}\in\Db(\I(k_X))$ and assume that
  $$ \SS(\me{F})\cap \supp(\mu\me{G})\subset T^*_XX$$
  Then there is a natural isomorphism
  $$ \Dr\me{H}om(\me{F},k_X)\otimes \me{G}
     \overset{\sim}{\lra} \Dr \me{IH}om(\me{F},\me{G}). $$
\end{prop}
The main theorem of \cite{K4} is the microlocal composition theorem that we included
in Section 8.3 (Theorem \ref{mainthm}).

%

\vspace*{1cm}

\noindent
Ingo Waschkies\\
Universit{\'e} Pierre et Marie Curie\\
Institut de Math{\'e}matiques\\
175, rue du Chevaleret, 75013 Paris, France\\
E-mail: ingo@math.jussieu.fr

\end{document}